\documentclass{article}
\begin{document}

\centerline{\bf\large $q-$Series Related with Higher Forms}
\vskip .4in

\centerline{\textbf{Nikos D. Bagis}}
\centerline{\textbf{Aristotele University of Thessaloniki-AUTH-Greece}}
\centerline{\textbf{nikosbagis@hotmail.gr}}

\begin{quote}

\[
\]

\centerline{\bf Abstract}
We examine $q-$series related to higher forms. These forms are cubics, quartics, etc. In some points, in the article we add parts from previous works, in such a way, the article be more complete and readable.

\[
\]
           				

\end{quote}

\section{Introduction.}

The study of quadratic forms have been built up by many great mathematicians such as Euler, Gauss, Dirichlet, Liouville, Eisenstein, Glaisher, Ramanujan among others. This theory has applications to a wide number of areas in modern mathematics including Gauss' circle problem in higher dimensions, class number theory, algebraic geometry, elliptic and theta functions, the Fermat-Wiles theorem, Eisenstein series and many other.\\ In this article using simple arguments we try to address the problem in higher dimensions such as cubes, quartics and as well to arbitrary degrees.\\  
\\  
We start with the well known quadratics. Let $K(x)$ be the complete elliptic integral of the first kind. This is given by
\begin{equation}
K(x)=\int^{\pi/2}_{0}\frac{d\theta}{\sqrt{1-x^2\sin^2(\theta)}}=\frac{\pi}{2}\cdot {}_2F_1\left(\frac{1}{2},\frac{1}{2};1;x^2\right),
\end{equation}
where ${}_2F_1$ is Gauss hypergeometric function.\\
In terms of Weber's $\lambda(\tau)-$modular function (see [3],[4]) 
\begin{equation}
\lambda(\tau)=16q\prod^{\infty}_{n=1}\left(\frac{1+q^{2n}}{1+q^{2n-1}}\right)^8,
\end{equation}
where $q=e^{i\pi \tau}$, $Im(\tau)>0$,
 $\tau=\sqrt{-r}$, the singular modulus  $k=k_r$, $r>0$ is  
\begin{equation}
k^2_{r}=\lambda(\tau)=\left(\frac{\theta_2(q)}{\theta_3(q)}\right)^4\textrm{, }
\end{equation}
with 
\begin{equation}
\theta_2(q)=\sum^{\infty}_{n=-\infty}q^{(n+1/2)^2}\textrm{ and }\theta_3(q)=\sum^{\infty}_{n=-\infty}q^{n^2}\textrm{, }|q|<1.
\end{equation} 
Also $k=k_r$, $0<k<1$ is the solution of the equation 
\begin{equation}
\frac{K(\sqrt{1-k_r^2})}{K(k_r)}=\sqrt{r}.
\end{equation}    
As usual we set $K=K(k_r)$ (the complete elliptic integral at singular values) and $K'=K(k'_r)$, where $k'_r=\sqrt{1-k_r^2}$ is the complementary singular modulus. The  Fourier expansion for the Jacobi elliptic function $\textrm{dn}$ (see [3] p.51-53) is\\
\begin{equation}
\textrm{dn}(q,u)=\frac{\pi}{2K}+\frac{2\pi}{K}\sum^{\infty}_{n=1}\frac{q^n}{1+q^{2n}}\cos(2nz),
\end{equation}  
where $z=\left(\frac{\pi}{2K}\right)u$ lies in the strip $|Im(z)|<\frac{\pi}{2}Im(\tau)$, $\tau=i\frac{K'}{K}$.\\

A very interesting connection between number theory and the theory of elliptic functions stems from  Jacobi's famous  theorem (see [3]):\\
\\
\textbf{Theorem 1.} (Jacobi)\\
If $q=e^{-\pi\sqrt{r}}$, $r>0$, then
\begin{equation}
\theta_3(q)=\sum^{\infty}_{n=-\infty}q^{n^2}=\sqrt{\frac{2K}{\pi}}.
\end{equation}

This theorem plays a key role in the theory of elliptic functions  and we shall use it here in our investigation of quadratic forms of general type.
It is very easy to see  by setting $u=0$ in (6), using $\textrm{dn}(q,0)=1$ and then multiplying both sides of (6) by $2K/\pi$, that  
$$
\frac{2K}{\pi}=1+4\sum^{\infty}_{m=1}\frac{q^m}{1+q^{2m}}=1+4\sum^{\infty}_{m=1}q^m\sum^{\infty}_{l=0}(-1)^lq^{2ml}=
$$
$$
=1+4\sum^{\infty}_{m=1}\sum^{\infty}_{l=0}(-1)^lq^{(2l+1)m}.
$$
Writing $n=(2l+1)m$, $d=2l+1$, so $l=(d-1)/2$, if $d$ runs through the odd divisors of $n$ we have 
\begin{equation}
\frac{2K}{\pi}=1+4\sum^{\infty}_{n=1}\left[\sum_{d-odd, d|n}(-1)^{\frac{d-1}{2}}\right]q^n.
\end{equation}
We define $\delta_0(n)=1$, if $n=0$ and $\delta_0(n)=4\sum_{d-odd, d|n}(-1)^{\frac{d-1}{2}}$, if $n\geq 1$.
If $r(n)$ denotes the number of representations of $n$ by the form  
$$
n=x^2+y^2\textrm{, }(x,y\in\bf Z\rm),
$$
then, if we consider the fact that
$$
\theta_3(q)^2=\sum^{\infty}_{n=-\infty}q^{n^2}\sum^{\infty}_{m=-\infty}q^{m^2}=\sum^{\infty}_{n,m=-\infty}q^{n^2+m^2}=\sum^{\infty}_{n=0}r(n)q^n
$$
and apply Jacobi's Theorem 1, we get \\
\\
\textbf{Theorem 2.} (Jacobi [8])\\
For $n=1,2,\ldots$ we have 
\begin{equation}
r(n)=4\sum_{d-odd,\textrm{ }d|n}(-1)^{\frac{d-1}{2}}
\end{equation}
and $r(0)=1$.\\ 

\section{Generalizations of Jacobi's two-square theorem}

Suppose we have two positive integers $A,B$, with $\gcd(A,B)=1$, and let $r_{A,B}(n)$ denote the number of representations of $n$ by the quadratic form
\begin{equation}
Ax^2+By^2.
\end{equation}
Then 
\begin{equation}
\theta_3\left(q^A\right)^2\theta_3\left(q^B\right)^2=\left(\sum^{\infty}_{n,m=-\infty}q^{An^2+Bm^2}\right)^2=\left(\sum^{\infty}_{n=0}r_{A,B}(n)q^n\right)^2.
\end{equation}
But, also
$$
\theta_3\left(q^A\right)^2\theta_3\left(q^B\right)^2=\left(\sum^{\infty}_{n=0}r(n)q^{nA}\right)\left(\sum^{\infty}_{m=0}r(m)q^{mB}\right)=
$$
\begin{equation}
=\sum^{\infty}_{n=0}\left(\sum_{kA+lB=n}r(k)r(l)\right)q^n.
\end{equation}
The linear Diophantine equation $kA+lB=n$ has solutions for all $n$ since $\gcd(A,B)=1|n$.\\

We now introduce the transformation $T$, which assigns the Taylor coefficient $f_n$ of a function $f(q)$ to the Taylor coefficient $\left(\sqrt{f}\right)_n$   of its square root $\sqrt{f(q)}$, i.e. 
\begin{equation}
\left(\sqrt{f}\right)_n=T(f_n).
\end{equation}
Hence
\begin{equation}
\sqrt{1+\sum^{\infty}_{n=1}a_nq^n}=1+\sum^{\infty}_{n=1}T(a_n)q^n.
\end{equation}

It is clear that if $T(a_n)=b_n$, then $b_n$ can evaluated from the recursion formula
\begin{equation}
a_n=\sum^{n}_{m=0}b_mb_{n-m}.
\end{equation} 
Also another way to evaluate the transform $T$ is using the Faa Di Bruno's Formula (see [1] p.823), which in this case is
\begin{equation}
T(f_n)=\sum_{m=0}^{n} h_m(f_0)\sum'\prod^{n}_{j=1}\frac{f_j^{a_j}}{a_j!},
\end{equation}
where the prime on the summation means that we sum over all non-negative integers $a_j$ such that $a_1+2a_2+3a_3+\ldots+na_n=n$ and $a_1+a_2+a_3+\ldots+a_n=m$. The function 
$h_m$ is $h_m(x)=(-1)^{m}x^{1/2-m}\left(\frac{-1}{2}\right)_m$, $(a)_m:=\frac{\Gamma\left(a+m\right)}{\Gamma(a)}$, where $m=1,2,\ldots$. But this evaluation is not so economical.\\

With the above notations we can proceed to\\ 
\\
\textbf{Theorem 3.}\\
Given two positive integers $A$,$B$ with $\textrm{gcd}(A,B)=1$, the number of the representations of $n=1,2,\ldots$ by the form $Ax^2+By^2$ is exactly 
\begin{equation}
r_{A,B}(n)=T\left(\sum_{kA+lB=n}r(k)r(l)\right).
\end{equation}
\\
Note that $r_{A,B}(0)$ is obviously 1.\\ 
\\
\textbf{Proposition 1.}\\
Given two positive integers $A$,$B$ with $\textrm{gcd}(A,B)=1$, the number of  representations of $n=1,2,\ldots$ by the form $Ax^2+By^2$ is exactly 
\begin{equation}
r_{A,B}(n)=\left[\frac{1}{n!}\frac{d^{n}}{dq^n}\sqrt{\sum^{n}_{t=0}\left(\sum_{kA+lB=t}r(k)r(l)\right)q^t}\right]_{q=0}.
\end{equation}

In the same way as above we can prove\\
\\
\textbf{Theorem 4.}\\
If $A_1,A_2,\ldots,A_N$ are positive integers such that $\textrm{gcd}(A_1,A_2,\ldots,A_N)=1$, the number of the representations of $n=1,2,\ldots$ by the form $\sum^{N}_{k=1}A_kx_k^2$ is exactly 
\begin{equation}
r_{2}(N,n)=T\left(\sum_{k_1A_1+k_2A_2+\ldots+k_NA_N=n}r(k_1)r(k_2)\ldots r(k_N)\right)
\end{equation}
and also
\begin{equation}
r_2(N,n)=\left[\frac{1}{n!}\frac{d^n}{dq^n}\sqrt{\sum^{n}_{t=0}\left(\sum_{k_1A_1+k_2A_2+\ldots+k_NA_N=t}r(k_1)r(k_2)\ldots r(k_N)\right)q^t}\right]_{q=0}.
\end{equation}
\\
\textbf{Proposition 2.}\\
Consider the non-homogeneous quadratic form  
\begin{equation}
Ax^2+By^2+Cx+Dy+E,
\end{equation}
with $A,B$ positive integers, $C,D,E$ general integers,   $gcd(A,B)=1$ and\\ $C\equiv0\textrm{ }\textrm{mod}(2A)$, $D\equiv0\textrm{ }\textrm{mod}(2B)$. Then $n$ has exactly 
\begin{equation}
r_{A,B}\left(n+\frac{C^2}{4A}+\frac{D^2}{4B}-E\right)
\end{equation}
representations in (21).\\
\\
\textbf{Proof.}\\
Write $C=-2L_1A$ and $D=-2L_2B$. Then $n=Ax^2+By^2+Cx+Dy+E$ is equivalent to $n=A(x-L_1)^2+B(y-L_2)^2-AL_1^2-BL_2^2+E$ and the number of representations of $n$ in (21) is equal to the number of representation of  $n+AL_1^2+BL_2^2-E=n+\frac{C^2}{4A}+\frac{D^2}{4B}-E$, by $Ax^2+By^2$. qed\\
\\
\textbf{Application 1.}\\
Let $A,B,C,D$ be as in Proposition 2, then
\begin{equation}
\sum^{\infty}_{n=-\infty}q^{An^2+Cn}\cdot\sum^{\infty}_{n=-\infty}q^{Bn^2+Dn}
=2\pi^{-1}q^{-n_0}K(k_r)\sqrt{m_{A,r}m_{B,r}},
\end{equation}
where $n_0=\frac{C^2}{4A}+\frac{D^2}{4B}$ and $q=e^{-\pi\sqrt{r}}$. The function $m_{n,r}=\frac{K(k_{n^2r})}{K(k_r)}$ is called a multiplier (see [4] pg.136) and takes algebraic values when, $n$ is a positive integer and $r$ is positive rational.\\
\\
\textbf{Proof.}\\
From Proposition 2 we have
$$
\sum^{\infty}_{n=-\infty}q^{An^2+Cn}\cdot\sum^{\infty}_{n=-\infty}q^{Bn^2+Dn}=\sum^{\infty}_{n,m=-\infty}q^{An^2+Bm^2+Cn+Dm}=
$$ 
$$
=\sum^{\infty}_{n=0}r_{A,B}(n)q^{n-n_0}
=q^{-n_0}\sum^{\infty}_{n=0}r_{A,B}(n)q^n
=q^{-n_0}\vartheta_3(q^A)\vartheta_3(q^B)=
$$
$$
=q^{-n_0}\sqrt{\frac{2K(k_{A^2r})2K(k_{B^2r})}{\pi^2}}
=q^{-n_0}\frac{2K}{\pi}\sqrt{m_{A,r}m_{B,r}}.
$$
\\
\textbf{Application 2.}\\
The equation
\begin{equation}
k(Ax^2+By^2+Cx+Dy+E)+l=n,
\end{equation}
have $r=r_{A,B}\left(\frac{n-l}{k}+\frac{C^2}{4A}+\frac{D^2}{4B}-E\right)$ solutions.\\
In general if $P_N(x)=\sum^{N}_{k=0}a_kx^k$ is a polynomial with integer coefficients
and if exists exactly one integer $n'$ such that $P_N(n')=n$, then
\begin{equation}
P_N\left(Ax^2+By^2+Cx+Dy+E\right)=n
\end{equation}
has 
\begin{equation} r_{A,B}\left(n'+\frac{C^2}{4A}+\frac{D^2}{4B}-E\right)
\end{equation}
integer solutions, (including 0).\\Furthermore, if the equation $P_N(n')=n$ has integer solutions $\{n'_1,n'_2,\ldots,n'_s\}$, then the number of representations of $n$ in (25) is
\begin{equation}
r=\sum_{i=1}^{s}r_{A,B}\left(n'_i+\frac{C^2}{4A}+\frac{D^2}{4B}-E\right).
\end{equation}
Any non-integer solution $n'$ of $P_N(n')=n$ leads to no representation (25) and hence makes no contribution to the sum (27). All the above hold of course, with the assumptions of Proposition 2.

\section{Representations in cubic forms}

In this section we give some formulas similar to Jacobi's Theorem 2 for the representations of a positive integer in cubic forms. The results of Section 2 can generalized to higher order terms under certain conditions. Historically, there are some results known regarding the cubic case. For example it is known that the Diophantine equation
\begin{equation}
ax^3-by^3=n,
\end{equation}      
for $a,b,n$ integers, have a finite number of solutions (see [6]).\\
From  the  Fermat-Wiles theorem it is known that
\begin{equation}
x^3+y^3=z^3,
\end{equation}      
have only trivial solutions i.e. $\{x,0,x\}$ and $\{0,x,x\}$.\\  
Also a result of Euler states that the equation
\begin{equation}
x^3+y^3=z^2,
\end{equation}
admits  a parametric solution in integers (see [5] p.578-579).\\
We proceed by stating and proving \\
\\  
\textbf{Theorem 5.}\\
The number representations of $n$ in the form $x^3+y^3$, ($x,y$ positive integers) is 
\begin{equation}
r^{+}_3(n)=\sum_{\scriptsize
\begin{array}{cc} 
	d|n\\
	d^3-4n=0
\end{array}\normalsize}1+2\cdot\sum_{\scriptsize
\begin{array}{cc} 
	d|n\\
	d^3-4n\neq0
\end{array}\normalsize
}X_{\textbf{N}}\left(\frac{d}{2}-\frac{1}{2}\sqrt{\frac{-d^2+4n/d}{3}}\right)
\end{equation} 
where $X_{\textbf{N}}(n)=1$ if $n$ is positive integer and 0 otherwise.\\ 
\\ 
\textbf{Proof.}\\
One has $x^3+y^3=(x+y)(x^2-xy+y^2)$. If we set $u=x+y$ and $v=x^2-xy+y^2$,  $x,y$ are given by 
$$
x=\frac{1}{6}\left(3u-\sqrt{-3u^2+12v}\right)\textrm{, }y=\frac{1}{6}\left(3u+\sqrt{-3u^2+12v}\right).
$$
Hence we get the necessary and sufficient conditions for $u,v$ to determine  an integer $n=uv$ that can be expressed as the sum of two cubes.\\ 
\\
\textbf{Remark.}
The quintic case is similar to the cubic. We have that the  number of representations of $n$ in the form $x^5+y^5$  ($x,y$ non-negative integers) is\\ 
$r_5(0)=1$ and if $n$ positive integer
\begin{equation}
r_5(n)=\sum_{\scriptsize
\begin{array}{cc} 
	d|n\\
	d^5-16n=0
\end{array}\normalsize}1
+2\cdot\sum_{\scriptsize
\begin{array}{cc} 
	d|n\\
	d^5-16n\neq0\\	
\end{array}\normalsize}X_{\textbf{N}}\left
(\frac{5d-\sqrt{-25d^2+10\sqrt{5d^4+20\frac{n}{d}}}}{10}\right), 
\end{equation}
where $X_{\textbf{N}}$ is the characteristic function on the non-negative integers.\\
\\

With the same arguments as in beginning of Section 2 we can use the same operator (14) and state the following\\
\\
\textbf{Theorem 6.}\\ 
Given two positive integers $A,B$ with gcd$(A,B)=1$, the number of the representations of $n\in\textbf{N}$ in the form $Ax^3+By^3$, $x,y>0$ is exactly 
\begin{equation}
s_{A,B}(n)=T\left(\sum_{kA+lB=n}r^{+}_3(k)r^{+}_3(l)\right),
\end{equation}  
where $r^{+}_3(n)$ is that of (31).\\
\\
\textbf{Proof.}\\
Assume the function $f(q)=\sum^{\infty}_{n=1}q^{n^3}$, then
$$
f\left(q^A\right)^2f\left(q^B\right)^2=\left(\sum^{\infty}_{n=1}s_{A,B}(n)q^n\right)^2.
$$ 
But also
$$
f\left(q^A\right)^2f\left(q^B\right)^2=\left(\sum^{\infty}_{n=1}r^{+}_3(n)q^{nA}\right)\left(\sum^{\infty}_{m=1}r^{+}_3(m)q^{mB}\right)=
$$ 
$$
=\sum^{\infty}_{n=1}\left(\sum_{kA+lB=n}r^{+}_3(k)r^{+}_3(l)\right)q^n.
$$ 
Since gcd$(A,B)=1|n$, we get the result. $qed$\\
\\
\textbf{Theorem 7.}\\
If $r_3(n)$ is the number of representations of the integer $n\geq 1$ in the form
\begin{equation}
x^3+y^3=n,
\end{equation}
where $x,y$ are any integers, then
\begin{equation}
r_3(n)=\sum_{\scriptsize
\begin{array}{cc}
	AB=n\\
	A^3=4n
\end{array}
\normalsize}1+2\sum_{\scriptsize
\begin{array}{cc}
	A B=n\\
	4B-A^2=3k^2>0\\
	A-k\equiv0(2)
\end{array}
\normalsize}1.
\end{equation}
Or using divisor sums is
\begin{equation}
r_3(n)=\sum_{\scriptsize 
\begin{array}{cc}
	0<d|n\\
	d^3=4n
\end{array}
\normalsize}1+2\sum_{\scriptsize 
\begin{array}{cc}
	0<d|n\\
	4n/d-d^2=3k^2>0\\
	d-k\equiv 0(2)
\end{array}
\normalsize}1.
\end{equation}
\\
\textbf{Proof.}\\
Assume now the representation
\begin{equation}
x^3+y^3=(x+y)(x^2-x y+y^2)=n
\end{equation}
If $x+y=u$, and $x^2-x y+y^2=v$, then $uv=n$ is the representation of $n$ in the form $uv$. Hence 
$$
x^2+y^2+2x y=u^2\textrm{ and }x^2+y^2-x y=v
$$ 
Hence if we write 
\begin{equation}
xy=s=\frac{1}{3}(u^2-v)\textrm{ and }x^2+y^2=m=\frac{1}{3}(u^2+2v),
\end{equation}
then
\begin{equation}
x+y=\pm\sqrt{m+2s}\textrm{ and }x-y=\pm\sqrt{m-2s}.
\end{equation}
The $x$ and $y$ are integers iff $m-2s=k^2$ and $m+2s=l^2$ and hence iff
\begin{equation}
2m=l^2+k^2\textrm{ and }4s=l^2-k^2,
\end{equation}
with $k,l$ both even or both odd. Then
\begin{equation}
x=\frac{l+k}{2}\textrm{ and }y=\frac{l-k}{2}.
\end{equation}
We denote
\begin{equation}
L(n,s,m):=\sum_{\scriptsize
\begin{array}{cc}
	-n\leq x,y\leq n\\
	xy=s\\
	x^2+y^2=m
\end{array}
\normalsize}1.
\end{equation} 
The arithmetic function $L(n,s,m)$, $n,s,m\in\textbf{N}$ takes only the values $0,2,4$ and
\begin{equation}
r_3(n)=\frac{1}{2}\sum_{d|n}L\left(n^2,\frac{d^2-n/d}{3},\frac{d^2+2n/d}{3}\right).
\end{equation}
Now if we assume that exist $l,k\in\textbf{Z}$ such that 
\begin{equation}
\frac{1}{3}(d^2-n/d)=\frac{1}{4}(l^2-k^2)\textrm{ and } \frac{1}{3}(d^2+2n/d)=\frac{1}{2}(l^2+k^2),
\end{equation}
then: If $k=0\Leftrightarrow m=2s\Leftrightarrow d^3=4n$, we have 2 solutions, and if relations (41) hold with $k\neq 0$, we have 4 solutions. Hence we can write:
$$ 
r_3(n)=\sum_{\scriptsize
\begin{array}{cc}
	d|n\\
	d^3=4n
\end{array}
\normalsize}1+2\sum_{\scriptsize
\begin{array}{cc}
	d|n, d^3\neq4n\\
	\frac{2}{3}(d^2+2n/d)=l^2+k^2\\
	\frac{4}{3}(d^2-n/d)=l^2-k^2\geq 0\\
	l-k\equiv0(2)
\end{array}
\normalsize}1.
$$
But since $\frac{2}{3}(d^2+2n/d)+\frac{4}{3}(d^2-n/d)=2d^2$, we have $l=d$ and we can write 
$$
r_3(n)=\sum_{\scriptsize
\begin{array}{cc}
	d|n\\
	d^3=4n
\end{array}
\normalsize}1+2\sum_{\scriptsize
\begin{array}{cc}
	d|n, d^2+3k^2=4n/d\\
	d-k\equiv0(2)
\end{array}
\normalsize}1
$$
Finally
$$
r_3(n)=\sum_{\scriptsize
\begin{array}{cc}
	AB=n\\
	A^3=4n
\end{array}
\normalsize}1+2\sum_{\scriptsize
\begin{array}{cc}
	A B=n\\
	4B-A^2=3k^2>0\\
	A-k\equiv0(2)
\end{array}
\normalsize}1.
$$
$qed$\\
\\

Next we will study third order theta functions of the form 
\begin{equation}
\phi(q):=\sum^{\infty}_{n=1}X(n)q^{n^3}\textrm{, }|q|<1.
\end{equation}
We will try to answer questions such: How these functions can be treated? Can these functions evaluated by simpler functions? Have these functions regularities such as the classical theta functions? Before proceeding to the theorems we give some notations.\\

In what will follow we shall use the function $X_{\textbf{N}}(n)$ of Theorem 5 and we will mean $r_3(n)=r^{+}_3(n)$. If we set
\begin{equation}
x^{\pm}_{d,n}:=\frac{d}{2}\pm\frac{\sqrt{\Delta_{d,n}}}{2},
\end{equation}
where $\Delta_{d,n}=p(d,n/d)$, $p(x,y):=\frac{4y-x^2}{3}$, 
then
$$
\sum^{\infty}_{n,m=1}X(n,m)q^{n^3+m^3}=\sum^{\infty}_{n=1}q^n\sum_{d|n\textrm{\scriptsize{,}}p(d,n/d)=0}X\left(\frac{d}{2},\frac{d}{2}\right)+
$$
\begin{equation}
+\sum^{\infty}_{n=1}q^n\sum^{*}_{d|n}X\left(x^{-}_{d,n},x^{+}_{d,n}\right)
+\sum^{\infty}_{n=1}q^n\sum^{*}_{d|n}X\left(x^{+}_{d,n},x^{-}_{d,n}\right).
\end{equation}
The asterisk $''*''$ in the sums means that only positive-integer values of $x^{\pm}_{d,n}$ give distribution in the sum, when $d$ runs through all positive divisors of $n$ i.e.
$$
\sum^{*}_{d|n}f(d,n/d):=\sum_{d|n\textrm{\scriptsize,\normalsize }p(d,n/d)\neq 0}X_{\textbf{N}}\left(x^{-}_{d,n}\right)X_{\textbf{N}}\left(x^{+}_{d,n}\right)f(d,n/d)=
$$
$$
=\sum_{\scriptsize
\begin{array}{cc}
d|n\\
p(d,n/d)\neq 0\\	
\end{array} \normalsize}X_{\textbf{N}}\left(x^{-}_{d,n}\right)
f(d,n/d)=
$$
\begin{equation}
=\sum_{\scriptsize
\begin{array}{cc}
d|n\\
p(d,n/d)=k^2\neq 0\\
2\leq d-k\equiv0(2)	
\end{array}\normalsize}f(d,n/d).
\end{equation}
\\
1) When the arithmetic function $X(n,m)$ is antisymmetric i.e. $X(n,m)=-X(m,n)$, then the second and third sums in the right side of (47) can be omitted.\\
2) When $X(n,m)$ is symmetric, then
$$
\sum^{\infty}_{n,m=1}X(n,m)q^{n^3+m^3}=\sum^{\infty}_{n=1}q^n\sum_{d|n\textrm{\scriptsize{,}}p(d,n/d)=0}X\left(\frac{d}{2},\frac{d}{2}\right)+
$$
\begin{equation}
+2\sum^{\infty}_{n=1}q^n\sum^{*}_{d|n}X\left(x^{-}_{d,n},x^{+}_{d,n}\right).
\end{equation}
Hence if we denote with
\begin{equation}
\textrm{Sym}X(n,m):=\frac{X(n,m)+X(m,n)}{2},
\end{equation}
the symmetric part of $X(n,m)$, then in general holds
$$
\sum^{\infty}_{n,m=1}X(n,m)q^{n^3+m^3}=\sum^{\infty}_{n=1}q^n\sum_{d|n\textrm{\scriptsize{,}}p(d,n/d)=0}X\left(\frac{d}{2},\frac{d}{2}\right)+
$$
$$
+2\sum^{\infty}_{n=1}q^n\sum^{*}_{d|n}\textrm{Sym}X\left(x^{-}_{d,n},y^{+}_{d,n}\right).
$$
Hence we get the next formula\\
\\
\textbf{Theorem 8.}\\
Assume that $X(n,m)$ is any double arithmetic function, then we have
$$
\sum^{\infty}_{n,m=1}X(n,m)q^{n^3+m^3}
=\sum^{\infty}_{n=1}q^n\sum_{\scriptsize
\begin{array}{cc}
	d|n\\
	p(d,n/d)=0
\end{array}
\normalsize}X\left(\frac{d}{2},\frac{d}{2}\right)
+
$$
$$
+2\sum^{\infty}_{n=1}q^n\sum_{d|n,p(d,n/d)\neq 0}X_{\textbf{N}}\left(x^{-}_{d,n}\right)\textrm{Sym}X\left(x^{-}_{d,n},x^{+}_{d,n}\right)=
$$
$$
=\sum^{\infty}_{n=1}q^n\sum_{\scriptsize
\begin{array}{cc}
	AB=n\\
	A^3=4n
\end{array}
\normalsize}X\left(\frac{A}{2},\frac{A}{2}\right)+
$$
\begin{equation}
+2\sum^{\infty}_{n=1}q^n\left(\sum_{\scriptsize
\begin{array}{cc}
	AB=n\\
	4B-A^2=3k^2\neq0\\
	2\leq A-k\equiv 0(2)
\end{array}\normalsize}
\textrm{Sym}X\left(\frac{A-k}{2},\frac{A+k}{2}\right)\right),
\end{equation}
where $\Delta_{d,n}=\frac{4}{3}n/d-\frac{1}{3}d^2=p(d,n/d)=k^2$.\\
\\

Now if $X(n,m)=n^2-(m+1)^2$, then $\textrm{Sym}X(n,m)=-(m+n+1)$ and we have
$$
\sum^{\infty}_{n,m=1}\left(n^2-(m+1)^2\right)q^{n^3+m^3}=
$$
$$
-\sum^{\infty}_{n=1}q^n\left(\sum_{d|n\textrm{\scriptsize{,}}4n/d-d^2=0}1\right)-\sum^{\infty}_{n=1}q^n\left(\sum_{d|n\textrm{\scriptsize{,}}4n/d-d^2=0}d\right)-
$$
$$
-2\sum^{\infty}_{n=1}q^n\sum^{*}_{d|n}\left(1+x^{-}_{d,n}+x^{+}_{d,n}\right)=
$$
$$
-\sum^{\infty}_{n=1}q^n\left(\sum_{d|n\textrm{\scriptsize{,}}4n/d-d^2=0}1\right)-\sum^{\infty}_{n=1}q^n\left(\sum_{d|n\textrm{\scriptsize{,}}4n/d-d^2=0}d\right)
-
$$
$$
-2\sum^{\infty}_{n=1}q^n\sum^{*}_{d|n}\left(1+d\right)=
$$
$$
-\sum^{\infty}_{n=1}q^n\left(\sum_{d|n\textrm{\scriptsize{,}}4n/d-d^2=0}1\right)-\sum^{\infty}_{n=1}q^n\left(\sum_{d|n\textrm{\scriptsize{,}}4n/d-d^2=0}d\right)-
$$
$$
-2\sum^{\infty}_{n=1}d^{*}(n)q^n-2\sum^{\infty}_{n=1}\sigma_1^{*}(n)q^n,
$$
where 
\begin{equation}
\sigma^{*}_{\nu}(n)=\sum^{*}_{d|n}d^{\nu}=\sum_{\scriptsize
\begin{array}{cc}
	d|n\\
	p(d,n/d)=k^2\neq0\\
	2\leq d-k\equiv0(2)
\end{array}
\normalsize}d^{\nu},
\end{equation} 
is the divisor sigma function of order 3 and $d^{*}(n):=\sigma_{0}^{*}(n)$. Moreover there holds
\begin{equation}
\sum^{\infty}_{n=1}q^n\left(\sum_{d|n\textrm{\scriptsize{,}}4n/d-d^2=0}1\right)=\sum^{\infty}_{n=1}q^{2n^3}.
\end{equation}
and 
\begin{equation}
\sum^{\infty}_{n=1}q^n\left(\sum_{d|n\textrm{\scriptsize{,}}4n/d-d^2=0}d\right)=\sum^{\infty}_{n=1}2nq^{2n^3}.
\end{equation}
Hence if we set
\begin{equation}
s_0(n):=\sum_{d|n\textrm{\scriptsize{,}}4n/d-d^2=0}1,
\end{equation}
and
\begin{equation}
s_1(n):=\sum_{d|n\textrm{\scriptsize{,}}4n/d-d^2=0}d,
\end{equation}
then
\begin{equation}
s_0(2n)=X_{3}(n)\textrm{ and }s_1(2n)=2\sqrt[3]{n}X_3(n),
\end{equation}
where $X_{\nu}(n)$ is the characteristic function on the positive powers of $\nu$ i.e. $X_{\nu}(n)=1$ if exists $m\in\textbf{N}:n=m^{\nu}$ and 0 else.\\ 
\\
\textbf{Theorem 9.}\\
If we denote
\begin{equation}
s_{\nu}(n):=\sum_{d|n\textrm{\scriptsize{,}}4n/d-d^2=0}d^{\nu}
\end{equation}
then
\begin{equation}
\sum^{\infty}_{n=1}2^{\nu}n^{\nu}q^{2n^3}=\sum^{\infty}_{n=1}s_{\nu}(n)q^n.
\end{equation}
\\

By these arguments also follows:
$$
d^{*}(n)=\sigma^{*}_0(n)=\sum^{*}_{d|n}1=\frac{1}{2}\left(r_3(n)-s_0(n)\right)=\frac{1}{2}\left(r_3(n)-X_3\left(\frac{n}{2}\right)\right)
$$
Hence for example if we write
$$
\sum^{\infty}_{n,m=1}(n^2-(m+1)^2)q^{n^3+m^3}=\sum^{\infty}_{n=1}A(n)q^{n},
$$
then
$$
A(n)=-s_0\left(n\right)-s_1(n)-2d^{*}(n)-2\sigma^{*}_1(n)=
$$
\begin{equation}
=-r_3(n)-s_1(n)-2\sigma^{*}_1(n).
\end{equation}
But also
$$
\sum^{\infty}_{n,m=1}\left(n^2-(m+1)^2\right)q^{n^3+m^3}=
$$
$$
=-\left(\sum^{\infty}_{n=1}q^{n^3}\right)^2-\left(\sum^{\infty}_{n=1}2nq^{n^3}\right)\left(\sum^{\infty}_{n=1}q^{n^3}\right)=
$$
\begin{equation}
=-\sum^{\infty}_{n=1}r_3\left(n\right)q^{n}+\sum^{\infty}_{n=1}\left(\sum^{2n}_{k=1}s_1(2n-k)s_0(k)\right)q^n.
\end{equation}
Hence from (60) and (61) we get
$$
-r_3(n)-s_1(n)-2\sigma_1^{*}(n)=-r_3(n)-\sum^{2n}_{k=1}s_1(2n-k)s_0(k),
$$
and finally\\
\\
\textbf{Theorem 10.}\\
If $n\in\textbf{N}$, then
\begin{equation}
\sigma^{*}_0(n)=\frac{1}{2}\left(r_3(n)-s_0(n)\right)
\end{equation}
and
\begin{equation}
\sigma_1^{*}(n)=-\frac{1}{2}s_1(n)+\frac{1}{2}\sum^{2n}_{k=1}s_1(2n-k)s_0(k).
\end{equation}
\\

In the case of
$$
\sum^{\infty}_{n,m=1}\left(n^3-(m+1)^3\right)q^{n^3+m^3}=\sum^{\infty}_{n=1}A(n)q^n,
$$
we have $X(n,m)=n^3-(m+1)^3$ and 
$$
\textrm{Sym}X(n,m)=-1-\frac{3}{2}(m^2+m)-\frac{3}{2}(n^2+n).
$$ 
Then 
$$
\sum^{\infty}_{n,m=1}\left(n^3-(m+1)^3\right)q^{n^3+m^3}=-\sum^{\infty}_{n=1}s_0(n)q^n-\frac{3}{2}\sum^{\infty}_{n=1}s_1(n)q^n-
$$
$$
-\frac{3}{4}\sum^{\infty}_{n=1}s_2(n)q^n-
2\sum^{\infty}_{n=1}\sigma^{*}_0(n)q^n-3\sum^{\infty}_{n=1}\sigma^{*}_1(n)q^n-
$$
\begin{equation}
-\sum^{\infty}_{n=1}\sigma^{*}_2(n)q^n-2n\sum^{\infty}_{n=1}\sigma^{*}_{-1}(n)q^n.
\end{equation}
But also
$$
\sum^{\infty}_{n,m=1}\left(n^3-(m+1)^3\right)q^{2n^3+2m^3}=\sum^{\infty}_{n,m=1}(n^3-m^3-3m^2-3m-1)q^{2n^3+2m^3}=
$$
$$
-\frac{3}{4}\left(\sum^{\infty}_{n=1}q^{2n^3}\right)\left(\sum^{\infty}_{n=1}4n^2q^{2n^3}\right)-\frac{3}{2}\left(\sum^{\infty}_{n=1}q^{2n^3}\right)\left(\sum^{\infty}_{n=1}2nq^{2n^3}\right)-
$$
$$
-\sum^{\infty}_{n=1}r_3(n)q^{2n}=
-\frac{3}{4}\sum^{\infty}_{n=1}\left(\sum^{2n}_{k=1}s_2(2n-k)s_0(k)\right)q^{2n}-
$$
\begin{equation}
-\frac{3}{2}\sum^{\infty}_{n=1}\left(\sum^{2n}_{k=1}s_1(2n-k)s_0(k)\right)q^{2n}-\sum^{\infty}_{n=1}r_3(n)q^{2n}.
\end{equation}
Hence equating both sides of (64) and (65), we get after simplifications
\begin{equation}
\sigma^{*}_2(n)+2n\sigma^{*}_{-1}(n)=-r_3(n)+s_0(n)-\frac{3}{4}s_2(n)-3\sigma^{*}_1(n).
\end{equation}
Assume now that $X(n,m)=nm$, then $\textrm{Sym}X(n,m)=nm$, $X\left(\frac{k}{2},\frac{k}{2}\right)=\frac{k^2}{4}$ and $\textrm{Sym}X\left(x^{-}_{d,n},x^{+}_{d,n}\right)=\frac{d^2}{3}-\frac{1}{3}\frac{n}{d}$. Therefore
\begin{equation}
\sum^{\infty}_{n,m=1}X(n,m)q^{2n^3+2m^3}=\frac{1}{4}\left(\sum^{\infty}_{n=1}2nq^{2n^3}\right)^2
=\frac{1}{4}\sum^{\infty}_{n=1}q^{2n}\sum^{2n}_{k=1}s_1(2n-k)s_1(k).
\end{equation}
Also
$$
\sum^{\infty}_{n,m=1}X(n,m)q^{2n^3+2m^3}=
$$
\begin{equation}
=\frac{1}{4}\sum^{\infty}_{n=1}q^{2n}\sum_{d|n,4n/d-d^2=0}d^2+\frac{2}{3}\sum^{\infty}_{n=1}\sigma^{*}_{2}(n)q^{2n}-\frac{2}{3}\sum^{\infty}_{n=1}n\sigma^{*}_{-1}(n)q^{2n}.
\end{equation}
Hence equating the second parts of (67) and (68), we get
\begin{equation}
\frac{1}{4}\sum^{2n}_{k=1}s_1(2n-k)s_1(k)=\frac{1}{4}s_2(n)+\frac{2}{3}\sigma^{*}_2(n)-\frac{2n}{3}\sigma^{*}_{-1}(n).
\end{equation}
Solving the system of (66),(69) we get the next\\
\\
\textbf{Theorem 11.}
$$
\sigma^{*}_2(n)=\frac{1}{4}\sum^{2n}_{k=1}s_{1}(2n-k)s_1(k)+\frac{1}{2}\sum^{2n}_{k=1}s_{1}(2n-k)s_0(k)+
$$
$$
+\frac{1}{4} \sum^{2n}_{k=1}s_{2}(2n-k)s_0(k)+\frac{1}{3}r_3(n)-\frac{1}{3}s_0(n)-\frac{1}{2}s_1(n)-
$$
\begin{equation}
-\frac{1}{2}s_2(n)-\frac{2}{3}\sigma^{*}_0(n)-\sigma^{*}_1(n).
\end{equation}
and
$$
\sigma^{*}_{-1}(n)=-\frac{1}{8n}\sum^{2n}_{k=1}s_1(2n-k)s_1(k)+\frac{1}{2n}\sum^{2n}_{k=1}s_1(2n-k)s_0(k)+
$$
$$
+\frac{1}{4n}\sum^{2n}_{k=1}s_2(2n-k)s_0(k)
+\frac{1}{3n}r_3(n)-\frac{1}{3n}s_0(n)-\frac{1}{2n}s_1(n)-
$$
\begin{equation}
-\frac{1}{8n}s_2(n)-\frac{2}{3n}\sigma^{*}_0(n)-\frac{1}{n}\sigma^{*}_1(n).
\end{equation}
\\

If $X(n,m)=(-1)^{n+m}$, then from Theorem 8 we get
$$
\sum^{\infty}_{n,m=1}(-1)^{n+m}q^{n^3+m^3}=\left(\sum^{\infty}_{n=1}(-1)^nq^{n^3}\right)^2=
$$
$$
=\sum^{\infty}_{n=1}q^n\left(\sum_{d|n\scriptsize{,}4n/d-d^2=0}(-1)^d\right)+2\sum^{\infty}_{n=1}q^n\sum^{*}_{d|n}(-1)^d=
$$
$$
=\sum^{\infty}_{n=1}q^nX_3\left(\frac{n}{2}\right)+2\sum^{\infty}_{n=1}q^n\left(\sum_{d|n,4n/d-d^2\neq0}X_{\textbf{N}}\left(x^{-}_{d,n}\right)(-1)^d\right)=
$$
$$
=\sum^{\infty}_{n=1}q^nX_3\left(\frac{n}{2}\right)+2\sum^{\infty}_{n=1}(-1)^nq^n\left(\sum_{d|n,4n/d-d^2\neq0}X_{\textbf{N}}\left(x^{-}_{d,n}\right)\right)=
$$
$$
=\sum^{\infty}_{n=1}d_3(n)q^n+ \sum^{\infty}_{n=1}(-1)^n\left(r_3(n)-d_3(n)\right)q^n=
$$
$$
=2\sum^{\infty}_{n=0}d_3(2n+1)q^{2n+1}+ \sum^{\infty}_{n=1}(-1)^nr_3(n)q^n=\sum^{\infty}_{n=1}(-1)^nr_3(n)q^n,
$$
since $d_3(2n+1)=0$ for all integers $n\geq 0$. Hence we get the next theorem.\\
\\
\textbf{Theorem 12.}\\
If $|q|<1$, then
\begin{equation}
\sum^{\infty}_{n,m=1}(-1)^{n+m}q^{n^3+m^3}
=\left(\sum^{\infty}_{n=1}(-1)^n q^{n^3}\right)^2
=\sum^{\infty}_{n=1}(-1)^nr_3(n)q^n
\end{equation}
and
\begin{equation}
\sum^{*}_{d|n}(-1)^d=(-1)^n\sum^{*}_{d|n}1=\frac{(-1)^n}{2}\left(r_3(n)-d_3(n)\right).
\end{equation}
In general
\begin{equation}
\sum^{*}_{d|n}(-1)^df(d)=(-1)^n\sum^{*}_{d|n}f(d),
\end{equation}
where 
\begin{equation}
d_{3}(n)=\sum^{*}_{d|n}1=\sum_{\scriptsize
\begin{array}{cc}
	d|n\\
	p(d,n/d)=k^2\neq 0\\
	2\leq d-k\equiv 0(2)
\end{array}
\normalsize}1.
\end{equation}
\\
\textbf{Theorem 13.}\\
If $f(x)$ is analytic in $\textbf{R}$, then
$$
f(x)=\sum^{\infty}_{n=0}\frac{f^{(n)}(0)}{n!}x^n
$$
and
\begin{equation}
\sum^{\infty}_{n=1}f(2n)q^{2n^3}=\sum^{\infty}_{n=1}s_f(n)q^n,
\end{equation}
where 
\begin{equation}
s_f(n)=\sum_{d|n\textrm{\scriptsize{,}}4n/d-d^2=0}f(d).
\end{equation}
\textbf{Proof.}\\
Easy.\\
\\
\textbf{Remark.} Since the above theorem exist for all $f(x)$ analytic in $\textbf{R}$ it will also hold for any arithmetic function $f(n):\textbf{N}\rightarrow\textbf{C}$.\\ 
\\
\textbf{Theorem 14.}\\
For all arithmetic functions $f:\textbf{N}\rightarrow\textbf{C}$, we have\\
1) If $p,n,m,l\in\textbf{N}$ and $(n,m)=1$, $p\neq 2l^3$, then 
\begin{equation}
s_f(pn)s_f(pm)=0.
\end{equation}
2) For all $n,m\in\textbf{N}$, we have 
\begin{equation}
s_f(2nm)=S_f(n,m),
\end{equation}
where 
\begin{equation}
S_f(n,m)=
\left\{\begin{array}{cc}
	f(2k)\textrm{, if }\exists k\in\textbf{N}:nm=k^3\\
	0\textrm{ else }
\end{array}\right\}.
\end{equation}
\\
\textbf{Theorem 15.}\\
Suppose that $f(n,w)$ is an arithmetic function of $n\in\textbf{N}$ and also a complex function of $w$, such that  
\begin{equation}
f(n,w)f(m,w)=f(n+m,w),
\end{equation} 
for all $n,m\in\textbf{N}:(n,m)=1$ and for all $w\in\textbf{C}$. Then if
\begin{equation}
\phi(w,q)=\sum^{\infty}_{n=1}q^{2n^3}f(2n,w),
\end{equation} 
we have\\
1)
\begin{equation}
\phi(w,q)^2=\sum^{\infty}_{l=1}q^{4l^3}f(4l,w)+2\sum^{\infty}_{l,t=1}q^{2l^3t}\sum_{\scriptsize\begin{array}{cc}
	d|t\\
	p(d,t/d)=k^2\neq 0\\
	2\leq d-k\equiv0(2)\\
\end{array}
\normalsize}f(2ld,w).
\end{equation}
2)
\begin{equation}
\left(\sum^{\infty}_{n=1}q^{n^3}e^{2\pi i n w}\right)^2
=\sum^{\infty}_{n=1}A(n,w)q^n,
\end{equation}
where
\begin{equation}
A(n,w)=\sum_{\scriptsize
\begin{array}{cc}
d|n\\
p(d,n/d)=0	
\end{array}
\normalsize}e^{2\pi i d w}
+2\sum_{\scriptsize\begin{array}{cc}
	d|n\\
	p(d,n/d)=k^2\neq 0\\
	2\leq d-k\equiv0(2)
\end{array}
\normalsize}e^{2\pi i d w}.
\end{equation}
\\
\textbf{Proof.}\\
1) It holds 
\begin{equation}
\phi(w,q):=\sum^{\infty}_{n=1}q^{2n^3}f(2n,w)=\sum^{\infty}_{n=1}q^{n}s_f(n,w).
\end{equation}  
Using Theorem 13 we have
$$
\phi(w,q)^2=\sum^{\infty}_{n,m=1}q^{n+m}s_f(n,w)s_f(m,w)
=
$$
$$
=\sum^{\infty}_{k=1}\sum_{\scriptsize\begin{array}{cc}
  n,m\geq1\\
  (n,m)=1
\end{array}
\normalsize}q^{(n+m)k}s_f(nk,w)s_f(mk,w)=
$$
$$
=\sum_{\scriptsize\begin{array}{cc}
	l,n,m\geq1\\
	(n,m)=1
\end{array}
\normalsize}q^{2l^3n+2l^3m}s_f(2l^3n,w)s_f(2l^3m,w)=
$$
$$
=\sum^{\infty}_{l=1}\sum_{\scriptsize\begin{array}{cc}
	n,m\geq1\\
	(n,m)=1
\end{array}
\normalsize}q^{2l^3(n+m)}s_f(2l^3n,w)s_f(2l^3m,w)=
$$
$$
=\sum^{\infty}_{l,t=1}q^{2l^3t}\sum_{\scriptsize\begin{array}{cc}
	n+m=t\\
	(n,m)=1
\end{array}
\normalsize}s_f(2l^3n,w)s_f(2l^3m,w)=
$$
$$
=\sum^{\infty}_{l,t=1}q^{2l^3t}\sum_{\scriptsize\begin{array}{cc}
	n^3+m^3=t\\
	(n,m)=1
\end{array}
\normalsize}f(2ln,w)f(2lm,w)=
$$
$$
=\sum^{\infty}_{l,t=1}q^{2l^3t}\sum_{\scriptsize\begin{array}{cc}
	n^3+m^3=t\\
	(n,m)=1
\end{array}
\normalsize}f(2l(n+m),w)=
$$
$$
=\sum^{\infty}_{l=1}q^{4l^3}f(4l,w)+2\sum^{\infty}_{l,t=1}q^{2l^3t}\sum_{\scriptsize\begin{array}{cc}
	d|t\\
	p(d,t/d)=k^2\neq 0\\
	2\leq d-k\equiv0(2)\\
	(d/2-k/2,d/2+k/2)=1
\end{array}
\normalsize}f(2ld,w)=
$$
$$
=\sum^{\infty}_{l=1}q^{4l^3}f(4l,w)+2\sum^{\infty}_{l,t=1}q^{2l^3t}\sum_{\scriptsize\begin{array}{cc}
	d|t\\
	p(d,t/d)=k^2\neq 0\\
	2\leq d-k\equiv0(2)\\
\end{array}
\normalsize}f(2ld,w).
$$
2) Applying Theorem 8 with $X(n,m)=f(n)f(m)$, $f(n)=e^{2\pi i n w}$, $X(n,m)=e^{2\pi i (n+m)w}$, $\textrm{Sym}X(n,m)=e^{2\pi i(n+m)w}$, $\textrm{Sym}X(x_{d,n}^{-},x_{d,n}^{+})=e^{2\pi i d w}$, we get
$$
\left(\sum^{\infty}_{n=1}q^{n^3}e^{2\pi i n w}\right)^2
=\sum^{\infty}_{n=1}q^n\sum_{d|n,p(d,n/d)=0}e^{2\pi i d w}
+
$$
$$
+2\sum^{\infty}_{n=1}q^n\sum_{\scriptsize\begin{array}{cc}
	d|n\textrm{, }p(d,n/d)\neq 0
\end{array}
\normalsize}X_{\textbf{N}}\left(x^{-}_{d,n}\right)e^{2\pi i d w}.
$$
$qed$\\
\\

Assume now the function
\begin{equation}
h(k,u,v):=\sum_{\scriptsize
\begin{array}{cc}
  1\leq n,m\leq k\\
	n+m=u\\
	n^2-nm+m^2=v
\end{array}\normalsize
}1.
\end{equation}
If $f(n,w)f(m,w)=f(n+m,w)$, for all $n,m\in\textbf{N}$ and
\begin{equation}
\phi(w,q)=\sum^{\infty}_{n=1}q^{n^3}f(n,w),
\end{equation}
then
$$
\phi(w,q)^2=\sum^{\infty}_{n,m=1}q^{n^3+m^3}f(n,w)f(m,w)=
$$
$$
=\sum^{\infty}_{n,m=1}q^{(n+m)(n^2-nm+m^3)}f(n,w)f(m,w)=
$$
$$
=\sum^{\infty}_{u,v=1}q^{uv}\sum_{\scriptsize
\begin{array}{cc}
  1\leq n,m\leq u v\\
	n+m=u\\
	n^2-nm+m^2=v
\end{array}\normalsize
}f(n,w)f(m,w)=
$$
$$
=\sum^{\infty}_{k=1}q^{k}\sum_{uv=k}\sum_{\scriptsize
\begin{array}{cc}
  1\leq n,m\leq u v\\
	n+m=u\\
	n^2-nm+m^2=v
\end{array}\normalsize
}f(n,w)f(m,w).
$$
Hence 
$$
\phi(w,q)^2=\sum^{\infty}_{k=1}q^{k}\sum_{uv=k}\sum_{\scriptsize
\begin{array}{cc}
  1\leq n,m\leq k\\
	n+m=u\\
	n^2-nm+m^2=v
\end{array}\normalsize
}f(n+m,w)=
$$
$$
=\sum^{\infty}_{k=1}q^{k}\sum_{uv=k}\sum_{\scriptsize
\begin{array}{cc}
  1\leq n,m\leq k\\
	n+m=u\\
	n^2-n m+m^2=v
\end{array}\normalsize
}f(u,w)=
$$
$$
=\sum^{\infty}_{k=1}q^{k}\sum_{d|k}f(d,w)\sum_{\scriptsize
\begin{array}{cc}
  1\leq n,m\leq k\\
	n+m=d\\
	n^2-n m+m^2=k/d
\end{array}\normalsize
}1=
$$
$$
=\sum^{\infty}_{k=1}q^{k}\sum_{d|k}f(d,w)h(k,d,k/d).
$$
Hence we get the next\\
\\
\textbf{Theorem 16.}\\
If $|q|<1$ and $w\in\textbf{C}$ and
\begin{equation}
\phi(w,q)=\sum^{\infty}_{n=1}q^{n^3}e^{2\pi i n w},
\end{equation}
then
\begin{equation}
\phi(w,q)^2=\sum^{\infty}_{n=1}q^{n}\sum_{d|n}e^{2\pi i d w}h(n,d,n/d),
\end{equation}
where 
\begin{equation}
h(k,u,v):=\sum_{\scriptsize
\begin{array}{cc}
  1\leq n,m\leq k\\
	n+m=u\\
	n^2-nm+m^2=v
\end{array}\normalsize
}1.
\end{equation}
\\
Also\\
\\
\textbf{Theorem 17.}\\
There holds
\begin{equation}
\sum_{d|n}h(n,d,n/d)=r_3(n).
\end{equation}
More general for any arithmetical function $f(n):\textbf{N}\rightarrow\textbf{C}$ we have
\begin{equation}
\sum_{d|n}f(d)h(n,d,n/d)=
\sum_{\scriptsize
\begin{array}{cc}
	AB=n\\
	A^3-4n=0
\end{array}
\normalsize}f(A)+2\sum_{\scriptsize 
\begin{array}{cc}
   AB=n\\
   4B-A^2=3k^2>0\\
	 2\leq A-k\equiv 0(2)
\end{array}
\normalsize}f(A)
\end{equation}
and if $|q|<1$, then
$$
\sum^{\infty}_{n,m=1}f(n+m)q^{n^3+m^3}=\sum^{\infty}_{n=1}q^n\sum_{d|n}f(d)h(n,d,n/d)=
$$
\begin{equation}
=\sum^{\infty}_{n=1}q^n
\sum_{\scriptsize\begin{array}{cc}
	AB=n\\
	A^3-4n=0
\end{array}
\normalsize}f(A)+2\sum^{\infty}_{n=1}q^n\sum_{\scriptsize 
\begin{array}{cc}
   AB=n\\
   4B-A^2=3k^2>0\\
	 2\leq A-k\equiv 0(2)
\end{array}
\normalsize}f(A).
\end{equation}
More general if $\chi_A(n,m)$ is a characteristic function on a subset $A$ of $\textbf{N}\times\textbf{N}$, then
\begin{equation}
\sum^{\infty}_{n,m=1}\chi_A(n,m)f(n+m)q^{n^3+m^3}=\sum^{\infty}_{n=1}q^n\sum_{d|n}f(d)h_{\chi_A}(n,d,n/d),
\end{equation}
where
\begin{equation}
h_{\chi_A}(k,u,v):=\sum_{\scriptsize
\begin{array}{cc}
  1\leq n,m\leq k\\
	n+m=u\\
	n^2-nm+m^2=v\\
	n,m\in A
\end{array}\normalsize
}1.
\end{equation}
\\
\textbf{Application 3.}\\
1) For any arithmetical function $f(n)$ we have 
\begin{equation}
\sum_{\scriptsize
\begin{array}{cc}
	n,m\geq 1\\
	(n,m)=1
\end{array}
\normalsize}f(n+m)q^{n^3+m^3}=\sum^{\infty}_{n=1}q^n\sum_{d|n}f(d)h^{*}(n,d,n/d), 
\end{equation}
where
\begin{equation}
h^{*}(k,u,v)=\sum_{\scriptsize
\begin{array}{cc}
  1\leq n,m\leq k\\
	n+m=u\\
	n^2-nm+m^2=v\\
	(n,m)=1
\end{array}\normalsize
}1.
\end{equation}
2) Assume that $f(n)$ is exponential arithmetic function in $\textbf{N}$ i.e. such that $f(0)=1$ and
\begin{equation}
f(n+m)=f(n)f(m)\textrm{, whenever }(n,m)=1\textrm{, }n,m\in\textbf{N}.
\end{equation}
Then
\begin{equation}
\sum_{\scriptsize
\begin{array}{cc}
	 n,m\geq 1\\
  (n,m)=1
\end{array}
\normalsize}f(n)f(m)q^{n^3+m^3}=\sum^{\infty}_{n=1}q^n\sum_{d|n}f(d) h^{*}(n,d,n/d),
\end{equation}
Actually $\sum_{d|n}h(n,d,n/d)=r_3(n)$ are the representations of $n$ as a sum of two natural cubes and $\sum_{d|n}h^{*}(n,d,n/d)$ is the representation of $n$ into two natural cubes prime each other. Of course the set of natural numbers is $\textbf{N}=\{1,2,3,\ldots\}$.\\
\\

Assume the theta function
\begin{equation}
\theta_{a,b}(\chi,q)=\sum^{\infty}_{n=1}\chi(n)q^{a n^2+b n}.
\end{equation}
We are going to evaluate its power series coefficients. For simplicity we assume $X(n,m)$ is such that $X(n,n)=\chi(n)$ and $\delta_{n,m}$ is 1, if $n=m$ and 0 else. Then 
$$
\theta_{a,b}(\chi,q)=\sum^{\infty}_{n,m=1}X(n,m)\delta_{n,m}q^{n(a m+b)}.
$$
Hence
\begin{equation}
\theta_{a,b}(\chi,q)=\sum^{\infty}_{n=1}q^n\sum_{d|n}X_{\textbf{\scriptsize{Z}}}\left(\frac{n/d-b}{a}\right)X\left(d,\frac{n/d-b}{a}\right)\delta_{d,(n/d-b)/a}.
\end{equation}
Consequently if $a,b$ are integers such $a>0$, then
\begin{equation}
\theta_{a,b}(\chi,q)=\sum^{\infty}_{n=1}q^n\sum_{\scriptsize 
\begin{array}{cc}
	d|n\\
	ad^2+bd=n
\end{array}
\normalsize}\chi(d).
\end{equation} 
In general\\
\\ 
\textbf{Theorem 18.}\\
Let $a_k$, $k=1,2,\ldots,\nu$ be integers and $a_{\nu}>0$ and the polynomial
$$
P(x)=\sum^{\nu}_{k=1}a_{k}x^k
$$ 
forms a non decreasing sequence in the set of positive integers. Then we have
\begin{equation}
\sum^{\infty}_{n=1}\chi(n)q^{P(n)}=\sum^{\infty}_{n=1}q^n\sum_{\scriptsize
\begin{array}{cc}
0<d|n\\
P(d)=n
\end{array}
\normalsize}\chi(d)\textrm{, }|q|<1.
\end{equation}
Hence\\ 
\\
\textbf{Theorem 19.}\\
Assume $P(x)$ is polynomial as in Theorem 18. The problem of finding, for a given positive integer $n$, the number of solutions of the equation 
\begin{equation}
P(x)=n,
\end{equation}
with $x$ positive integer, reduces to evaluate
\begin{equation}
\sum_{\scriptsize
\begin{array}{cc}
0<d|n\\
P(d)=n
\end{array}
\normalsize}1.
\end{equation}
If for the given integer $n>0$, the sum (106) is zero, then we have no positive integer solutions of $P(x)=n$.\\  
\\
\textbf{Theorem 20.}\\
Assume that $a_j\in\textbf{Z}$, $j=1,2,\ldots,\nu$ with $\nu$ even and $a_{\nu}>0$. Also assume that
\begin{equation}
P(x)=\sum^{\nu}_{k=1}a_kx^k
\end{equation}
is decreasing in $\textbf{Z}_{\leq c}$ and increasing in $\textbf{Z}_{>c}$, with $c\in\textbf{Z}$ and $P(c)=m$. We define $R_1(n)$ such that 
\begin{equation}
R_1(n)=\sum_{\scriptsize
\begin{array}{cc}
	0<d|n\\
	P(d)=n
\end{array}
\normalsize}\chi(d)+\sum_{\scriptsize
\begin{array}{cc}
	0>d|n\\
	P(d)=n
\end{array}
\normalsize}\chi(d)=\sum_{\scriptsize
\begin{array}{cc}
	d\neq 0\textrm{\scriptsize, }\textrm{\scriptsize abs}(d)|n\\
	P(d)=n
\end{array}
\normalsize}\chi(d).
\end{equation}
Then\\
\textbf{1)} If $m=P(c)=0$, have solutions $c=0\textrm{ and }c=c_1$, then\\
\textbf{i)} If $c_1$ is non zero integer, then
\begin{equation}
\sum^{\infty}_{n=-\infty}\chi(n)q^{P(n)}=\chi(0)+\chi(c_1)+\sum^{\infty}_{n=1}R_1(n)q^{n}.
\end{equation}
\textbf{ii)} If $c_1=0$,  or not an integer, then
\begin{equation}
\sum^{\infty}_{n=-\infty}\chi(n)q^{P(n)}=\chi(0)+\sum^{\infty}_{n=1}R_1(n)q^{n}.
\end{equation}
\textbf{2)}\\ 
\textbf{i)} If $m<0$ and $P(c)=m$  and $c\neq0$, $c\in\textbf{Z}$, then
\begin{equation}
\sum^{\infty}_{n=-\infty}\chi(n)q^{P(n)}=\chi(0)+\chi(c)+\sum^{l}_{j=1}\chi(r_j)q^{P(r_j)}+\sum^{\infty}_{n=1}R_1(n)q^{n},
\end{equation}  
where $r_j$, $j=1,\ldots,l$ are all integer points such that $P(r_j)<0$.\\
\textbf{ii)} If $m<0$ and $c$ not integer, then
\begin{equation}
\sum^{\infty}_{n=-\infty}\chi(n)q^{P(n)}=\chi(0)+\sum^{l}_{j=1}\chi(r_j)q^{P(r_j)}+\sum^{\infty}_{n=1}R_1(n)q^{n},
\end{equation}  
where $r_j$, $j=1,\ldots,l$ are all integer points such that $P(r_j)<0$.\\
\\
\textbf{Corollary 21.}\\
Assume that $P(x)=ax^2+bx$, with $a>0$. Assume also that $P(x)$ is decreasing in $\textbf{Z}_{\leq c}$ and increasing in $\textbf{Z}_{>c}$.\\
\textbf{1)} If $P(c)=0$ and\\
\textbf{i)} If $P(c)=0\Leftrightarrow c=0$, then
\begin{equation}
\sum^{\infty}_{n=-\infty}\chi(n)q^{an^2+bn}
=\chi(0)+\sum^{\infty}_{n=1}q^n\sum_{\scriptsize\begin{array}{cc}
d\neq 0\textrm{, }\textrm{abs}(d)|n\\
ad^2+bd=n
\end{array}\normalsize}\chi(d).
\end{equation}
\textbf{ii)} If $(P(c)=0\Leftrightarrow c=0,-1)$ or $(P(c)=0\Leftrightarrow c=0,1)$, then
 \begin{equation}
\sum^{\infty}_{n=-\infty}\chi(n)q^{an^2+bn}
=\chi(0)+\chi(\pm 1)+\sum^{\infty}_{n=1}q^n\sum_{\scriptsize\begin{array}{cc}
d\neq 0\textrm{, }\textrm{abs}(d)|n\\
ad^2+bd=n
\end{array}\normalsize}\chi(d).
\end{equation}
\textbf{2)} If $P(c)<0$, then\\ 
\textbf{i)} If $a|b$ and $P(k_j)$, $j=1,\ldots,l$ are all the negative values of $P(x)$ in $\textbf{Z}$, we have 
$$
\sum^{\infty}_{n=-\infty}\chi(n)q^{an^2+bn}=
$$
\begin{equation}
=\chi(0)+\chi(-b/a)+\sum^{l}_{j=1}\chi(k_j)q^{ak_j^2+bk_j}+\sum^{\infty}_{n=1}q^n\sum_{\scriptsize\begin{array}{cc}
d\neq 0\textrm{, }\textrm{abs}(d)|n\\
ad^2+bd=n
\end{array}\normalsize}\chi(d).
\end{equation}
\textbf{ii)} If $b/a$ is not an integer and $P(k_j)$, are all the negative values of $P(x)$ in $\textbf{Z}$, we have
$$
\sum^{\infty}_{n=-\infty}\chi(n)q^{an^2+bn}=
$$
\begin{equation}
=\chi(0)+\sum^{l}_{j=1}\chi(k_j)q^{ak_j^2+bk_j}+\sum^{\infty}_{n=1}q^n\sum_{\scriptsize\begin{array}{cc}
d\neq 0\textrm{, }\textrm{abs}(d)|n\\
ad^2+bd=n
\end{array}\normalsize}\chi(d).
\end{equation} 
\\
\textbf{Theorem 22.}\\
Assume $P_1(x),P_2(x)$ are polynomials as in Theorem 20. Then the number of representations of a positive integer $n$ in the form
\begin{equation}
P_1(x)+P_2(y)=n\textrm{, }x,y\in\textbf{Z}
\end{equation} 
is
\begin{equation}
R_{12}(n)=\sum^{n}_{l=0}R_1(l)R_2(n-l)\textrm{, }n=1,2,\ldots,
\end{equation}
where
\begin{equation}
R_{\{1,2\}}(n)=\sum_{\scriptsize
\begin{array}{cc}
	d\neq 0\textrm{\scriptsize, }\textrm{\scriptsize abs}(d)|n\\
	P_{\{1,2\}}(d)=n
\end{array}
\normalsize}1\textrm{, when }n=1,2,\ldots
\end{equation}
and $R_1(0)=r_{1}$, $R_2(0)=r_{2}$, where $r_{1}$, $r_{2}$ are the number of integer solutions of $P_1(n)=0$ and $P_2(n)=0$ respectively.\\
\\
\textbf{Example.}\\
i) Assume the form
\begin{equation}
ax^4+ay^4+b x^3+by^3=n\textrm{, }x,y\in\textbf{Z}
\end{equation}
with $a,b$ both integers and $a>0$ and $a-|b|\geq 0$. We set
\begin{equation}
R_1(n)=\sum_{\scriptsize
\begin{array}{cc}
	d\neq 0\textrm{\scriptsize, }\textrm{\scriptsize abs}(d)|n\\
	ad^4+bd^3=n
\end{array}
\normalsize}1,\textrm{ when }n=1,2,\ldots.
\end{equation}
If $R_1(0)$ are the number of integer roots of $an^4+bn^3=0$, then the representations of the non negative integer $n$ in (120) is 
\begin{equation}
R_2(n)=\sum^{n}_{l=0}R_1(n-l)R_1(l)\textrm{, }n=0,1,2,\ldots
\end{equation} 
ii) We now assume the form
\begin{equation}
a_1x^4+b_1x^3+a_2y^4+b_2y^3\textrm{, }x,y\in\textbf{Z}
\end{equation}
with $a_1,a_2,b_1,b_2$ integers, $a_1,a_2$ positive and $a_1-|b_1|,a_2-|b_2|$ non negative and set $P_1(x)=a_1x^4+b_1x^3$, $P_2(x)=a_2x^4+b_2x^3$. Then the number of representations of $n$ in (123) is
\begin{equation}
\sum^{n}_{l=0}R_1(l)R_2(n-l),
\end{equation}
where
\begin{equation}
R_{\{1,2\}}(n)=\sum_{\scriptsize
\begin{array}{cc}
	d\neq 0\textrm{\scriptsize, }\textrm{\scriptsize abs}(d)|n\\
	P_{\{1,2\}}(d)=n
\end{array}
\normalsize}1,\textrm{ when }n=1,2,\ldots,
\end{equation}
and $R_{\{1,2\}}(0)=r_{0\{1,2\}}$ is the number of solutions of $a_{\{1,2\}}n^4+b_{\{1,2\}}n^3=0$ respectively.\\
\\
\textbf{Remark.} We can assume also different type of polynomials $P_1$, $P_2$. For  example we can take $P_1(x)=a_1x^4+b_1x^3+c_1 x$ and $P_2(x)=a_2x^6+b_2x^3$. Then the representation theorem will holds if $P_1(x),P_2(x)$ are polynomials as in Theorem 20. Also  $R_1(0)=r_{01}$ are the number of distinct roots of $a_1n^4+b_1n^3+c_1n=0$ and $R_2(0)=r_{02}$ the number of distinct roots of $a_2n^6+b_2n^3=0$. The representation form is
\begin{equation}
a_1x^4+b_1x^3+c_1x+a_2y^6+b_2y^3\textrm{, }x,y\in\textbf{Z}.
\end{equation}

For more restricted forms we have evaluations like:\\
\\
\textbf{Theorem 23.}\\
For $n>0$ define
\begin{equation}
L_1(n):=\sum_{\scriptsize
\begin{array}{cc}
	d\neq 0\textrm{\scriptsize, }\textrm{\scriptsize abs}(d)|n\\
	ad^{\nu}=n
\end{array}
\normalsize}1
\textrm{ and }
L_2(n):=\sum_{\scriptsize
\begin{array}{cc}
	d\neq 0\textrm{\scriptsize, }\textrm{\scriptsize abs}(d)|n\\
	bd^{k}=n
\end{array}
\normalsize}1.
\end{equation}
For $n=0$ we set $L_1(0)=L_2(0)=1$. Then the number of representations of $n$ into
\begin{equation}
ax^{\nu}+by^{k}=n\textrm{, }x,y\textrm{ integers, }
\end{equation}
with $\nu,k$ positive even and $a,b>0$ is
\begin{equation}
R_2(n)=\sum^{n}_{l=0}L_1(l)L_2(n-l).
\end{equation}
\\
\textbf{Theorem 24.}\\
Assume the form
\begin{equation}
P_1(x)+P_2(y)\textrm{, }x,y\in\textbf{N}.
\end{equation}
When $P_{\{1,2\}}(x)$ are polynomials of degree $\nu_{\{1,2\}}>0$ not necessary even and $a_{\nu_{\{1,2\}}}>0$, $P_{\{1,2\}}(0)=0$, $P_{\{1,2\}}(x)$ non decreasing, we can find the number of representations of $n\in\textbf{N}$ in (130) when $x,y$ positive integers. Let
\begin{equation}
R_{\{1,2\}}(n)=\sum_{\scriptsize 
\begin{array}{cc}
	d>0\textrm{\scriptsize, }d|n\\
	P_{\{1,2\}}(d)=n
\end{array}
\normalsize}1.
\end{equation}
Then the representations are
\begin{equation}
R_{12}(n)=\sum^{n}_{l=1}R_{1}(l)R_{2}(n-l).
\end{equation}

\section{A fundamental formula for the complex integration}

\textbf{Proposition 3.}\\
If $x$ is positive real number and $f$ is analytic in $(-1,1)$ with $f(0)=0$, then
\begin{equation}
\exp\left(\int^{x}_{+\infty}f(e^{-t})dt\right)=\prod^{\infty}_{n=1}(1-e^{-nx})^{\frac{1}{n}\sum_{d|n}\frac{f^{(d)}(0)}{d!}\mu\left(\frac{n}{d}\right)},
\end{equation}
where $\mu$ is the Moebius-$\mu$ arithmetic function (see [2]) and take the values $(-1)^r$ when $n$ square free and product of $r$ primes, else is $0$. Also $\mu(1)=1$.\\
\\
\textbf{Proof.}\\
Because $f(0)=0$ and $f$ analytic in $(-1,1)$, the integral $\int^{x}_{+\infty}f(e^{-t})dt$ exists for every $x>0$. We assume that exists arithmetic function $X(n)$ such that: 
\begin{equation}
\exp\left(\int^{x}_{+\infty}f(e^{-t})dt\right)=\prod^{\infty}_{n=1}(1-e^{-nx})^{X(n)}.
\end{equation}
We will determinate this function $X$.\\
Taking logarithms in both sides of (134) we have
$$
\int^{x}_{+\infty}f(e^{-t})dt=\sum^{\infty}_{n=1}X(n)\log(1-e^{-nx})=-\sum^{\infty}_{n=1}X(n)\sum^{\infty}_{m=1}\frac{e^{-mx}}{m}=
$$
$$
=-\sum^{\infty}_{n,m=1}X(n)n\frac{e^{-mnx}}{mn}=-\sum^{\infty}_{n=1}\frac{e^{-nx}}{n}\sum_{d|n}X(d)d.\eqno{:(A)}$$
Derivating (A) we get
$$
f(x)=\sum^{\infty}_{n=1}e^{-nx}\sum_{d|n}X(d)d.\eqno{:(B)}
$$
But from analytic property of $f$ in $(-1,1)$ we have 
$$
f(x)=\sum^{\infty}_{n=1}\frac{f^{(n)}(0)}{n!}x^n
$$ 
and consequently
$$
f(e^{-x})=\sum^{\infty}_{n=1}\frac{f^{(n)}(0)}{n!}e^{-nx}.
$$
Therefore from (B) and the above relation it must be 
$$
\frac{f^{(n)}(0)}{n!}=\sum_{d|n}X(d)d.
$$
By applying the Moebius inversion theorem (see [2]) we get
$$
X(n)=\frac{1}{n}\sum_{d|n}\frac{f^{(d)}(0)}{d!}\mu\left(\frac{n}{d}\right).
$$
This completes the proof.\\ 
\\
\textbf{Corollary 25.}\\
Let $|q|<1$, then 
\begin{equation}
e^{-f(q)}=\prod^{\infty}_{n=1}\left(1-q^n\right)^{\frac{1}{n}\sum_{d|n}\frac{f^{(d)}(0)}{\Gamma(d)}\mu\left(\frac{n}{d}\right)}.
\end{equation}
\\
\textbf{Proof.}\\
Setting where $\frac{f^{(n)}(0)}{n!}=\frac{f_1^{(n)}(0)}{n!}n$ and using Proposition 3, we get immediately the result.\\
\\
\textbf{Theorem 26.}\\
If $f:\textbf{N}\rightarrow\textbf{N}$, is a strictly increasing arithmetical function, then
\begin{equation}
\exp\left(\sum^{\infty}_{n=1}\chi(n)q^{f(n)}\right)=\prod^{\infty}_{n=1}\left(1-q^n\right)^{-X(n)}\textrm{, }|q|<1,
\end{equation}
where
\begin{equation}
X(n)=\frac{1}{n}\sum_{\scriptsize
\begin{array}{cc}
1\leq k\leq n\\
f(k)|n
\end{array}
\normalsize}\chi(k)f(k)\mu\left(n/f(k)\right).
\end{equation}
Provided that both parts of (136) converge.\\
\\
\textbf{Proof.}\\
Assume an arbitrary arithmetical function $\chi(m)$ and the characteristic function $X_{f}(m)$, such that $X_f(m)=1$ if $m\in\textbf{N}$ is of the form $m=f(n)$ and otherwise $X_f(m)=0$. Then 
$$
\exp\left(\sum^{\infty}_{n=1}\chi\left(f(n)\right)q^{f(n)}\right)=\exp\left(\sum^{\infty}_{m=1}\chi(m)X_f(m)q^m\right)
=\prod^{\infty}_{n=1}(1-q^n)^{-Y(n)},
$$
where 
$$
Y(n)=\frac{1}{n}\sum_{d|n}\chi(d)X_f(d)d\mu\left(n/d\right)=\frac{1}{n}\sum_{f(d)|n}\chi(f(d))f(d)\mu\left(n/f(d)\right).
$$
Hence
\begin{equation}
\exp\left(\sum^{\infty}_{n=1}\chi\left(f(n)\right)q^{f(n)}\right)=\prod^{\infty}_{n=1}\left(1-q^n\right)^{-1/n\sum_{f(d)|n}\chi(f(d))f(d)\mu(n/f(d))}.
\end{equation}
Since $\chi(n)$ is arbitrary, we can replace $\chi(f(n))$ with $\chi(n)$. Hence the result follows, provided that both parts of (138) are convergent.\\
\\
\textbf{Notes.}\\
Let $P(x)$ be any polynomial with $\textrm{deg}(P)=\nu>0$ and leading coefficient $a_{\nu}$ positive integer. Let also $a_k$ are the remaining coefficients in $\textbf{Z}$. Then we can write according to Theorem 18
\begin{equation}
\exp\left(\sum^{\infty}_{n=1}\chi(n)q^{P(n)}\right)=\prod^{\infty}_{n=1}\left(1-q^n\right)^{-X(n)}\textrm{, }|q|<1,
\end{equation}
where
\begin{equation}
X(n)=\frac{1}{n}\sum_{d|n}\left(\sum_{\scriptsize
\begin{array}{cc}
0<\delta|d\\
P(\delta)=d
\end{array}
\normalsize}\chi(\delta)\right)d\mu\left(\frac{n}{d}\right).
\end{equation}
Also the Dirichlet series related with $X(n)$ is
\begin{equation}
L(X,s)=\sum^{\infty}_{n=1}\frac{X(n)}{n^s}=\frac{1}{\zeta(s+1)}\sum^{\infty}_{n=1}\frac{\chi(n)}{P(n)^s}.
\end{equation}
In case $P(n)=n^{\nu}$, then
\begin{equation}
L(X,s)=\frac{L(\chi,\nu s)}{\zeta(s+1)},
\end{equation}
where $L(\chi,s)=\sum^{\infty}_{n=1}\chi(n)n^{-s}$.\\
\\
\textbf{Theorem 27.}\\
The inversion formula for the divisor sum (131) read as:\\
If $f:\textbf{N}\rightarrow\textbf{N}$, $f$ increasing  and 
\begin{equation}
X(n)=\frac{1}{n}\sum_{\scriptsize
\begin{array}[pos]{cc}
	1\leq d\leq n\\
	f(d)|n
\end{array}
\normalsize}\chi(d)f(d)\mu\left(n/f(d)\right),
\end{equation}
then
\begin{equation}
\frac{1}{n}\sum_{d|n}X(d)d=\left\{
\begin{array}[pos]{cc}
	\chi\left(f^{(-1)}(n)\right)\textrm{, if }f^{(-1)}(n)\in\textbf{N}\\
	0\textrm{, otherwise }
\end{array}\right\}
\end{equation}
and the opposite. Also
\begin{equation}
X(n)=\frac{1}{n}\sum_{\scriptsize
\begin{array}[pos]{cc}
	0<d|n\\
	f^{(-1)}(d)\in\textbf{N}
\end{array}
\normalsize}\chi\left(f^{(-1)}(d)\right)d\mu(n/d).
\end{equation}
\\
\textbf{Theorem 28.}\\ 
Assume that $f(n):\textbf{N}\rightarrow\textbf{N}$, with $f$ strictly increasing and $\chi(n)$ any arithmetical function, then
\begin{equation}
\exp\left(\sum^{\infty}_{n=1}\chi(n)q^{f(n)}\right)=\prod^{\infty}_{n=1}\left(1-q^n\right)^{-X(n)},
\end{equation}
where
\begin{equation}
X(n)=\frac{1}{n}\sum_{\scriptsize 
\begin{array}[pos]{cc}
	0<d|n\\
	f^{(-1)}(d)\in\textbf{N}
\end{array}
\normalsize}\chi\left(f^{(-1)}(d)\right)d\mu(n/d).
\end{equation}
\\
More generally we have the next\\
\\
\textbf{Theorem 29.}\\
If $f(n):\textbf{N}\rightarrow\textbf{N}$ is any integer value arithmetical function strictly increasing and $f^{(-1)}(n)$ is its inverse, in the sense $f\left(f^{(-1)}(n)\right)=n$, then if also $|q|<1$:
$$
\exp\left(\sum^{\infty}_{n=1}\sum^{\infty}_{l=1}\frac{q^{lf(n)}}{l}\chi(n)\sum_{d|l}g(d)\right)=
$$
$$
=\exp\left(\sum^{\infty}_{n=1}\frac{q^{n}}{n}\sum_{\scriptsize
\begin{array}{cc}
	0<d|n\\
	f^{(-1)}(d)\in\textbf{N}
\end{array}
\normalsize}\chi\left(f^{(-1)}(d)\right)d\sum_{d_1|(n/d)}g(d_1)\right)=
$$
$$
=\exp\left(\sum^{\infty}_{n=1}\frac{q^{n}}{n}\sum_{\scriptsize
\begin{array}{cc}
1\leq d\leq n\\
	f(d)|n\\
\end{array}
\normalsize}\chi\left(d\right)f(d)\sum_{d_1|(n/f(d))}g(d_1)\right)=
$$
\begin{equation}
=\prod^{\infty}_{n=1}\left(1-q^n\right)^{-X(n)},
\end{equation}
where
$$
X(n)=\frac{1}{n}\sum_{\scriptsize
\begin{array}{cc}
	1\leq d\leq n\\
	f(d)|n
\end{array}
\normalsize}\chi(d)f(d)g\left(n/f(d)\right)=
$$
\begin{equation}
=\frac{1}{n}\sum_{\scriptsize
\begin{array}{cc}
	0<d|n\\
	f^{(-1)}(d)\in\textbf{N}
\end{array}
\normalsize}\chi\left( f^{(-1)}(d)\right)dg\left(n/d\right).
\end{equation}  
\\
\textbf{Proposition 4.}\\
\begin{equation}
\exp\left(\sum^{\infty}_{n=1}\chi(n)\frac{q^{f(n)}}{1-q^{f(n)}}\right)=\prod^{\infty}_{n=1}\left(1-q^n\right)^{-X(n)}\textrm{, }|q|<1,
\end{equation}
where
\begin{equation}
X(n)=\frac{1}{n}\sum_{\scriptsize
\begin{array}{cc}
	1\leq k\leq n\\
	f^{(-1)}\left(\textrm{gcd}(n,k)\right)\in\textbf{N}
\end{array}
\normalsize}\chi\left( f^{(-1)}\left(\textrm{gcd}\left(n,k\right)\right)\right)\textrm{gcd}(n,k).
\end{equation}
\\
\textbf{Proof.}\\
Set $g(n)=\phi(n)$ in Theorem 29, where $\phi(n)$ is the Euler-Totient function, then use the identities $\sum_{d|n}\phi(d)=n$ and 
\begin{equation}
\sum^{n}_{k=1}f(\textrm{gcd}(n,k))=\sum_{d|n}f(d)\phi(n/d).
\end{equation}
\\
\textbf{Example.}\\
If $\chi(n)=\left(\frac{n}{5}\right)$ and $f(n)=n^2+n$, then 
$$
f^{(-1)}(n)=\frac{1}{2}\left(-1+\sqrt{1+4n}\right)
$$
and
$$
\exp\left(\sum^{\infty}_{n=1}\chi(n)\frac{q^{n^2+n}}{1-q^{n^2+n}}\right)=\prod^{\infty}_{n=1}\left(1-q^n\right)^{-X(n)},
$$
where
$$
X(n)=\frac{1}{n}\sum^{n}_{k=1}
\textrm{gcd}(n,k)\left(\frac{1}{2}(-1+\sqrt{1+4\textrm{gcd}(n,k)})|5\right)
$$
and for $k$ integer
$$
(n|k):=\left\{
\begin{array}{cc}
	\left(\frac{n}{k}\right)\textrm{, if }n\textrm{ integer }\\
0\textrm{, if }n\textrm{ not integer }
\end{array}
\right\}.
$$
\\
\textbf{Application 4.}\\
i) If $g(n)=n$, then
\begin{equation}
\exp\left(\sum^{\infty}_{n=1}q^n\sum_{d|n}\frac{\sigma_1(d)}{d}\chi\left(f^{(-1)}\left(n/d\right)\right)\right)
=\prod^{\infty}_{n=1}\left(1-q^n\right)^{-X(n)},
\end{equation}
where
\begin{equation}
X(n)=\sum_{f(d)|n}\chi(d).
\end{equation}
In particular if $f(x)=x^2$, then
\begin{equation}
\exp\left(\sum^{\infty}_{n,l=1}\frac{q^{ln^2}}{l}\chi(n)\sigma_1(l)\right)=\prod^{\infty}_{n=	1}\left(1-q^n\right)^{-X(n)},
\end{equation}
where
\begin{equation}
X(n)=\sum_{d^2|n}\chi(d).
\end{equation}
In case of $\chi(n)=\mu(n)$, we have 
\begin{equation}
\sum_{d^2|n}\mu(d)=|\mu(n)|.
\end{equation}
This leads to
\begin{equation}
\exp\left(\sum^{\infty}_{n,l=1}\frac{q^{n^2l}}{l}\mu(n)\sigma_1(l)\right)=\prod^{\infty}_{n=1}\left(1-q^n\right)^{-|\mu(n)|}.
\end{equation}
Hence
\begin{equation}
\exp\left(\sum^{\infty}_{n,l=1}\frac{q^{n^2l}}{l}\mu(n)\sigma_1(l)\right)=\prod_{\scriptsize
\begin{array}{cc}
	a\geq 1\\
	a=square free
\end{array}
\normalsize}\frac{1}{1-q^a}\textrm{, }|q|<1.
\end{equation}
ii) Moreover if we set $f(n)=n$, $\sum_{d|n}g(d)=n^{\nu}$, $\chi(n)=-C(n^2)$ in Theorem 29 then we have
$$
\exp\left(-\sum^{\infty}_{n=1}\sum^{\infty}_{l=1}q^{nl}C(n^2)l^{\nu-1}\right)
=\exp\left(-\sum^{\infty}_{n=1}\frac{q^n}{n}\sum_{0<d|n}C(d^2)d\left(n/d\right)^{\nu}\right)=
$$
\begin{equation}
=\exp\left(-\sum^{\infty}_{n=1}q^n\sum_{d|n}d^{\nu-1}C(n^2/d^2)\right)=
\prod^{\infty}_{n=1}(1-q^n)^{X(n)},
\end{equation}
where (see definition (182) below)
\begin{equation}
X(n)=\frac{1}{n}\sum_{d|n}C(d^2)d\sum_{d_1|(n/d)}d_1^{\nu}\mu\left(\frac{n/d}{d_1}\right)=\sum_{d|n}C(n^2/d^2)\frac{h_{\nu}(d)}{d}.
\end{equation}
According to Borcherds if $f(z)=\sum^{\infty}_{n=0}C(n)q^n$ is a modular form in $\Gamma_0(4)$, with $C(n)$ zero at $n\equiv 2,3(\textrm{mod}4)$ and weight $\nu+1/2$, $\nu$ is even positive integer, then  
\begin{equation}
\Psi(z)=-C(0)\frac{B_{\nu}}{2\nu}+\sum^{\infty}_{n=1}\left(\sum_{0<d|n}C\left(n^2/d^2\right)d^{\nu-1}\right)e^{2\pi i n z},
\end{equation}
is a modular form of integer weight $2\nu$ when $\nu>0$. If $\nu=0$ and all $C(n)$ are integers then exists rational number $h$, such that (we omit the first term in (162))
\begin{equation}
q^{-h}\exp\left(\Psi(z)\right)=q^{-h}\prod^{\infty}_{n=1}\left(1-q^n\right)^{C(n^2)}
\end{equation}
is a modular form of weight $C(0)$ of $\textrm{SL}_2\left(\textbf{Z}\right)$, where $h$ is some rational.\\
iii) Assume that $p$ is positive integer and $a>0$, then
\begin{equation}
\vartheta(a,p;z)=q^{p/8+a^2/(2p)-a/2}\sum^{\infty}_{n=-\infty}(-1)^nq^{pn^2/2+(p/2-a)n}
\end{equation}
is a modular form of weight $1/2$ in $\Gamma(p)$. That is if $a_1,b_1,c_1,d_1$ are integers such that $a_1,d_1\equiv1(\textrm{mod}p)$, $b_1,c_1\equiv0(\textrm{mod}p)$ and $a_1d_1-b_1c_1=1$, we get 
\begin{equation}
\vartheta\left(\frac{a_1z+b_1}{c_1z+d_1}\right)=\epsilon\sqrt{c_1z+d_1}\vartheta(z),
\end{equation} 
where $\epsilon$ depends only on $a_1,b_1,c_1,d_1$ and  $\epsilon^{24}=1$. Set $a=16$, $p=128$, then
$$
\phi(z)=q^9\sum^{\infty}_{n=-\infty}(-1)^nq^{64n^2+48n}=\sum^{\infty}_{n=-\infty}(-1)^nq^{(8n+3)^2}=\sum^{\infty}_{n=-\infty}\chi_0(n)q^{n^2}=
$$
$$
=\sum_{n\equiv 3(\textrm{\scriptsize mod\normalsize})8}(-1)^{\frac{n-3}{8}}q^{n^2}=
$$
\begin{equation}
=q^9-q^{25}-q^{121}+q^{169}+q^{361}-q^{441}-q^{729}+q^{841}+\ldots
\end{equation}
is a modular form in $\Gamma_1(128)$ of weight $1/2$ and all coefficients are non zero only at  $1(\textrm{mod}4)$. Also $\chi_0(n)=(-1)^{\frac{n-3}{8}}$. More general
$$
\vartheta(a,p;8pz)=\sum^{\infty}_{n=-\infty}(-1)^nq^{(2np+p-2a)^2}=\sum_{\scriptsize \begin{array}{cc}
n\in\textbf{\scriptsize Z}\\
n\equiv p-2a(\textrm{\scriptsize mod\normalsize}2p)
\end{array}
\normalsize}(-1)^{\frac{n-(p-2a)}{2p}}q^{n^2}.
$$
Hence
\begin{equation}
\vartheta(a,p;z)=\sum_{\scriptsize \begin{array}{cc}
n\in\textbf{\scriptsize Z}\\
n\equiv p-2a(\textrm{\scriptsize mod\normalsize}2p)
\end{array}\normalsize
}(-1)^{\frac{n-(p-2a)}{2p}}q^{n^2/(8p)}.
\end{equation}
If we assume that $X(n,m)$ is any bouble arithmetical function set 
$$
\sum^{\infty}_{n,m=-\infty}X(n,m)q^{n^2+m^2}=\sum^{\infty}_{n=0}R(n)q^{n}
$$
and
\begin{equation}
\textrm{Sym}^{*}X\left(n,m\right):=\frac{1}{2}\left(X(n,m)+X(m,n)+X(-n,-m)+X(-m,-n)\right)
\end{equation}
and
\begin{equation}
A(k,n)=\left\{\begin{array}{cc}
\frac{1}{2}\textrm{Sym}^*X\left(-\frac{\sqrt{2n}}{2},\frac{\sqrt{2n}}{2}\right)\textrm{, if }\sqrt{2n}\in\textbf{N}\textrm{ and }k=\sqrt{2n}\\
\textrm{Sym}^*X\left(x^{-}_{k,n},x^{+}_{k,n}\right)\textrm{, if }\sqrt{2n}\textrm{ not in }\textbf{N}
\end{array}\right\},
\end{equation}
where 
\begin{equation}
x^{\pm}_{k,n}=\frac{1}{2}\left(\pm k-\sqrt{2n-k^2}\right).
\end{equation}
Then 
$$
R(n)=\sum_{\scriptsize\begin{array}{cc} 
0\leq |k|\leq\left[\sqrt{2n}\right]\\
2n-k^2=l^2\geq0
\end{array}\normalsize}A(k,n).
$$
Hence\\
\\
\textbf{Theorem 29.1}\\
In general holds
\begin{equation}
\sum^{\infty}_{n,m=-\infty}X(n,m)q^{n^2+m^2}=\sum^{\infty}_{n=0}q^n\sum_{\scriptsize\begin{array}{cc} 
0\leq |k|\leq\left[\sqrt{2n}\right]\\
2n-k^2=l^2\geq0
\end{array}\normalsize}A(k,n).
\end{equation}
\\

Setting in (171) 
$$
X(n,m)=\chi_{p-2a,2p}(n)(-1)^{\frac{n-(p-2a)}{2p}}\chi_{p-2a,2a}(m)(-1)^{\frac{m-(p-2a)}{2p}},
$$
where $\chi_{a,b}(n)=1$ if $n$ is integer of the form $n\equiv a(\textrm{mod}b)$, and 0 else, we get
$$
A(k,n)=A_0(k,n)=
$$
$$
-i^{\frac{4a-\sqrt{2n-k^2}}{p}}\chi_{p-2a,2p}\left(\frac{-k-\sqrt{2n-k^2}}{2}\right)\chi_{p-2a,2p}\left(\frac{k-\sqrt{2n-k^2}}{2}\right)-
$$
\begin{equation}
-i^{\frac{4a+\sqrt{2n-k^2}}{p}}\chi_{p-2a,2p}\left(\frac{-k+\sqrt{2n-k^2}}{2}\right)\chi_{p-2a,2p}\left(\frac{k+\sqrt{2n-k^2}}{2}\right).
\end{equation}
Hence we have 
$$
\left(\sum^{\infty}_{n=-\infty}(-1)^nq^{pn^2/2+(p/2-a)n}\right)^2=
$$
\begin{equation}
=q^{-a^2/p-p/4+a}\sum^{\infty}_{n=0}q^{n/(8p)}\sum_{\scriptsize\begin{array}{cc} 
0\leq |k|\leq\left[\sqrt{2n}\right]\\
2n-k^2=l^2\geq0
\end{array}\normalsize}A_0(k,n),
\end{equation}
where 
$$
A_0(k,n)=-i^{\frac{4a-l}{p}}\chi_{p-2a,2p}\left(-\frac{k+l}{2}\right)\cdot\chi_{p-2a,2p}\left(\frac{k-l}{2}\right)-
$$
\begin{equation}
-i^{\frac{4a+l}{p}}\chi_{p-2a,2p}\left(-\frac{k-l}{2}\right)\cdot\chi_{p-2a,2p}\left(\frac{k+l}{2}\right)
\end{equation}
and $l=\sqrt{2n-k^2}$. Hence we get the next:\\
\\
\textbf{Theorem 29.2}\\
If $a,p$ are integers with $p>0$ and $p>2|a|$, then
\begin{equation}
\left(\sum^{\infty}_{n=-\infty}(-1)^nq^{pn^2/2+(p/2-a)n}\right)^2
=q^{-a^2/p-p/4+a}\sum^{\infty}_{n=0}C(a,p,n)q^{n/(8p)},
\end{equation}
where 
\begin{equation}
C(a,p,n)=-\sum_{\scriptsize\begin{array}{cc} 
0\leq |k|\leq\left[\sqrt{2n}\right]\\
2n-k^2=l^2\geq1\\
k\equiv 0(\textrm{mod}2p)\\
l-k\equiv(2p\pm 4a)(\textrm{mod}4p)
\end{array}\normalsize}(-1)^{\frac{4a\mp l}{2p}}.
\end{equation}
\\
\textbf{Theorem 29.3}\\
If $\chi(n)$ is full multiplicative function in $\textbf{Z}$ with $\chi(0)=0$, $\chi(1)=1$, then
\begin{equation}
\left(\sum^{\infty}_{n=-\infty}\chi(n)q^{n^2}\right)^2=\sum^{\infty}_{n=0}C_{\chi}(n)q^{n/16},
\end{equation}
where 
\begin{equation}
C_{\chi}(n)=\sum_{\scriptsize\begin{array}{cc} 
0\leq |m|\leq\left[\sqrt{2n}\right]\\
2n-m^2=0\\
m\equiv 0(\textrm{mod}8)
\end{array}\normalsize}\chi\left(\frac{-m^2}{64}\right)+2\sum_{\scriptsize\begin{array}{cc} 
0\leq |m|\leq\left[\sqrt{2n}\right]\\
2n-m^2=l^2\geq1\\
m\equiv 0(\textrm{mod}4)\\
l-m\equiv 0(\textrm{mod}8)
\end{array}\normalsize}\chi\left(\frac{l^2-m^2}{64}\right).
\end{equation}
\\
iv) Also one can easily get
\begin{equation}
\exp\left(\sum^{\infty}_{n=1}\frac{q^n}{n^{2-\nu}}\sum_{d|n}\chi(n/d)(n/d)^{2-\nu}\right)=\prod^{\infty}_{n=1}\left(1-q^n\right)^{-X(n)},
\end{equation}
where
\begin{equation}
X(n)=\frac{1}{n}\sum_{d|n}\chi(d)d\sum_{d_1|(n/d)}d_1^{\nu-1}\mu\left(\frac{n/d}{d_1}\right).
\end{equation}
Hence
$$
\exp\left(\sum^{\infty}_{n=1}\frac{q^{n}}{n^{\nu}}\sum_{d|n}\chi(d)d^{\nu-1}\right)=
$$
\begin{equation}
=\prod^{\infty}_{n=1}\left(1-q^n\right)^{-1/n\sum_{d|n}\chi(d)\sum_{d_1|(n/d)}d_1^{1-\nu}\mu\left(\frac{n/d}{d_1}\right)}.
\end{equation}
Setting the arithmetical function $h_{a}(n)$ to be
\begin{equation}
h_{a}(n)=\sum_{d|n}d^{a}\mu\left(\frac{n}{d}\right),
\end{equation}
we get\\
\\
\textbf{Theorem 30.}\\
If $|q|<1$, we have
\begin{equation}
\exp\left(\sum^{\infty}_{n=1}\frac{q^{n}}{n^{\nu}}\sum_{d|n}\chi(d)d^{\nu-1}\right)=\prod^{\infty}_{n=1}\left(1-q^n\right)^{-1/n\sum_{d|n}\chi(d)h_{1-\nu}(n/d)}.
\end{equation}
A more general result which can be derived from the above formula is
\begin{equation}
\exp\left(\sum^{\infty}_{n=1}q^n\sum_{d|n}d^{-1}f(d,n/d)\right)=\prod^{\infty}_{n=1}\left(1-q^n\right)^{-X(n)},
\end{equation}
where
\begin{equation}
X(n)=\frac{1}{n}\sum_{d|n}\sum_{\delta|(n/d)}\delta f(d,\delta)\mu\left(\frac{n/d}{\delta}\right).
\end{equation} 
Also
\begin{equation}
\exp\left(\sum^{\infty}_{n=1}q^n\sum_{d|n}f(n,d)\right)=\prod^{\infty}_{n=1}\left(1-q^n\right)^{-X(n)},
\end{equation}
where
\begin{equation}
X(n)=\frac{1}{n}\sum_{d|n}\sum_{\delta|(n/d)}d\delta f(d\delta,d)\mu\left(\frac{n/d}{\delta}\right).
\end{equation}
\\

But it is well known that\\
\\
\textbf{Theorem 31.}\\
If $\chi(n)$ is Dirichlet character modulo $N$ and $k$ a positive integer such $\chi(-1)=(-1)^k$, then the function
\begin{equation}
G_{k,\chi}(z)=c_k(\chi)+\sum^{\infty}_{n=1}q^n\sum_{d|n}\chi(d)d^{k-1}\textrm{, }q=e(z)\textrm{, }Im(z)>0
\end{equation}
is a modular form of weight $k$ in $M_k\left(\Gamma_0(N),\chi\right)$. Here $c_k(\chi)=\frac{1}{2}L(1-k,\chi)$ and $L(s)$ is the analytic continuation of  $L(s)=\sum^{\infty}_{n=1}\chi(n)n^{-s}$.\\
\\

Assuming the above two theorems we can get applications.

\[
\]

If $E_{\nu}(z)$ denotes the Eisenstein series, then
\begin{equation}
E_{-\nu}(z)=1+\frac{2}{\zeta(\nu+1)}\sum^{\infty}_{l=1}\frac{\sigma_{\nu+1}(l)}{l^{\nu+1}}q^l
\end{equation}
and
\begin{equation}
\frac{1}{2\pi i}\sum^{\infty}_{l=1}\frac{\sigma_{\nu-1}(l)}{l^{\nu}}q^l=\frac{\zeta(\nu-1)}{2}\int^{z}_{i\infty}\left(E_{2-\nu}(w)-1\right)dw\textrm{, }\nu\neq 1
\end{equation}
Hence
$$
\exp\left(\frac{1}{4\pi i}\sum^{\infty}_{n,l=1}\mu(n)\frac{q^{lf(n)}}{n^s}\frac{\sigma_{\nu-1}(l)}{l^{\nu}}\right)=
$$
$$
=\exp\left(\frac{1}{2}\zeta\left(\nu-1\right)\sum^{\infty}_{n=1}\frac{\mu(n)}{n^s}\int^{zf(n)}_{i\infty}\left[E_{2-\nu}\left(w\right)-1\right]dw\right)=
$$
$$
=\exp\left(\frac{1}{2}\zeta\left(\nu-1\right)\sum^{\infty}_{n=1}\frac{\mu(n)f(n)}{n^s}\int^{z}_{i\infty}\left[E_{2-\nu}\left(wf(n)\right)-1\right]dw\right).
$$
Now we use Theorem ? with $\chi(m)=\mu(m) m^{-s}$ and $g(m)=m^{-\nu+1}$, to get
$$
\exp\left(\frac{1}{2}\zeta\left(\nu-1\right)\sum^{\infty}_{n=1}\frac{\mu(n)}{n^s}\int^{z f(n)}_{i\infty}\left[E_{2-\nu}(w)-1\right]dw\right)
=
$$
\begin{equation}
=\prod^{\infty}_{n=1}\left(1-q^n\right)^{-\frac{1}{2\pi i}X(n)},
\end{equation}
where
\begin{equation}
X(n)=\sum_{\scriptsize
\begin{array}{cc}
	1\leq d\leq n\\
	f(d)|n
\end{array}
\normalsize}\frac{\mu(d)}{d^{s}}\left(n/f(d)\right)^{-\nu}=n^{-\nu}\sum_{\scriptsize
\begin{array}{cc}
	1\leq d\leq n\\
	f(d)|n
\end{array}
\normalsize}\frac{\mu(d)}{d^s}f(d)^{\nu}.
\end{equation}
Also with $f(n)=n^k$, we get
$$
\Pi=\exp\left(\sum^{\infty}_{n,l=1}\frac{q^{ln^k}}{l}\frac{\mu(n)}{n^s}\sigma_{-\nu+1}(l)\right)
=\prod^{\infty}_{n=1}\left(1-q^n\right)^{-X(n)},
$$
where
$$
X(n)=\sum_{\scriptsize
\begin{array}{cc}
	1\leq d\leq n\\
	d^{k}|n
\end{array}
\normalsize}\frac{\mu(d)}{d^{s}}\left(n/d^k\right)^{-\nu}=n^{-\nu}\sum_{\scriptsize
\begin{array}{cc}
	1\leq d\leq n\\
	d^{k}|n
\end{array}
\normalsize}\mu(d)d^{k\nu-s}.
$$
Set now
\begin{equation}
h_{\nu,k}(n)=\sum_{d|n}d^{\nu}\mu_k\left(\frac{n}{d}\right),
\end{equation}
with
\begin{equation}
\mu_{k}(n)=\sum_{d^k|n}\mu(d).
\end{equation}
Also set
\begin{equation}
\mu_{k,v}(n)=\sum_{d^k|n}\frac{\mu(d)}{d^{v}}.
\end{equation}
Then using the identity $\sigma_{\nu}(n)=n^{\nu}\sigma_{-\nu}(n)$, we get with $s=k\nu+v$:
$$
\Pi=\exp\left(\sum^{\infty}_{n,l=1}\frac{q^{n^kl}}{l}\frac{\mu(n)}{n^s}l^{1-\nu}\sigma_{\nu-1}(l)\right)=\exp\left(\sum^{\infty}_{n,l,y=1}\frac{q^{n^kly}}{l}\frac{\mu(n)}{n^{\nu k}n^v}l^{1-\nu}y^{\nu-1}\right)=
$$
$$
=\exp\left(\sum^{\infty}_{n,x,y=1}\frac{q^{n^kxy}}{(xy)^{\nu}}\frac{\mu(n)}{n^{k\nu}n^v}x^{\nu-1}\right)=\exp\left(\sum^{\infty}_{n,x,y=1}\frac{q^{(n^k y) x}}{(n^ky)^{\nu}}\frac{\mu(n)}{xn^v}\right)=
$$
$$
=\exp\left(\sum^{\infty}_{A,x=1}\frac{q^{Ax}}{A^{\nu}x}\sum_{d^{k}|A}\mu(d)d^{-v}\right).
$$
Using definition relations (193)-(195), we get 
$$
\Pi=\exp\left(\sum^{\infty}_{A,x=1}\frac{q^{Ax}}{A^{\nu}x}\mu_{k,v}(A)\right)=\exp\left(\sum^{\infty}_{n=1}\frac{q^n}{n^{\nu}}\sum_{d|n}(n/d)^{\nu-1}\mu_{k,v}(d)\right).
$$
Hence
$$
\exp\left(\sum^{\infty}_{n=1}\frac{q^n}{n}\sum_{d|n}d\frac{1}{d^{\nu}}\mu_{k,v}\left(\frac{n}{d}\right)\right)=\prod^{\infty}_{n=1}\left(1-q^n\right)^{-\mu_{k,v}(n)/n^{\nu}}.
$$
Finally we get\\
\\
\textbf{Lemma 32.}\\
If $|q|<1$, then for any $k,\nu\in\textbf{N}$, $s=k\nu+v$ holds
\begin{equation}
\exp\left(\sum^{\infty}_{n,l=1}\frac{q^{n^kl}}{l^{\nu}}\frac{\mu(n)}{n^s}\sigma_{\nu-1}(l)\right)=\prod^{\infty}_{n=1}\left(1-q^n\right)^{-\mu_{k,v}(n)/n^{\nu}}.
\end{equation}
And as a generalization
\begin{equation}
\exp\left(\sum^{\infty}_{n,l=1}\frac{q^{n^kl}}{l^{\nu}}\frac{h(n)}{n^s}\sigma_{\nu-1}(l)\right)=\prod^{\infty}_{n=1}\left(1-q^n\right)^{-\widehat{h}_{k,v}(n)/n^{\nu}},
\end{equation}
where
\begin{equation}
\widehat{h}_{k,v}(n)=\sum_{d^k|n}\frac{h(d)}{d^v}.
\end{equation}
\\
Hence the next general theorem rises\\
\\
\textbf{Theorem 33.}\\
If $|q|<1$, then for any arithmetical function $f(n)$ holds
\begin{equation}
\exp\left(\sum^{\infty}_{n,l=1}\frac{q^{n^kl}}{n^{v}l}h(n)\sum_{d|l}f(n^kd)d\right)=\prod^{\infty}_{n=1}\left(1-q^n\right)^{-\widehat{h}_{k,v}(n)f(n)}.
\end{equation}
\\
This theorem gives rise to the next more general\\
\\
\textbf{Theorem 34.}\\
If $|q|<1$ and if $g(n):\textbf{N}\rightarrow\textbf{N}$ strictly increasing, then
\begin{equation}
\exp\left(\sum^{\infty}_{n,l=1}\frac{q^{g(n)l}}{l}h(n)\sum_{d|l}f\left(g(n)d\right)d\right)=\prod^{\infty}_{n=1}\left(1-q^n\right)^{-\widetilde{h}_{g}(n)f(n)},
\end{equation} 
where
\begin{equation}
\widetilde{h}_{g}(n)=\sum_{g(d)|n}h(d).
\end{equation}
\\
\textbf{Theorem 35.}\\
There exists a function $\mu_k^*(n)$ such that 
\begin{equation}
\sum_{d|n}\mu_k^*(d)\mu_k(d)=\textbf{1}_n\textbf{, }\forall n,k  \in\textbf{N}.
\end{equation}
Then
\begin{equation}
\mu_{k}^*(n)\mu_k(n)=\mu(n)\textrm{, }\forall n,k\in\textbf{N}
\end{equation}
and $\mu_k^*(n)$ can chosen to be independed of $k$ i.e. $\mu_k^*(n)=\mu(n)$.\\ 
\\
\textbf{Remarks.}\\
In general holds
$$
\exp\left(\sum^{\infty}_{n=1}\frac{q^n}{n}\sum_{d|n}f(d)\right)=\prod^{\infty}_{n=1}\left(1-q^n\right)^{-f(n)/n}.
$$
Assume that $f^*(n)$ is the arithmetic inverse of $f(n)$ in the sense
\begin{equation}
\sum_{d|n}f^*(d)f(d)=\textbf{1}_n\textrm{, }\forall n\in\textbf{N},
\end{equation}
then
\begin{equation}
f^*(n)f(n)=\mu(n).
\end{equation}
For example in case $f(n)=\lambda(n)=\lambda_2(n)$, then $f^*(n)=\mu(n)^2$. Hence 
$$
\lambda_2(n)\mu(n)^2=\mu(n).
$$
In general if (see Theorem 43 below)
$$
\lambda_{\nu}(n)=\sum_{d^{\nu}|n}\mu\left(\frac{n}{d^{\nu}}\right)\textrm{, }\nu=2,3,4,\ldots,
$$
then
$$
\lambda^*_{\nu}(n)=\mu\left(n\right)^2
$$
and
$$
\lambda_{\nu}(n)\mu(n)^2=\mu(n).
$$
\\
\textbf{Proposition 6.}\\
If $f:\textbf{N}\rightarrow\textbf{N}$ is any arithmetical function with the property 
$$
\forall l,n\in\textbf{N}\textrm{ with }l|n\Leftrightarrow f(l)|f(n),
$$
then from (137) setting where $n\rightarrow f(n)$, we get
\begin{equation}
X\left(f(n)\right)=f(n)^{-\nu}\sum_{d|n}\frac{\mu(d)}{d^s}f(d)^{\nu}.
\end{equation}
\\
\textbf{Lemma 36.}\\
Suppose $f,\chi,g,g_1$ are arithmetical functions with $f,g_1:\textbf{N}\rightarrow\textbf{N}$ and $f$ strictly increasing, $g_1$ non-decreasing. Then 
$$
P_1=\exp\left(\sum_{n,l=1}^{\infty}\frac{q^{f(n)(l+g_1(n))}}{l}\chi(n)g(l)\right)=
$$
\begin{equation}
\exp\left(\sum^{\infty}_{n=1}q^n\sum_{\scriptsize
\begin{array}{cc}
	0<d|n\\
	f_i(d)=integer\geq1\\
	n/d-g_1(f_i(d))\geq 1\\
\end{array}
\normalsize}\frac{\chi\left(f_i(d)\right)}{n/d-g_1\left(f_i(d)\right)} g\left(n/d-g_1\left(f_i(d)\right)\right)\right).
\end{equation}
\\

Then also
$$
P_1=\exp\left(\sum^{\infty}_{n=1}\frac{q^n}{n}\sum_{\scriptsize
\begin{array}{cc}
	0<d|n\\
	f_i(d)\in\textbf{\scriptsize N}\\
	n/d-g_1\left(f_i(d)\right)\geq 1
\end{array}
\normalsize}\chi\left(f_i(d)\right)d g\left(n/d-g_1\left(f_i(d)\right)\right)\right)\times
$$
$$
\times\exp\left(\sum^{\infty}_{n=1}\frac{q^n}{n}\sum_{\scriptsize
\begin{array}{cc}
	0<d|n\\
	f_i(d)\in\textbf{\scriptsize N}\\
	n/d-g_1(f_i(d))\geq 1
\end{array}
\normalsize}\frac{\chi\left(f_i(d)\right)g_1\left(f_i\left(d\right)\right)d}{n/d-g_1\left(f_i(d)\right)}g\left(n/d-g_1\left(f_i(d)\right)\right)\right).
$$
Set now
\begin{equation}
h(d,n)=\frac{n}{d}-g_1\left(f_i(d)\right).
\end{equation}
Then
$$
P_1=\exp\left(\sum^{\infty}_{n=1}\frac{q^n}{n}\sum_{\scriptsize
\begin{array}{cc}
0<d|n\\
f_i(d)\in\textbf{\scriptsize N}\\
	h(d,n)\geq 1
\end{array}
\normalsize}\chi\left(d\right)f(d)g\left(h(d,n)\right)\right)\times
$$
$$
\times\exp\left(\sum^{\infty}_{n=1}\frac{q^n}{n}\sum_{\scriptsize
\begin{array}{cc}
	0<d|n\\
	f_i(d)\in\textbf{\scriptsize N}\\
	h(d,n)\geq 1
\end{array}
\normalsize}\frac{\chi\left(d\right)f(d)g_1\left(d\right)}{h(d,n)}g\left(h(d,n)\right)\right).
$$
Or equivalently
$$
P_1=\exp\left(\sum^{\infty}_{n=1}\frac{q^n}{n}\sum_{\scriptsize
\begin{array}{cc}
1\leq k\leq n\\
	h(k,n)\in\textbf{\scriptsize N}
\end{array}
\normalsize}\chi\left(k\right)f(k)g\left(h(k,n)\right)\right)\times
$$
$$
\times\exp\left(\sum^{\infty}_{n=1}\frac{q^n}{n}\sum_{\scriptsize
\begin{array}{cc}
1\leq k\leq n\\
	h(k,n)\in\textbf{\scriptsize N}
\end{array}
\normalsize}\frac{\chi\left(k\right)f(k)g_1\left(k\right)}{h(k,n)}g\left(h(k,n)\right)\right).
$$
Even thought the cases $\chi(n)=(\pm1)^n$, $g(l)=(\pm1)^l l$, $f(n)=An+B$, $g_1(n)=Cn+D$, $A,B,C,D$ integers and $A,C>0$, lead us to Lerch series which are connected to Ramanujan's mock theta functions, the cases of $g(l)=l$ and $\chi(n)$ any arithmetical function is a quite general and provide us with considerably simplified forms. Hence in this case we have\\
\\
\textbf{Theorem 37.}\\If $\chi,f,g_1$ are as in Lemma 36 and $|q|<1$, then 
\begin{equation}
\sum_{n=1}^{\infty}\chi(n)\frac{q^{f(n)g_1(n)}}{1-q^{f(n)}}=\sum^{\infty}_{n=1}q^n\sum_{\scriptsize
\begin{array}{cc}
	0<d|n\\
	f_i(d)=integer\geq1\\
	n/d-g_1(f_i(d))\geq 0\\
\end{array}
\normalsize}\chi\left(f_i(d)\right).
\end{equation}
Hence we can also write
\begin{equation}
\sum_{n=1}^{\infty}\chi(n)\frac{q^{f(n)g_1(n)}}{1-q^{f(n)}}=\sum^{\infty}_{n=1}q^n\sum_{\scriptsize
\begin{array}{cc}
  1\leq k\leq n\\
	f(k)|n\\
	n\geq f(k)g_1(k)\\
\end{array}
\normalsize}\chi\left(k\right).
\end{equation}
\\

Now assume that $f(x)=a_1x+b_1$ and $g_1(x)=c_1x+d_1$ are integer dinomials with positive leading coefficients. Then\\ 
\\
\textbf{Theorem 38.}\\If $\chi(n)$ is any arithmetical function and $a_1,b_1,c_1,d_1$ are integers with $a_1,c_1>0$, $a_1+b_1\geq 1$, $c_1+d_1\geq1$ and $|q|<1$, then
$$
\phi(\chi;q):=\sum^{\infty}_{n=1}\chi(n)\frac{q^{(a_1n+b_1)(c_1n+d_1)}}{1-q^{a_1n+b_1}}=
$$
\begin{equation}
=\sum^{\infty}_{n=1}q^n\sum_{\scriptsize
\begin{array}{cc}
  1\leq k\leq n\\
	(a_1k+b_1)|n\\
	n\geq (a_1k+b_1)(c_1k+d_1)\\
\end{array}
\normalsize}\chi\left(k\right).
\end{equation}
Also in view of Theorem 39 below we have
\begin{equation}
\exp\left(\phi(\chi;q)\right)=\prod^{\infty}_{n=1}\left(1-q^n\right)^{-X(n)},
\end{equation}
where
\begin{equation}
X(n)=\frac{1}{n}\sum_{d|n}\sum_{\delta|(n/d)}d\delta\chi\left(\frac{d-b_1}{a_1}\right)X_{\phi}(d\delta,d)\mu\left(\frac{n/d}{\delta}\right)
\end{equation}
and
\begin{equation}
X_{\phi}(n,k):=\left\{
\begin{array}{cc}
1\textrm{, if }\frac{k-b_1}{a_1}=integer\geq 1\textrm{ and }\frac{n}{k}-
\frac{k-b_1}{a_1}\geq 1\\
0\textrm{, otherwise }
\end{array}
\right\}.
\end{equation}
Provided that sum and product converge.\\
\\
\textbf{Theorem 38.1}\\
$$
\sum^{\infty}_{n,l=1}\frac{q^{f(n)(l+g_1(n))}}{l}\chi(n)g(l)=
$$
$$
=\sum^{\infty}_{n=1}q^n\sum_{\scriptsize
\begin{array}{cc}
	0<d|n\\
	f_i(d)=integer\geq1\\
	n/d-g_1(f_i(d))\geq 1\\
\end{array}
\normalsize}\frac{\chi\left(f_i(d)\right)}{n/d-g_1\left(f_i(d)\right)} g\left(n/d-g_1\left(f_i(d)\right)\right).\eqno{(214.1)}
$$
Hence the number of representations of any positive integer $n$ in the form
$$
f(x)(y+g_1(x))=n,\eqno{(214.2)}
$$
where $x,y$ posotive integers is
$$
r(n)=\sum_{\scriptsize
\begin{array}{cc}
	0<d|n\\
	f_i(d)=integer\geq1\\
	n/d-g_1(f_i(d))\geq 1\\
\end{array}
\normalsize}1.\eqno{(214.3)}
$$
Also we have the next Broecherds product: 
$$
\exp\left(\sum^{\infty}_{n,l=1}\frac{q^{f(n)(l+g_1(n))}}{l}\chi(n)g(l)\right)=\prod^{\infty}_{n=1}(1-q^n)^{-X(n)},\eqno{(214.4)}
$$
where
$$
X(n)=\frac{1}{n}\sum_{\scriptsize \begin{array}{cc}
0<d|n\textrm{\scriptsize, }0<\delta|(n/d)\\
f_i(d)=integer\geq1\\
\delta-g_1(f_i(d))\geq 1
\end{array}\normalsize}d\delta \chi(f_i(d))\frac{g(\delta-g_1(f_i(d)))}{\delta-g_1(f_i(d))}\mu\left(\frac{n/d}{\delta}\right).\eqno{(214.5)}
$$
\\
\textbf{Corollary 38.2}\\
The number of representations of positive integer $n$ in the form
$$
x^{\nu}+xy=n,\eqno{(214.6)}
$$
where $x,y,\nu$ are positive integers, is
$$
r(n)=\sum_{\scriptsize
\begin{array}{cc}
	0<d|n\\
	n/d-d^{\nu-1}\geq 1\\
\end{array}
\normalsize}1.\eqno{(214.7)}
$$
\\
\textbf{Theorem 39.}\\If $\chi,f,g_1$ are as in Lemma 36 and $|q|<1$, then
\begin{equation}
\exp\left(\sum^{\infty}_{n=1}\chi(n)\frac{q^{f(n)g_1(n)}}{1-q^{f(n)}}\right)=\prod^{\infty}_{n=1}\left(1-q^n\right)^{-X(n)},
\end{equation}
where
\begin{equation} X(n)=\frac{1}{n}\sum_{d|n}\sum_{\delta|(n/d)}d\delta\chi(f_i(d))X_0(d\delta,d)\mu\left(\frac{n/d}{\delta}\right)
\end{equation}
and
\begin{equation}
X_0(n,k):=\left\{
\begin{array}{cc}
	1\textrm{, if }f_i(k)=integer\geq1\textrm{ and }n\geq kg_1(f_i(k))\\
0\textrm{, otherwise }
\end{array}
\right\}.
\end{equation}
\\
\textbf{Proof.}\\
We use Lemma 36 with $g(l)\rightarrow l$ and $\chi(n)\rightarrow \chi(n)n$. After that we use Theorem 30 relations (186),(187).\\
\\
\textbf{Theorem 40.}\\
If $z_1,z_2$ are complex numbers in $\textbf{H}$ (the upper half plane) and $f(z)=\sum^{\infty}_{n=0}a_nq^n$, $q=e(z)$, homomorphic also in $\textbf{H}$, then
\begin{equation}
\exp\left(2\pi i \int^{z_2}_{z_1}f(z)dz\right)=\left(\frac{q_2}{q_1}\right)^{a_0}\prod^{\infty}_{n=1}\left(\frac{1-q_2^n}{1-q_1^n}\right)^{-X(n)}, 
\end{equation} 
where $q_j=e(z_j)$, $j=1,2$ and 
\begin{equation}
X(n)=\frac{1}{n}\sum_{d|n}a_d\mu\left(\frac{n}{d}\right).
\end{equation}
\\
\textbf{Example.}\\
It is known that if $q=e^{2\pi i z}$, $Im(z)>0$ then
$$
\sum^{\infty}_{n,m=-\infty}q^{n^2+n m+m^2}=1+6\sum^{\infty}_{n=1}q^n\sum_{d|n}\left(\frac{d}{-3}\right).
$$
Hence from Theorem 40, we have
$$
\exp\left(2\pi i \int^{z_2}_{z_1}\sum^{\infty}_{n,m=-\infty}e^{2\pi i z(n^2+n m+m^2)}dz\right)=
$$
$$
=\frac{q_2}{q_1}\exp\left(6\sum^{\infty}_{n=1}\frac{q_2^n}{n}\sum_{d|n}\left(\frac{d}{-3}\right)-6\sum^{\infty}_{n=1}\frac{q_1^n}{n}\sum_{d|n}\left(\frac{d}{-3}\right)\right)=
$$
$$
=\frac{q_2}{q_1}\prod^{\infty}_{n=1}\left(\frac{1-q_2^n}{1-q_1^n}\right)^{\frac{6}{n}\left(\frac{n}{-3}\right)}.
$$
Hence
\begin{equation}
\exp\left(2\pi i \int^{z_2}_{z_1}\sum^{\infty}_{n,m=-\infty}e^{2\pi i z(n^2+n m+m^2)}dz\right)=\frac{q_2}{q_1}\prod^{\infty}_{n=1}\left(\frac{1-q_2^n}{1-q_1^n}\right)^{\frac{6}{n}\left(\frac{n}{-3}\right)}.
\end{equation}
\\
\textbf{Corollary 41.}\\
Assume that exists a function $f(z)$ and constants $k,N,\epsilon$ such that for all $z\in\textbf{H}$ hold
\begin{equation}
\exp\left(2\pi i \int^{-1/(Nz)}_{z}f\left(w\right)dw\right)=\epsilon z^k.
\end{equation}
Further if
\begin{equation}
\int^{z+1}_{z}f\left(w\right)dw=0,
\end{equation}
then 
\begin{equation}
\exp\left(2\pi i\int^{z}_{c_0}f(w)dw\right)
\end{equation}
is a modular form of weight $k$ in $\Gamma(N)$. Also (from Theorem 40) the function
\begin{equation}
\phi(z)=\prod^{\infty}_{n=1}(1-q^n)^{-X(n)},
\end{equation}  
is modular form of weight $k$ in $\Gamma(N)$.\\
\\

Taking the logarithms and derivating both sides of (221), we can write 
$$
f\left(\frac{-1}{Nz}\right)\frac{1}{Nz^2}-f\left(z\right)=\frac{k}{2\pi i z}.
$$
Hence
$$
f\left(\frac{-1}{Nz}\right)\frac{-1}{N z}+zf\left(z\right)=-\frac{k}{2\pi i}.
$$
If we set 
\begin{equation}
g(z)=-\frac{2\pi i}{k}zf\left(z\right),
\end{equation}
then
\begin{equation}
g\left(\frac{-1}{Nz}\right)+g(z)=1.
\end{equation}
We can write
\begin{equation}
g\left(-\frac{1}{\sqrt{N}z}\right)+g\left(\frac{z}{\sqrt{N}}\right)=1
\end{equation}
and if
\begin{equation}
h(z)=g\left(\frac{z}{\sqrt{N}}\right),
\end{equation}
then
\begin{equation}
h\left(-\frac{1}{z}\right)+h(z)=1.
\end{equation}
Hence we get the next corollary.\\
\\
\textbf{Corollary 42.}\\
Let $h(z)$ be a function such that for all $z\in\textbf{H}$ we have
\begin{equation}
h\left(\frac{-1}{z}\right)+h(z)=1\textrm{,  }\frac{h\left(z+\sqrt{N}\right)}{z+\sqrt{N}}=\frac{h\left(z\right)}{z}.
\end{equation}
Then the function
\begin{equation}
f(z)=-\frac{k}{2\pi i z}h\left(z\sqrt{N}\right)
\end{equation}
have Fourier expansion
\begin{equation}
f(z)=\sum^{\infty}_{n=0}a_nq^n\textrm{, }q=e(z)\textrm{, }Im(z)>0.
\end{equation}
Also if
\begin{equation}
X(n)=\frac{1}{n}\sum_{d|n}a_d\mu\left(\frac{n}{d}\right),
\end{equation}
then the function 
\begin{equation}
\phi(z)=\prod^{\infty}_{n=1}\left(1-q^n\right)^{-X(n)},
\end{equation}
is a modular form of weight $k$ in a certain group $\Gamma(N)$. Also holds the next representation
\begin{equation}
\phi(z)=\exp\left(2\pi i\int^{z}_{i\infty}f(w)dw\right).
\end{equation}
\\
\textbf{Example.}\\
Assume that
\begin{equation}
E_4(z)=q^A\prod^{\infty}_{n=1}\left(1-q^n\right)^{-X(n)}.
\end{equation}
For to evaluate $A$ and $X(n)$, we write $M(q)=E_4(z)=\exp\left(f(q)\right)$, then from (see [22] pg.330):
$$
q\frac{dM}{dq}=\frac{LM-N}{3}
$$
and
$$
\frac{J'(q)}{J(q)}=-\frac{N}{qM},
$$
we get
$$
\frac{M'(q)}{M(q)}=\frac{L}{3q}-\frac{N}{3qM}=\frac{L}{3q}+\frac{J'(q)}{3J(q)},
$$
where $j(z)=J(q)$, $q=e(z)$, $Im(z)>0$ is the $j-$invariant. Hence writing 
$$
\frac{M'(q)}{M(q)}=\sum^{\infty}_{n=1}A_{n-1}q^n,
$$
we have
\begin{equation}
A_{n}=\frac{1}{3}c_n-8\sigma_1(n+1),
\end{equation}
where $c_n$ are the series coefficients of $J'(q)/J(q)$. Hence
$$
M(q)=\exp\left(\sum^{\infty}_{n=1}A_{n-1}\frac{q^n}{n}\right)=\prod^{\infty}_{n=1}\left(1-q^n\right)^{-X(n)},
$$
where
$$
-X(n)=\frac{1}{n}\sum_{d|n}A_{d-1}\mu\left(\frac{n}{d}\right)=8-\frac{1}{3n}\sum_{d|n}c_{d-1}\mu\left(\frac{n}{d}\right).
$$
Hence
\begin{equation}
E_4(z)=\prod^{\infty}_{n=1}\left(1-q^n\right)^{8-\frac{1}{3n}\sum_{d|n}c_{d-1}\mu(n/d)}.
\end{equation}
\\
\textbf{Example.}\\
Assume that $X(n)=1$ and $A=1/24$, then $\phi(z)=\eta(z)$, where $\eta(z)=q^{1/24}\prod^{\infty}_{n=1}\left(1-q^n\right)$ is the Dedekind eta function. This function have modular properties. i.e If $ad-bc=1$, then exist $\epsilon=\epsilon(a,b,c,d)$ and $\epsilon^{24}=1$ such that 
$$
\eta(\sigma(z))=\epsilon(a,b,c,d)(cz+d)^{1/2}\eta(z)\textrm{, }\forall z\in\textbf{H}.
$$
Hence if we assume the function
$$
f(z)=\sum^{\infty}_{n=1}\sigma_1(n)q^n\textrm{, }q=e(z)\textrm{, }Im(z)>0,
$$
then
$$
X(n)=1=\frac{1}{n}\sum_{d|n}\sigma_1(d)\mu(n/d)
$$
and $\psi(z)=\exp\left[2\pi i(z/24-F(z))\right]$, where $F'(z)=f(z)$, behaves exactly as $\eta(z)$ i.e. is a modular form of weight 1/2 and  
$$
\psi(\sigma(z))=\epsilon(a,b,c,d)(cz+d)^{1/2}\psi(z)\textrm{, }\forall z\in\textbf{H}.
$$
Actually it is
$$
\eta(z)=\exp\left[2\pi i\left(z/24-F(z)\right)\right].
$$
Hence in better detail
\begin{equation}
\exp\left(2\pi i\int^{z_2}_{z_1}E_2(z)dz\right)=\frac{\Delta(z_2)}{\Delta(z_1)},
\end{equation}
where $\Delta(z)=\eta(z)^{24}$.\\
More generally if $\nu$ is even positive integer and   
$$
F_{\nu}(z)=-\frac{1}{2\pi i}\frac{2\nu}{B_{\nu}}\sum^{\infty}_{n=1}\sigma_{\nu-1}(n)\frac{q^n}{n},
$$
then
$$
z_2+F_{2\nu}(z_2)-\left(z_1+F_{2\nu}(z_1)\right)=\int^{z_2}_{z_1}E_{2\nu}(z)dz
$$
and
\begin{equation}
\frac{q_2}{q_1}\prod^{\infty}_{n=1}\left(\frac{1-q_2^n}{1-q_1^n}\right)^{4n^{2\nu-2}\nu/B_{2\nu}}=\exp\left(2\pi i\int^{z_2}_{z_1}E_{2\nu}(z)dz\right),
\end{equation}
since 
$$
\frac{1}{n}\sum_{d|n}\sigma_{\nu-1}(d)\mu\left(\frac{n}{d}\right)=n^{\nu-2}.
$$
\\
\textbf{Theorem 43.}\\
If $\phi(z)$ is a modular form of weight $k$ and index $N$ i.e. in $\Gamma(N)$. Then for the function
\begin{equation}
g(z)=-\frac{z}{k}\frac{\phi'(z)}{\phi(z)}
\end{equation}
holds
\begin{equation}
g\left(\frac{-1}{Nz}\right)+g\left(z\right)=1\textrm{, }\forall z\in\textbf{H}
\end{equation}
and $g(z)/z$ is $1-$periodic. The opposite is also true.\\
\\
\textbf{Example.}\\
Let $E_4(z)$ be the Eisenstein series of wight 4. Then set $f(w)$ be such that
$$
E_4(z)=\exp\left(2\pi i\int^{z}_{c_0}f(w)dw\right),
$$
where $c_0$ suitable constant. Solving with respect to $f(z)$ the above equation we get
$$
\frac{d}{dz}\log\left(E_4(z)\right)=2\pi i f(z). 
$$
Hence
$$
f(z)=\frac{1}{2\pi i}\frac{E_4'(z)}{E_4(z)}
$$
and indeed as one can see numerically for the function
$$
g(z)=-\frac{z}{4}\frac{E_4'(z)}{E_4(z)}
$$
holds
$$
g\left(\frac{-1}{z}\right)+g(z)=1
$$
and
$$
\frac{g(z+1)}{z+1}=\frac{g(z)}{z}\textrm{, }\forall z\in\textbf{H}.
$$
\\
\textbf{Example.}\\
If $\lambda(n)$ is Liouville's lambda arithmetical function and if $q_j=e(z_j)$, $Im(z_j)>0$, $j=1,2$, then
\begin{equation}
\exp\left(2\pi i\int^{z_2}_{z_1}\theta_3(z)dz\right)=\frac{q_2}{q_1}\prod^{\infty}_{n=1}\left(\frac{1-q_2^n}{1-q_1^n}\right)^{-2\lambda(n)/n},
\end{equation}
where
\begin{equation}
\theta_3(z)=\sum^{\infty}_{n=-\infty}q^{n^2}\textrm{, }q=e(z)\textrm{, }Im(z)>0.
\end{equation}
\\
\textbf{Proof.}\\
Use Theorem 45 below and the identities
$$
\frac{\theta_{3}(z)-1}{2}=\sum^{\infty}_{n=1}X_2(n)q^n=\sum^{\infty}_{n=1}q^{n^2}\textrm{, }|q|<1,
$$
$$
\lambda(n)=\sum_{d^2|n}\mu\left(\frac{n}{d^2}\right),
$$
where $\mu(n)$ is the Moebious $\mu$ arithmetical function.\\

More generally one can see that\\
\\
\textbf{Theorem 44.}\\
If $q=e(z)$, $Im(z)>0$ and define the next generalization of theta function as
\begin{equation}
\phi_{\nu}(z):=\sum^{\infty}_{n=1}q^{n^{\nu}}\textrm{, }\nu=2,3,4,\ldots
\end{equation}
Then holds
\begin{equation}
\exp\left(2\pi i\int^{z}_{i\infty}\phi_{\nu}(w)dw\right)=\prod^{\infty}_{n=1}\left(1-q^n\right)^{-\lambda_{\nu}(n)/n},
\end{equation}
where 
\begin{equation}
\lambda_{\nu}(n):=\sum_{d^{\nu}|n}\mu\left(\frac{n}{d^{\nu}}\right),
\end{equation}
is the generalized Liouville function.\\
For this function also holds
\begin{equation}
\sum_{d|n}\lambda_{\nu}(d)=X_{\nu}(n):=\left\{
\begin{array}{cc}
	1\textrm{, if }\exists m\in\textbf{N}:n=m^{\nu}\\
	0\textrm{, else }
\end{array}
\right\}.
\end{equation}
Also if $(n,m)=1$, then
\begin{equation}
\lambda_{\nu}(nm)=\lambda_{\nu}(n)\lambda_{\nu}(m).
\end{equation}
\begin{equation}
\lambda_{\nu}\left(n^{\nu}\right)=1.
\end{equation}
\\
\textbf{Theorem 45.}\\
If $\chi(n)$ is any arithmetical function, then
\begin{equation}
\exp\left(\sum^{\infty}_{n=1}\chi(n)q^{n^{\nu}}\right)=\prod^{\infty}_{n=1}\left(1-q^n\right)^{-X(n)},
\end{equation}
where 
\begin{equation}
X(n)=\frac{1}{n}\sum_{d^{\nu}|n}\chi(d)d^{\nu}\mu\left(n/d^{\nu}\right).
\end{equation}
\\
\textbf{Theorem 46.}\\
If $|q|<1$, then
\begin{equation}
\phi_{\nu}(z)=\sum^{\infty}_{n=1}q^{n^{\nu}}=\sum^{\infty}_{n=1}X_{\nu}(n)q^n\textrm{, }|q|<1.
\end{equation}
Then
\begin{equation}
X_{\nu}(nm)X_{\nu}\left(\textrm{gcd}(n,m)\right)=X_{\nu}(n)X_{\nu}(m)\textrm{, }\forall n,m\in\{1,2,\ldots\}
\end{equation}
and $X_{\nu}(n)$ have Dirichlet series
\begin{equation}
L\left(X_{\nu},s\right)=\sum^{\infty}_{n=1}\frac{X_{\nu}(n)}{n^s}=\zeta(\nu s),
\end{equation}
where $\zeta(s)$ is the Riemann's zeta function.\\
\\
\textbf{Proof.}\\
Easy.\\
\\
\textbf{Theorem 47.}\\
If $q=e(z)$, $Im(z)>0$ and $\nu=2,3,\ldots$, then
\begin{equation}
\phi_{\nu}(z)=\sum^{\infty}_{n=1}q^{n^{\nu}}=\sum^{\infty}_{n=1}\frac{\lambda_{\nu}(n)q^n}{1-q^n}.
\end{equation}
In general
\begin{equation}
\sum^{\infty}_{n=1}\chi(n)q^{f(n)}=\sum^{\infty}_{n=1}\frac{X(n)nq^n}{1-q^n},
\end{equation}
where
\begin{equation}
X(n)=\frac{1}{n}\sum_{f(d)|n}\chi(d)\mu\left(\frac{n}{f(d)}\right).
\end{equation}
Also
\begin{equation}
X(n)=\frac{1}{n}\sum_{\scriptsize
\begin{array}{cc}
	d|n\\
f_i(d)=integer\geq1
\end{array}
\normalsize}\chi(f_i(d))\mu\left(\frac{n}{d}\right).
\end{equation} 
\\
\textbf{Proof.}
$$
\phi_{\nu}(z)=\sum^{\infty}_{n=1}X_{\nu}(n)q^n=\sum^{\infty}_{n=1}q^n\sum_{d|n}\lambda_{\nu}(d)=
$$
$$
=\sum^{\infty}_{n,m=1}q^{nm}\lambda_{\nu}(n)=\sum^{\infty}_{n=1}\frac{\lambda_{\nu}(n)q^n}{1-q^n}.
$$
\\
\textbf{Theorem 48.}\\
For every positive integer $n$ we define the function $c_3(n)$ which is 0 when $n$ is cube-free, and $\mu(r)/r$, where $r$ is the cube-free power part of $n$. For example $c_3(2)=0$, $c_3(6)=0$, $c_3(8)=1$, $c_3(12)=0$, $c_3(16)=\mu(2)/2=-1/2$, $c_3(24)=\mu(3)/3=-1/3$, $c_3(2^3 3^3 4)=\mu(4)/4=0,\ldots$etc. Then it holds 
\begin{equation}
\exp\left(\sum^{\infty}_{n=1}q^{n^3}\right)=e^q\prod^{\infty}_{n=1}(1-q^n)^{-c_3(n)}\textrm{, }|q|<1.
\end{equation}
The function $c_3(n)$ is not multiplicative. Also holds
\begin{equation}
c_3(n)=\frac{1}{n}\sum_{d^3|n,d>1}d^3\mu\left(\frac{n}{d^3}\right)
\end{equation}
and the arithmetical function $c_3(n)+\frac{\mu(n)}{n}$ is multiplicative.\\
\\
\textbf{Theorem 49.}\\
For every positive integer $n$ we define the function $c_4(n)$ which is 0 when $n$ is quartic-free power, and $\mu(r)/r$, where $r$ is the quartic-free part of $n$. For example $c_4(2)=0$, $c_4(6)=0$, $c_4(8)=0$, $c_4(12)=0$, $c_4(16)=\mu(1)/1=1$, $c_4(24)=0$, $c_4(2^4\cdot2\cdot3)=\mu(6)/6=1/6$, $c(2^4 3^4 4)=\mu(4)/4=0,\ldots$etc. Then it holds 
\begin{equation}
\exp\left(\sum^{\infty}_{n=1}q^{n^4}\right)=e^q\prod^{\infty}_{n=1}(1-q^n)^{-c_4(n)}\textrm{, }|q|<1.
\end{equation}
The function $c_4(n)$ is not multiplicative. Also holds
\begin{equation}
c_4(n)=\frac{1}{n}\sum_{d^4|n,d>1}d^4\mu\left(\frac{n}{d^4}\right)
\end{equation}
and the arithmetical function $c_4(n)+\frac{\mu(n)}{n}$ is multiplicative.\\
In general
\begin{equation}
\exp\left(\sum^{\infty}_{n=1}q^{n^{\nu}}\right)=e^q\prod^{\infty}_{n=1}(1-q^n)^{-c_{\nu}(n)}\textrm{, }|q|<1.
\end{equation}
Also holds
\begin{equation}
c_{\nu}(n)=\frac{1}{n}\sum_{d^{\nu}|n,d>1}d^{\nu}\mu\left(\frac{n}{d^{\nu}}\right)
\end{equation}
and the arithmetical function $c_{\nu}(n)+\frac{\mu(n)}{n}$ is multiplicative.\\ 
\\

If $|q|<1$ and $\chi(n)$ any arithmetical function, then 
\begin{equation}
\exp\left(\sum^{\infty}_{n=1}\chi(n)q^{n^{\nu}}\right)=e^{\chi(1)q}\prod^{\infty}_{n=1}\left(1-q^n\right)^{-\chi(n_{\nu}(n))\epsilon(n^{*}_{\nu}(n))},
\end{equation}
provided that both parts are convergent and $\epsilon(n)$ is the function defined in (253) below. The function $n^*_{\nu}(n)$ is the $\nu-$th free power part of the number $n$ i.e. if the prime decomposition of $n$ is
\begin{equation}
n=p_1^{a_1}p_2^{a_2}\ldots p_{t}^{a_t}=\left(p_1^{b_1}p_2^{b_2}\ldots p_t^{b_{t}}\right)^{\nu}p_1^{k_1}p_2^{k_2}\ldots p_t^{k_{t}},
\end{equation}
with $p_1<p_2<\ldots<p_{t}$ and $a_i=b_i\nu+k_i$, $0\leq k_1,k_2,\ldots ,k_{t}<\nu$, then 
\begin{equation}
n_{\nu}(n)=p_1^{b_1}p_2^{b_2}\ldots p_t^{b_{t}}\textrm{ and }n^{*}_{\nu}(n)=p_1^{k_1}p_2^{k_2}\ldots p_t^{k_{t}}.
\end{equation}
By convection we set $n_{\nu}(n):=0$ iff $p_1^{b_1}p_2^{b_2}\ldots p_t^{b_{t}}=1$. Also when some $k_{j}$ is greater than 1, we have $\mu\left(n^{*}_{\nu}(n)\right)=0$. Hence
$$
P=\prod^{\infty}_{n=1}\left(1-q^n\right)^{-\chi(n_{\nu}(n))\epsilon(n^*_{\nu}(n))}=
$$
$$
=\prod_{\scriptsize
\begin{array}{cc}
n=p_1^{a_1}p_2^{a_2}\ldots p_t^{a_t}\\
a_i\equiv0,1(\textrm{mod}\nu)	
\end{array}
\normalsize}\left(1-q^n\right)^{-\chi\left(n_{\nu}(n)\right)\epsilon(\overline{n^*_{\nu}}(n))},
$$
where $\overline{n}$ is the free power part of $n$ i.e. if $n=p_1^{a_1}p_2^{a_2}\ldots p_{r}^{a_r}$, $p_1<p_2<\ldots p_r$, is the prime decomposition of $n$, then $\overline{n}=p_1p_2\ldots p_r$.
Also if $a_i=b_i\nu+v_i$, $0\leq v_i\leq \nu-1$, then set $i=i_k$, $k=1,2,\ldots,r$ for those $i$'s we have $v_{i}=v_{i_k}=1$ and $i=j_{k}$, $k=1,2,\ldots,m$, for those $j$'s we have $v_{i}=v_{j_k}=0$. Then
$$
P=\prod_{\scriptsize 
\begin{array}{cc}
n=p_1^{a_1}p_2^{a_2}\ldots p^{a_{t}}_{t}\\
a_i\equiv 0,1(\textrm{mod}\nu)	
\end{array}
\normalsize}\left(1-q^n\right)^{-\chi\left(p_{i_1}^{b_{i_1}}p_{i_2}^{b_{i_2}}\ldots p_{i_r}^{b_{i_{r}}}p_{j_1}^{b_{j_1}}p_{j_2}^{b_{j_2}}\ldots p_{j_m}^{b_{j_m}}\right)\epsilon\left(p_{i_1}p_{i_2}\ldots p_{i_r}\right)}=
$$
$$
=\prod_{\scriptsize 
\begin{array}{cc}
n=p_1^{a_1}p_2^{a_2}\ldots p^{a_{t}}_{t}\\
a_i\equiv 0,1(\textrm{mod}\nu)	
\end{array}
\normalsize}\left(1-q^n\right)^{-\chi\left(\sqrt[\nu]{n/(p_{i_1}p_{i_2}\ldots p_{i_r})}\right)(-1)^{r}/(p_{i_1}p_{i_2}\ldots p_{i_r})}.
$$
Note also that if $n_{\nu}(n)=0$, then we define  $n_{\nu}^{*}(n):=0$. Moreover if $n_{\nu}(n)\neq 0$, then  $n_{\nu}^{*}(n)=n/(n_{\nu}(n))^{\nu}$. Hence we get the next\\
\\
\textbf{Theorem 50.}\\
If $\nu$ is integer greater than 1 and $\chi(n)$ any arithmetical function, then\\
\textbf{i)}
\begin{equation}
\exp\left(\sum^{\infty}_{n=1}\chi(n)q^{n^{\nu}}\right)=e^{\chi(1)q}\prod^{\infty}_{n=1}\left(1-q^n\right)^{-Y_{\nu}(n)}.
\end{equation} 
The arithmetical function $Y_{\nu}(n)$ is 
$$
Y_{\nu}(n)=\frac{1}{n}\sum_{d>1,d^{\nu}|n}\chi(d)d^{\nu}\mu\left(\frac{n}{d^{\nu}}\right)=
$$
\begin{equation}
=\left\{
\begin{array}{cc}
\chi\left(p^{b_{i_1}}_{i_1}p^{b_{i_2}}_{i_2}\ldots p^{b_{i_k}}_{i_k}\right)\frac{(-1)^r}{p_{j_1}p_{j_2}\ldots p_{j_r}}\textrm{ , if }n\in S_{\nu}\\
\chi(m)\textrm{ , if }n=m^{\nu}>1\\
0\textrm{ , else }
\end{array}\right\}=\chi(n_{\nu}(n))\epsilon\left(n^*_{\nu}(n)\right),
\end{equation}
where $S_{\nu}$ is the set of all positive integers $n$ with prime number factorization $n=p_1^{a_1}p_2^{a_2}\ldots p_{t}^{a_t}=\left(p_{i_1}^{b_{i_1}}p_{i_2}^{b_{i_2}}\ldots p_{i_k}^{b_{i_k}}\right)^{\nu}p_{j_1}p_{j_2}\ldots p_{j_r}$. The function $\epsilon(n)$ is defined as 
\begin{equation}
\epsilon(n):=\left\{\begin{array}{cc}
\frac{\mu(n)}{n}\textrm{, if }n\neq 0\\
0\textrm{, else.}
\end{array}\right\}.
\end{equation}
Also (see [23] pp.29 Example and Theorem 66 below). If we define $Z_{N}(z)$ as
$$
\exp\left(\sum^{\infty}_{n=1}\chi(n)e\left(n^{\nu}Z_{N}(z)\right)\right)=
$$
$$
=e^{\chi(1)e\left(Z_{N}(z)\right)}\prod^{\infty}_{n=1}\left(1-e\left(nZ_{N}(z)\right)\right)^{-\chi\left(n_{\nu}(n)\right)\epsilon\left(n^{*}_{\nu}(n)\right)}=
$$
$$
=\exp\left(C+8\pi^2\int^{z\sqrt{N}}_{i\infty}\eta(w)^4dw\right),\eqno{(271.1)}
$$ 
for a certain constant $C$. Then 
$$
Z_N\left(-\frac{1}{N z}\right)=\theta^{{\{\nu\}}{(-1)}}_{\chi}\left(C_0-\theta^{\{\nu\}}_{\chi}\left(Z_{N}(z)\right)\right),\eqno{(271.2)}
$$
where $C_0$ is constant. The function $\theta^{\{\nu\}}_{\chi}(z)=\sum^{\infty}_{n=1}\chi(n)q^{n^{\nu}}$, $q=e(z)$, $Im(z)>0$.\\
\textbf{ii)} If $\chi(n)$ is multiplicative, then $\frac{\mu(n)}{n}+Y_{\nu}(n)$ is multiplicative.\\
\textbf{iii)} The Dirichlet series of $Y_{\nu}(n)$ are
\begin{equation}
L\left(Y_{\nu},s\right)=-\frac{\chi(1)}{\zeta(s+1)}+\frac{L(\chi,s\nu)}{\zeta(s+1)}.
\end{equation}
\\
\textbf{Theorem 51.}\\
If $\nu$ is integer greater than 1 and $\chi(n)$ any arithmetical function, then\\
\textbf{i)}
\begin{equation}
\sum^{\infty}_{n=1}\chi(n)q^{n^{\nu}}=\chi(1)q+\sum^{\infty}_{n=1}\frac{A_{\nu}(n)q^n}{1-q^n}\textrm{, }|q|<1.
\end{equation}
The arithmetical function $A_{\nu}(n)$ is 
$$
A_{\nu}(n)=\sum_{d>1,d^{\nu}|n}\chi(d)\mu\left(\frac{n}{d^{\nu}}\right)=
$$
\begin{equation}
=\left\{
\begin{array}{cc}
(-1)^r\chi\left(p^{b_{i_1}}_{i_1}p^{b_{i_2}}_{i_2}\ldots p^{b_{i_k}}_{i_k}\right)\textrm{ , if }n\in S_{\nu}\\
\chi(m)\textrm{ , if }n=m^{\nu}>1\\
0\textrm{ , else }
\end{array}\right\}=\chi(n_{\nu}(n))\mu(n^{*}_{\nu}(n)),
\end{equation}
where $S_{\nu}$ is the set of all positive integers $n$ with prime number factorization $n=p_1^{a_1}p_2^{a_2}\ldots p_{t}^{a_t}=\left(p_{i_1}^{b_{i_1}}p_{i_2}^{b_{i_2}}\ldots p_{i_k}^{b_{i_k}}\right)^{\nu}p_{j_1}p_{j_2}\ldots p_{j_r}$.\\
\textbf{ii)} When $\chi(n)$ is multiplicative, then the functions 
\begin{equation}
\mu(n)+\chi\left(n_{\nu}(n)\right)\mu\left(n^{*}_{\nu}(n)\right)\textrm{ and }\mu(n)/n+\chi\left(n_{\nu}(n)\right)\epsilon\left(n^{*}_{\nu}(n)\right)
\end{equation}
are also multiplicative.\\
\textbf{iii)} If $\chi(n)$ is even in $\textbf{Z}$, then it holds
\begin{equation}
\sum^{\infty}_{n=-\infty}\chi(n)q^{|n|^{\nu}}=\chi(0)+2\chi(1)q+2\sum^{\infty}_{n=1}\frac{\chi(n_{\nu}(n))\mu(n^{*}_{\nu}(n))q^n}{1-q^n}\textrm{, }|q|<1.
\end{equation}
\textbf{iv)} The Dirichlet series of $A_{\nu}(n)$ is 
\begin{equation}
L\left(A_{\nu},s\right)=\sum^{\infty}_{n=1}\frac{A_{\nu}(n)}{n^s}=-\frac{\chi(1)}{\zeta(s+1)}+\frac{L(\chi,\nu s+\nu)}{\zeta(s+1)}.
\end{equation}
\\
\textbf{Theorem 51.1}\\
If $n,\nu\in\textbf{N}$ and $\nu>1$, then
$$
n^{(1)}_{\nu}(n)=\sum_{d|n}\varphi_{\nu}(d),\eqno{(277.1)}
$$
where $\varphi_{\nu}(n)=\varphi(\sqrt[\nu]{n})$ if $n=m^{\nu}$ for some $m\in\textbf{N}$ and $0$ else. $\varphi(n)$ is Euler's totient arithmetical function.\\
\textbf{Remarks.}\\
\textbf{i)} For $n^{(1)}_{\nu}(n)$ see also relation (294).\\
\textbf{ii)} $n^{(1)}_{\nu}(n)$ is multiplicative function.\\
\\
\textbf{Proof.}\\
If $\nu\in\textbf{N}$ and $\nu\geq 2$, we have the next Euler product 
$$
\sum^{\infty}_{n=1}\frac{n^{(1)}_{\nu}(n)}{n^s}=\prod_{p-prime}\left(1+\frac{n^{(1)}_{\nu}(p)}{p^s}+\frac{n^{(1)}_{\nu}(p^2)}{p^{2s}}+\ldots\right)=\ldots=
$$
$$
=\prod_{p-prime}\left(1+\frac{1}{p^s}+\frac{1}{p^{2s}}+\ldots+\frac{1}{p^{(\nu-1)s}}\right)\frac{1}{1-\frac{p}{p^{\nu s}}}=\frac{\zeta(s)\zeta(s\nu-1)}{\zeta(s\nu)},
$$
since $n^{(1)}_{\nu}\left(p^m\right)=n^{(1)}_{\nu}(p^k)$, where $m=\nu k+v$, $0\leq v\leq \nu-1$. Hence using
$$
\frac{\zeta(s\nu-1)}{\zeta(s\nu)}=\sum^{\infty}_{n=1}\frac{\varphi_{\nu}(n)}{n^s},
$$
we get the result.\\
\\
\textbf{Corollary 51.2}\\
We have
$$
n_{\nu}(n)=\left\{\begin{array}{cc}\sum_{d|n}\varphi_{\nu}(d)\textrm{, if }n^{(1)}_{\nu}(n)\neq 1\\0\textrm{, if }n^{(1)}_{\nu}(n)=1\end{array}\right\}
$$
and
$$
n^{*}_{\nu}(n)=\left\{\begin{array}{cc}\frac{n}{(n_{\nu}(n))^{\nu}}\textrm{, if }n_{\nu}(n)\neq0\\0\textrm{, if }n_{\nu}(n)=0\end{array}\right\}.
$$
\\
\textbf{Corollary 52.}
$$
\exp\left(\sum^{\infty}_{n=1}\chi(n)q^{n^2}\right)=e^{\chi(1)q}\prod^{\infty}_{n=1}\left(1-q^n\right)^{-\chi\left(n_2(n)\right)\epsilon\left(n^{*}_2(n)\right)}.
$$
Also
\begin{equation}
\exp\left(\sum^{\infty}_{n=1}\chi(n)q^{n^2}\right)=e^{\chi(1)q}\prod^{\infty}_{
n=1}\left(1-q^n\right)^{-Y(n)}.
\end{equation}
The function $Y(n)$ is 
$$
Y(n)=\frac{1}{n}\sum_{d>1,d^2|n}\chi(d)d^2\mu\left(\frac{n}{d^2}\right)=
$$
\begin{equation}
=\left\{\begin{array}{cc}
\chi\left(p^{b_{i_1}}_{i_1}p^{b_{i_2}}_{i_2}\ldots p^{b_{i_k}}_{i_k}\right)\frac{(-1)^r}{p_{j_1}p_{j_2}\ldots p_{j_r}}\textrm{ , if }n\in S\\
\chi(m)\textrm{ , if }n=m^2>1\\
0\textrm{ , else }\\
\end{array}
\right\},
\end{equation}
where $S$ is the set of all positive integers $n$ with prime number factorization $n=p_1^{a_1}p_2^{a_2}\ldots p_{t}^{a_t}=\left(p_{i_1}^{b_{i_1}}p_{i_2}^{b_{i_2}}\ldots p_{i_k}^{a_{i_k}}\right)^2p_{j_1}p_{j_2}\ldots p_{j_r}$. Also if $\chi(n)$ is multiplicative, then $\frac{\mu(n)}{n}+Y(n)$ is multiplicative.\\
\\
\textbf{Corollary 52.1}\\
If $\nu\in\textbf{N}$ and $\nu>1$, $q=e(z)$, $Im(z)>0$, then
$$
\exp\left(\sum^{\infty}_{n=1}\chi(n)q^{n^{\nu}}\right)=e^{\chi(1)q}\prod^{\infty}_{n=1}\left(1-q^n\right)^{-\chi\left(n_{\nu}(n)\right)\epsilon\left(n^{*}_{\nu}(n)\right)},
$$
where 
$$
\epsilon(n)=\left\{\begin{array}{cc}\frac{\mu(n)}{n}\textrm{, if }n\neq 0\\0\textrm{, if }n=0\end{array}\right\}.
$$
\\
\textbf{Corollary 53.}\\
Assume that $\nu=2$ and $\chi(n)=2$, then
\begin{equation}
2\sum^{\infty}_{n=1}q^{n^2}=\sum^{\infty}_{n=-\infty}q^{n^2}-1=\theta_3(q)-1=2q+\sum^{\infty}_{n=1}\frac{A_{2}(n)q^n}{1-q^n},
\end{equation}
where $\chi(n)=2$. Hence $A_{2}(n)=2(-1)^r$, if $r=r(n)$ is such that: $n=p_1^{a_1}p_2^{a_2}\ldots p_t^{a_t}$ and $n=\left(p^{b_{i_1}}_{i_1}p_{i_2}^{b_{i_2}}\ldots p_{i_k}^{b_{i_k}}\right)^2p_{j_1}p_{j_2}\ldots p_{j_r}$ and 0 otherwise. Hence
\begin{equation}
\theta_3(q)=1+2q+2\sum^{\infty}_{n=1}\frac{\mu_2(n)q^n}{1-q^n},
\end{equation}
where 
\begin{equation}
\mu_2(n):=\left\{\begin{array}{cc}
(-1)^r\textrm{ , if }n_2(n)>1\textrm{ and when }n_2^{*}(n)=p_{j_1}p_{j_2}\ldots p_{j_r}\\
0\textrm{ , else }
\end{array}\right\}.
\end{equation} 
Of course $\mu(n_2^{*}(n))\neq 0$ iff $n_2^*(n)$ is square free and $n_{2}(n)>1$. Also $n=(n_2(n))^2n_2^{*}(n)$, if $n_2(n)>1$. The same happens and with the other higher values of $\nu$.\\
\\
In the same way as above we have the next\\
\\ 
\textbf{Theorem 54.}\\
In general for any integer $\nu\geq 2$ we denote
\begin{equation}
\phi_{\nu}(q):=\sum^{\infty}_{n=-\infty}q^{|n|^{\nu}}\textrm{, }|q|<1.
\end{equation}
Then
\begin{equation}
\phi_{\nu}(q)=1+2q+2\sum^{\infty}_{n=1}\frac{\mu_{\nu}(n)q^n}{1-q^n},
\end{equation}
where
\begin{equation}
\mu_{\nu}(n):=\left\{\begin{array}{cc}
(-1)^{r}\textrm{ , if }n_{\nu}(n)>1\textrm{ and }n_{\nu}^{*}(n)=p_{j_1}p_{j_2}\ldots p_{j_r}\\
0\textrm{ , else }
\end{array}\right\}=\mu(n^{*}_{\nu}(n)).
\end{equation} 
Here $n_{\nu}^{*}(n)$ is the $\nu-$th power free part of $n$ and $n=(n_{\nu}(n))^{\nu}n_{\nu}^{*}(n)$, if $n_{\nu}(n)>1$.\\
\\
\textbf{Remarks.}\\ 
\textbf{i)} If the $\nu-$th power free part of $n$ is not  square free, then from the definition of $\mu_{\nu}(n)$, we have $\mu_{\nu}(n)=0$.\\
\textbf{ii)} The function 
\begin{equation}
\mu(n)+\mu_{\nu}(n)
\end{equation} 
is multiplicative.\\
\\
\textbf{Corollary 55.}\\Assume $\nu=2$ and $\chi(n)=2(-1)^n$, then
\begin{equation}
\theta_4(q)=1-2q+2\sum^{\infty}_{n=1}\frac{\mu^{*}_2(n)q^{n}}{1-q^{n}},
\end{equation}
where
\begin{equation}
\mu^{*}_2(n)=\left\{\begin{array}{cc}
(-1)^{p_{i_1}p_{i_2}\ldots p_{i_k}}\mu_2(n)\textrm{ , if }n=\left(p^{b_{i_1}}_{i_1}p^{b_{i_2}}_{i_2}\ldots p^{b_{i_k}}_{i_k}\right)^2p_{j_1}p_{j_2}\ldots p_{j_r}\\
0\textrm{ , else  }
\end{array}\right\}.
\end{equation}
Also $\mu_2^{*}(n)=(-1)^{n_2(n)}\mu(n^{*}_2(n))$.\\
\\
\textbf{Theorem 56.}\\
Set
\begin{equation}
\phi^{*}_{\nu}(q)=\sum^{\infty}_{n=-\infty}(-1)^n q^{|n|^{\nu}}.
\end{equation}
Then
\begin{equation}
\phi^{*}_{\nu}(q)=1-2q+2\sum^{\infty}_{n=1}\frac{\mu^{*}_{\nu}(n)q^n}{1-q^n},
\end{equation}
where
\begin{equation}
\mu^{*}_{\nu}(n):=\left\{\begin{array}{cc}
(-1)^{p_{i_1}p_{i_2}\ldots p_{i_k}}\mu_{\nu}(n)\textrm{ , if }n=\left(p^{b_{i_1}}_{i_1}p^{b_{i_2}}_{i_2}\ldots p^{b_{i_k}}_{i_k}\right)^{\nu}p_{j_1}p_{j_2}\ldots p_{j_{r}}\\
0\textrm{ , else  }
\end{array}
\right\}
\end{equation}
and $\mu_{\nu}(n)$ as defined in Theorem 54.\\
\\
\textbf{Remarks.}\\It is clear that if $2|n_{\nu}(n)$, then $\mu^{*}_{\nu}(n)=\mu_{\nu}(n)$ and if 2 does not divides  $n_{\nu}(n)$, then $\mu^{*}_{\nu}(n)=-\mu_{\nu}(n)$. Hence
\begin{equation}
\mu^{*}_{\nu}(n)=(-1)^{n_{\nu}(n)}\mu_{\nu}(n)
\end{equation}
and we have
$$
\sum^{\infty}_{n=-\infty}(-1)^nq^{|n|^{\nu}}=1-2q+2\sum^{\infty}_{n=1}\frac{(-1)^{n_{\nu}(n)}\mu(n^{*}_{\nu}(n))q^n}{1-q^n}=
$$
\begin{equation}
=1-2q+2\sum_{\scriptsize\begin{array}{cc}
n\geq1\\
n_{\nu}=even
\end{array}\normalsize}\frac{\mu_{\nu}(n)q^n}{1-q^n}-2\sum_{\scriptsize\begin{array}{cc}
n\geq 1\\
n_{\nu}=odd
\end{array}\normalsize}\frac{\mu_{\nu}(n)q^n}{1-q^n}.
\end{equation}
\\

Assume now that $\chi(n)$ is any arithmetical function. Assume also that $n^{(1)}_{\nu}(n)=1$, if $n_{\nu}(n)=0$ and $n^{(1)}_{\nu}(n)=n_{\nu}(n)$, else. Also $n^{(2)}_{\nu}(n)=n/n^{(1)}_{\nu}(n)$ Then
\begin{equation}
A_{\nu}(n)=\chi\left(n^{(1)}_{\nu}(n)\right)\mu\left(n^{(2)}_{\nu}(n)\right).
\end{equation} 
The functions $n^{(1)}_{\nu}(n)$ and $n^{(2)}_{\nu}(n)$, are multilicative.\\
Hence
\begin{equation}
\sum^{\infty}_{n=1}\chi(n)q^{n^{\nu}}=\chi(1)q+\sum^{\infty}_{n=1}\frac{\chi\left(n_{\nu}(n)\right)\mu\left(n^{*}_{\nu}(n)\right)q^n}{1-q^n}.
\end{equation}
Hence  
$$
\sum^{\infty}_{n,m=1}\chi(n)\chi(m)q^{n^{\nu}+m^{\nu}}=\sum^{\infty}_{l=1}\left(\sum_{n^{\nu}+m^{\nu}=l}\chi(n)\chi(m)\right)q^l=
$$
$$
=\left(\chi(1)q+\sum^{\infty}_{n,m=1}\chi\left(n_{\nu}(n)\right)\mu\left(n^{*}_{\nu}(n)\right)q^{mn}\right)^2=
$$
$$
=\left(\sum^{\infty}_{l=1}q^l\sum_{d|l}\chi\left(n^{(1)}_{\nu}(d)\right)\mu\left(n^{(2)}_{\nu}(d)\right)\right)^2.
$$
Hence we get the next\\
\\
\textbf{Theorem 57.}\\
If $\chi(n)$ is any arithmetical function, then
\begin{equation}
\sum_{\scriptsize\begin{array}{cc}
n^{\nu}+m^{\nu}=l\\
n\geq 1,m\geq 1\\
\end{array}
\normalsize}\chi(n)\chi(m)=\sum^{l}_{t=0}A^{*}_{\nu}(t)A^{*}_{\nu}(l-t),
\end{equation}
where
\begin{equation}
A^{*}_{\nu}(t)=\sum_{d|t}\chi\left(n^{(1)}_{\nu}(d)\right)\mu\left(n^{(2)}_{\nu}(d)\right)\textrm{, }A^{*}_{\nu}(0)=0.
\end{equation}
\\

Assuming $\chi(n)=1$, we get
\begin{equation}
\sum_{\scriptsize\begin{array}{cc}
n^{\nu}+m^{\nu}=l\\
n\geq 1,m\geq 1\\
\end{array}
\normalsize}1=\sum^{l}_{t=0}A^{*}_{\nu}(t)A^{*}_{\nu}(l-t),
\end{equation}
where
\begin{equation}
A^{*}_{\nu}(t)=\sum_{d|t}\mu\left(n^{(2)}_{\nu}(d)\right)\textrm{, }A^{*}_{\nu}(0)=0.
\end{equation}
The Fermat$-$Wiles theorem states that
\begin{equation}
\sum^{l^{\nu}}_{t=0}A^{*}_{\nu}(t)A^{*}_{\nu}(l^{\nu}-t)=0.
\end{equation} 
\\
\textbf{Corollary 58.}\\
If $\chi(n)$ is completely multiplicative, then
\begin{equation}
\sum_{\scriptsize\begin{array}{cc}
n^{\nu}+m^{\nu}=l\\
n\geq 1,m\geq 1\\
\end{array}
\normalsize}\chi(nm)=\sum^{l}_{t=0}A^{*}_{\nu}(t)A^{*}_{\nu}(l-t),
\end{equation}
where
\begin{equation}
A^{*}_{\nu}(t)=\sum_{d|t}\chi\left(n^{(1)}_{\nu}(d)\right)\mu\left(n^{(2)}_{\nu}(d)\right)\textrm{, }A^{*}_{\nu}(0)=0.
\end{equation}
\\
\textbf{Theorem 58.1.}\\
For given positive integers $l,\nu$, with $\nu>1$, the number of solutions of
$$
n^{\nu}+m^{\nu}=l\textrm{, }n,m\in \textbf{N},
$$
is 
$$
\sum^{l}_{t=0}A^{*}_{\nu}(t)A^{*}_{\nu}(l-t),
$$
where 
$$
A^{*}_{\nu}(t)=\sum_{d|t}\mu\left(\frac
{d}{\left(\sum_{\delta|d}\varphi_{\nu}(\delta)\right)^{\nu}}\right).
$$
\\

Assume the arithmetical function $f:\textbf{N}\rightarrow\textbf{N}$, with $f(1)=1$ and the  equation $f(x)y=n:(eq)$, $x,y,n$ integers greater than 1. Assume also that $n$, $n'$ are integers greater or equal to 1, such that equations $f(x)y=n$ and $f(x)y=n'$ have solutions
$$
Sol(n)=\{(x_1,y_1),(x_2,y_2),\ldots,(x_k,y_k)\}
$$ 
and
\begin{equation}
Sol(n')=\{(x'_1,y'_1),(x'_2,y'_2),\ldots,(x'_m,y'_m)\}
\end{equation}
respectively. Then 
\begin{equation}
A_f(n)=\sum^{k}_{l=1}\chi(x_l)\mu(y_l)
\end{equation}
and
\begin{equation}
A_f(n')=\sum^{m}_{l=1}\chi(x'_l)\mu(y'_l),
\end{equation}
where
\begin{equation}
A_f(n):=\sum_{d>1\textrm{\scriptsize, \normalsize}f(d)|n}\chi(d)\mu\left(\frac{n}{f(d)}\right).
\end{equation}
\\
\textbf{Theorem 59.}\\
Given an integer number $n\geq 1$ we assume the arithmetical functions $f(n),g(n)\in\textbf{N}$ such that $f,g$ are increasing and $f(1)=g(1)=1$. Then if the solutions of $f(x)y=n$ and $g(x)y=n$, $x,y$ positive integers, are $x=[x_f]=[x_f]_k(n)$, $y=[x_f]^{*}=[x_f]^{*}_k(n)$, $k=1,2,\ldots,r$ and $x=[x_g]=[x_g]_k(n)$, $y=[y_g]^*=[x_g]^{*}_k(n)$, $k=1,2,\ldots,m$, we have  
\begin{equation}
A_f(n):=\sum_{d>1\textrm{\scriptsize, \normalsize} f(d)|n}\chi(d)\mu\left(\frac{n}{f(d)}\right)=\sum^{r}_{k=1}\chi\left([x_f]_k(n)\right)\mu\left([x_f]^{*}_k(n)\right)
\end{equation}
and
\begin{equation}
A_g(n):=\sum_{d>1\textrm{\scriptsize, \normalsize} g(d)|n}\psi(d)\mu\left(\frac{n}{g(d)}\right)=\sum^{m}_{k=1}\psi\left([x_g]_k(n)\right)\mu\left([x_g]^{*}_k(n)\right),
\end{equation}
where $\chi(n)$, $\psi(n)$ are any arithmetical functions. Also then
\begin{equation}
\sum_{\scriptsize\begin{array}{cc}
f(n)g(m)=l\\
n\geq 1,m\geq 1
\end{array}\normalsize}\chi(n)\psi(m)=\sum_{d|l}A^{*}_{f}(d)A^{*}_{g}(l/d),
\end{equation}
where
\begin{equation}
A^{*}_{f}(t)=\sum_{d|t}A_f(d)
\end{equation}
and
\begin{equation}
A^{*}_{g}(t)=\sum_{d|t}A_g(d).
\end{equation}
In case of $f(n)=n^{\nu}$, $g(n)=n^{\lambda}$, then
\begin{equation}
A_{x^{\nu}}(n)=A_{\nu}(n)=\chi\left(n^{(1)}_{\nu}(n)\right)\mu\left(n^{(2)}_{\nu}(n)\right)
\end{equation}
\\
\textbf{Example.}\\
Suppose that $f(n)=2n^2-n$, $g(n)=n^3$ and $\chi(n)=\psi(n)=1$, then
$$
\sum_{\scriptsize\begin{array}{cc}
(2n^2-n)m^3=l\\
n\geq 1,m\geq 1
\end{array}
\normalsize}1=\sum_{d|l}A^{*}_{f}(d)A^{*}_{g}(l/d),
$$
where 
$$
A^{*}_{f}(t)=\sum_{d|t}A_f(d)
$$
and
$$
A^{*}_{g}(t)=\sum_{d|t}A_g(d),
$$
$$
A_f(n)=\sum^{r}_{k=1}\mu([x_f]^*_k(k))
$$
$$
A_{g}(n)=\mu(n^{(1)}_{3}(n))
$$
where $y=[y_f]=[x_f]^{*}_k(n)$, $y=[y_g]=[x_{g}]^{*}_k(n)=n^{(1)}_{3}(n)$ are defined as solutions of $(2X^2-X)Y=n$ and $X^3Y=n$, with $[y_f]\geq 1$, $[y_g]\geq 1$ respectively ($n$ being greater or equal to 1).\\
\\
\textbf{Theorem 60.}\\
Let $\nu$ be integer greater than 1 and $|q|<1$. Let also $f:\textbf{N}\rightarrow\textbf{N}$ be strictly increasing arithmetical function, with $f(0)=0$ and $\chi(n)$ any arithmetical function. Then 
\begin{equation}
\sum^{\infty}_{n=1}\chi(n)q^{f(n)}=\sum^{\infty}_{n=1}q^n\sum_{d|n}\sum_{f(\delta)|d}\chi(\delta)\mu\left(\frac{d}{f(\delta)}\right).
\end{equation}
\begin{equation}
\sum^{\infty}_{n=1}\chi(n)q^{n^{\nu}}=\sum^{\infty}_{n=1}q^n\sum_{d|n}\chi(n^{(1)}_{\nu}(d))\mu(n^{(2)}_{\nu}(d)).
\end{equation}
Also
\begin{equation}
\sum_{\scriptsize\begin{array}{cc}
d|n\\
f(d)|n
\end{array}\normalsize}\chi(d)\mu\left(\frac{n}{f(d)}\right)=\sum_{d|n}\mu\left(\frac{n}{d}\right)\sum_{\scriptsize\begin{array}{cc}
\delta|d\\ 
f(\delta)=d
\end{array}\normalsize}\chi(\delta)
\end{equation}
and
\begin{equation}
\sum_{d|n}\chi(n^{(1)}_{\nu}(d))\mu(n^{(2)}_{\nu}(d))=\left\{\begin{array}{cc}
\chi(m)\textrm{, if }n=m^{\nu}\\
0\textrm{, else }
\end{array}\right\}.
\end{equation}
Also if we set   
\begin{equation}
\phi_{\nu}(\chi;z):=\sum^{\infty}_{n=1}\chi(n)q^{n^{\nu}}\textrm{, }q=e(z)\textrm{, }Im(z)>0,
\end{equation}
then
\begin{equation}
\exp\left(2\pi i\int_{z}^{i\infty}\phi_{\nu}(\chi;w)dw\right)=e^{(1-\chi(1))q}\prod^{\infty}_{n=1}(1-q^n)^{C_{\nu}\left(\chi;n\right)/n},
\end{equation}
where 
\begin{equation}
C_{\nu}\left(\chi;n\right)=\mu(n)+\chi(n_{\nu}(n))\mu(n^{*}_{\nu}(n))
\end{equation}
and $C_{\nu}\left(\chi;n\right)$ is multiplicative function if $\chi(n)$ is multiplicative.\\
\\
\textbf{Theorem 61.}\\
Suppose that $R(n)$ are the number of representations of a positive integer in the form 
\begin{equation}
f(x)+g(y)=n\textrm{, }x,y\in\textbf{N}.
\end{equation} 
If $f,g$ are positive integer valued increasing functions, with $f(0)=g(0)=0$, then
\begin{equation}
R(n)=\sum^{n}_{k=1}A_{f}(k)A_{g}(n-k),
\end{equation}
where
\begin{equation}
A_{f}(n)=\sum_{d|n}\sum_{f(\delta)|d}\mu\left(\frac{d}{f(\delta)}\right)=\sum_{\scriptsize\begin{array}{cc}
d|n\\
f(d)=n
\end{array}\normalsize}1
\end{equation}
and
\begin{equation}
A_{g}(n)=\sum_{d|n}\sum_{g(\delta)|d}\mu\left(\frac{d}{g(\delta)}\right)=\sum_{\scriptsize\begin{array}{cc}
d|n\\
g(d)=n
\end{array}\normalsize}1.
\end{equation}
\\
\textbf{Theorem 62.}\\
Assume the polynomial functions $f,g:\textbf{N}\rightarrow\textbf{N}$, with $f(0)=g(0)=0$ and increasing in $\textbf{N}$. Then the equation
\begin{equation}
f(x)g(y)=n\textrm{, where }x,y\in \textbf{N},
\end{equation}
have
\begin{equation}
\sum_{d|n}A_f(d)A_g\left(\frac{n}{d}\right)
\end{equation}
solutions. Where 
\begin{equation}
A_f(n)=\sum_{\scriptsize\begin{array}{cc}
d|n\\
f(d)=n
\end{array}\normalsize}1\textrm{, }A_g(n)=\sum_{\scriptsize\begin{array}{cc}
d|n\\
g(d)=n
\end{array}\normalsize}1.
\end{equation}
\\
\textbf{Theorem 63.}\\
Assume the polynomial $f(x,y)$ with coefficients in $\textbf{Z}$ and $\textrm{deg}f>0$. Assume also that $f(x,0)=0$, for all $x\in\textbf{N}$ and $f(0,y)=0$, for all $y\in\textbf{N}$. Then\\
$\textbf{1)}$ The equation
\begin{equation}
f(x,y)=n\textrm{, where }x,y\in\textbf{N}, 
\end{equation}
have
\begin{equation}
\sum_{\scriptsize\begin{array}{cc}
d|n\textrm{\scriptsize, }\delta|n\\
f(d,\delta)=n
\end{array}\normalsize
}1
\end{equation}
solutions.\\
$\textbf{2)}$ It holds
\begin{equation}
\sum^{\infty}_{n,m=1}\chi(n,m)q^{f(n,m)}=\sum^{\infty}_{n=1}q^n\sum_{\scriptsize\begin{array}{cc}
d|n\textrm{\scriptsize , \normalsize}\delta|n\\
f(d,\delta)=n
\end{array}\normalsize}\chi(d,\delta).
\end{equation}
\\
\textbf{Theorem 64.}\\
Let $\chi(n),\psi(n)$ be arithmetical functions. Let also $f,g:\textbf{N}\rightarrow\textbf{N}$ are any increasing arithmetical functions with $f(0)=g(0)=0$. We set
\begin{equation}
A(\chi,f;n):=\sum_{d|n}\sum_{f(\delta)|d}\chi(\delta)\mu\left(\frac{d}{f(\delta)}\right)
\end{equation}
and
\begin{equation}
B(\psi,g;n):=\sum_{d|n}\sum_{g(\delta)|d}\psi(\delta)\mu\left(\frac{d}{g(\delta)}\right),
\end{equation}
where $\mu(n)$ is the Moebius function. Then for $|q|<1$ holds
$$
\left(\sum^{\infty}_{n=1}\chi(n)q^{f(n)}\right)\left(\sum^{\infty}_{n=1}\psi(n)q^{g(n)}\right)=
$$
\begin{equation}
=\sum^{\infty}_{n=1}q^n\sum_{\scriptsize\begin{array}{cc}
n_1+n_2=n\\
n_1,n_2>0
\end{array}\normalsize}A(\chi,f;n_1)B(\psi,g;n_2).
\end{equation}  
\\
\textbf{Notes.}\\
The function $n^{(2)}_{\nu}(n)$, have the following properties:\\ 
\textbf{i)} Is multiplicative i.e. If $n,m$ are positive integrs with $(n,m)=1$, then
\begin{equation}
n^{(2)}_{\nu}(nm)=n^{(2)}_{\nu}(n)n^{(2)}_{\nu}(m).
\end{equation}
\textbf{ii)} If $n_{\nu}(n)>1$, then $n^{(2)}_{\nu}(n)=n^{*}_{\nu}(n)$ and when $n_{\nu}(n)=1$, then $n^{(2)}_{\nu}(n)=n$.\\
\textbf{iii)}
\begin{equation}
\sum_{d>1, d^{\nu}|n}\chi(d)\mu\left(\frac{n}{d^{\nu}}\right)=\chi\left(n^{(1)}_{\nu}(n)\right)\mu\left(n^{(2)}_{\nu}(n)\right)\textrm{, }\forall n\in\textbf{N}.
\end{equation}
\textbf{iv)} The number of representations $r_{\nu}(n)$ of a positive integer $n$ in the form
\begin{equation}
ax^{\nu}+by^{\nu}=n\textrm{, }x,y\geq1,
\end{equation}
where $a,b$ positive integers with $(a,b)=1$, is
\begin{equation}
r_{\nu}(n)=\sum_{\scriptsize\begin{array}{cc}
n_1,n_2>0\\
n_1+n_2=n\\
n^{(2)}_{\nu}(n_1)=a\textrm{\scriptsize,  \normalsize}n^{(2)}_{\nu}(n_2)=b
\end{array}
\normalsize}1.
\end{equation} 
\textbf{v)} The case of the represetation of $n\in\textbf{N}$ in a form
$$P_1(x)+P_2(y)=n\textrm{, }x,y\in\textbf{Z},$$ with $P_1(0),P_2(0)\neq0$ can be settled as follows: We assume that $P_1(x),P_2(x)$ are polynomials as in Theorems 20,22 except, (perhaps), that $P_1(0),P_2(0)\neq0$. We write     
\begin{equation}
R_{\{1,2\}}(n)=\sum_{\scriptsize
\begin{array}{cc}
	d\neq 0\textrm{\scriptsize, }\textrm{\scriptsize abs}(d)|n\\
	P^*_{\{1,2\}}(d)=\textrm{\scriptsize abs}(n)
\end{array}
\normalsize}1\textrm{, when }n=1,2,\ldots
\end{equation}
and $R_1(0)=r_{1}$, $R_2(0)=r_{2}$, where $r_{1}$, $r_{2}$ are the number of integer solutions of $P^*_1(n),P^*_2(n)=0$ in $\textbf{Z}$, where  $P^*_1(n)=P_1(n)-P_1(0)$ and $P^*_2(n)=P_2(n)-P_2(0)$ respectively. Then the number of representations is, setting
\begin{equation}
R_{12}(n)=\sum^{n}_{l=0}R_1(l)R_2(n-l)\textrm{, }n=1,2,\ldots,
\end{equation}
exactly $R_{12}\left(n-P_1(0)-P_2(0)\right)$.\\
\textbf{vi)} In the case of the represetation of $n\in\textbf{N}$ in the form
$$P_1(x)P_2(y)=n\textrm{, }x,y\in\textbf{Z},$$ with $P_1(x),P_2(x)$ polynomials as in Theorems 20,22, except, (perhaps), that $P_1(0),P_2(0)\neq0$: Write again        
\begin{equation}
R_{\{1,2\}}(n)=\sum_{\scriptsize
\begin{array}{cc}
	d\neq 0\textrm{\scriptsize, }\textrm{\scriptsize abs}(d)|n\\
	P^*_{\{1,2\}}(d)=n
\end{array}
\normalsize}1\textrm{, when }n=1,2,\ldots
\end{equation}
and $R_1(0)=r_{1}$, $R_2(0)=r_{2}$, where $r_{1}$, $r_{2}$ are the integer solutions of $P^*_1(n),P^*_2(n)=0$ in $\textbf{Z}$, where  $P^*_1(n)=P_1(n)-P_1(0)$ and $P^*_2(n)=P_2(n)-P_2(0)$ respectively. Then the number of representations is
\begin{equation}
R_{12}(n)=\sum_{0<\textrm{\scriptsize abs}(d)|n}R_1\left(d-P_1(0)\right)R_2\left(\frac{n}{d}-P_2(0)\right)\textrm{, }n=1,2,\ldots
\end{equation}
\textbf{vii)} The equation
$$
P_1(x)+P_2(y)+P_3(z)P_4(w)=n
$$ 
Have
$$
R_{1234}(n)=\sum_{\scriptsize \begin{array}{cc}
k+l=n\\
k\geq 1\textrm{\scriptsize, }l\geq 1 
\end{array}\normalsize}R_{12}(k)R_{34}(l)
$$
solutions, where
$$
R_{12}(n)=\sum_{\scriptsize \begin{array}{cc}
k+l=n\\
k\geq 1\textrm{\scriptsize, }l\geq 1 
\end{array}\normalsize}R_{1}(k)R_{2}(l)\textrm{ and }R_{34}(n)=\sum_{0<d|n}R_{3}(d)R_{4}\left(\frac{n}{d}\right)
$$
and for $j=1,2,3,4$, we have 
$$
R_{j}(n)=\sum_{\scriptsize \begin{array}{cc}
d\neq 0\textrm{\scriptsize, \normalsize}\textrm{\scriptsize abs\normalsize}(d)|n\\
P^{*}_{j}(d)=n
\end{array}\normalsize}1.
$$
The functions $P_j(x)$ are all polynomials of positive degree, involving also absolute values of $x$.\\
\\
\textbf{Example.}\\
Assume
$$
\theta_1(q)=\sum^{\infty}_{n=-\infty}q^{2n^2+n}\textrm{ and }\theta_2(q)=\sum^{\infty}_{n=-\infty}q^{n^2-3n}.
$$
Then
$$
A_1(n)=\left\{\begin{array}{cc}
1\textrm{, if }n=0\\
\sum_{\scriptsize\begin{array}{cc}
d\neq 0\textrm{\scriptsize, abs}(d)|n\\
2d^2+d=n
\end{array}}1\textrm{, if }n>0
\end{array}\right\}.
$$
$$
\sum^{\infty}_{n=-\infty}q^{2n^2+n}=\sum^{\infty}_{n=0}A_1(n)q^n
$$
and
$$
A_2(n)=\left\{\begin{array}{cc}
2\textrm{, if }n=-2\\
1\textrm{, if }n=0\\
\sum_{\scriptsize\begin{array}{cc}
d\neq 0\textrm{\scriptsize, abs}(d)|n\\
d^2-3d=n
\end{array}\normalsize}1\textrm{, if }n>0
\end{array}\right\}.
$$ 
Hence
$$
\sum^{\infty}_{n=-\infty}q^{n^2-3n}=2q^{-2}+\sum^{\infty}_{n=0}A_2(n)q^n.
$$
Then easily we get that equation
$$
2x^2+x+y^2-3y=n\textrm{, }x,y\in\textbf{Z}\textrm{, }n\in\textbf{N},
$$ 
have
$$
r(n)=2A_2(n+2)+\sum_{\scriptsize\begin{array}{cc}
k\geq 0,l\geq 0\\
k+l=n
\end{array}\normalsize}A_1(k)A_2(l).
$$
solutions.\\
\\
\textbf{Example.}\\
Assume the form
$$
|x|^3+|y|^3-x-y+z^2w^2=n\textrm{, }x,y,z,w,\in\textbf{Z}
$$
where $z\neq 0$ or $w\neq 0$. Then
$$
R_1(n)=\left\{\begin{array}{cc}
2\textrm{, if }n=0\\
\sum_{\scriptsize\begin{array}{cc}
d\neq 0\textrm{, \scriptsize abs}(d)|n\\
\textrm{\scriptsize abs}(d)^3-d=n
\end{array}
\normalsize}1
\textrm{, if }n>0
\end{array}\right\}.
$$
$$
A_1(n)=\sum^{n}_{l=0}R_1(l)R_1(n-l).
$$
$$
R_2(n)=\left\{\begin{array}{cc}
1\textrm{, if }n=0\\
\sum_{\scriptsize\begin{array}{cc}
d\neq 0\textrm{, \scriptsize abs}(d)|n\\
d^2=n
\end{array}
\normalsize}1
\textrm{, if }n>0
\end{array}\right\}.
$$
$$
A_2(n)=\sum_{d|n}R_2(d)R_2\left(\frac{n}{d}\right).
$$
Hence
$$
r(n)=\sum^{n}_{l=0}A_1(l)A_2(n-l).
$$

\section{The functional equation of $\sum^{\infty}_{n=1}\chi(n) n^{-\nu s}$}

For $\nu$ integer greater than or equal to 2, we define the function 
$$
\phi_{\nu}(z):=\sum^{\infty}_{n=-\infty}q^{|n|^{\nu}}=1+2\phi^{(0)}_{\nu}(z),
$$
where
$$ \phi^{(0)}_{\nu}(z)=\sum^{\infty}_{n=1}X_{\nu}(n)q^n\textrm{, }q=e(z):=e^{2\pi i z}\textrm{, }Im(z)>0
$$
and $X_{\nu}(n)$ is such that $X_{\nu}(n)=1$, if $n$ is a power of $\nu$, otherwise 0. It is easy to see someone that $X_{\nu}(n)$ is multiplicative arithmetic function on $\textbf{N}$ i.e. if $n,m\in\textbf{N}$ and the greatest common divisor $(n,m)$ of $n,m$ is $(n,m)=1$ i.e. coprime, then 
$$
X_{\nu}(nm)=X_{\nu}(n)X_{\nu}(m).
$$
More generaly it holds
$$
X_{\nu}(n m)X_{\nu}\left((n,m)\right)=X_{\nu}(n)X_{\nu}(m)\textrm{, }\forall n,m\in\textbf{N}.
$$
Then also the function $\phi_{\nu}(z)$ can be written as 
$$
\phi_{\nu}(z)=1+2\sum^{\infty}_{n=1}X_{\nu}(n)q^{n}\textrm{, }q=e(z)\textrm{, }Im(z)>0.
$$
The Dirichlet series attached with $\phi_{\nu}(z)$ is
$$
L\left(X_{\nu};s\right)=\sum^{\infty}_{n=1}\frac{X_{\nu}(n)}{n^s}=\zeta(\nu s).
$$
This function rises if we consider the Mellin transform of $
\phi^{(0)}_{\nu}(it)$. Hence if 
$$
\Lambda(\nu;s):=\int^{\infty}_{0}\phi^{(0)}_{\nu}(it)t^{s-1}dt,
$$  
then
$$
\Lambda(\nu;s)=(2\pi)^{-s} \Gamma(s)\sum^{\infty}_{n=1}\frac{X_{\nu}(n)}{n^s}=(2\pi)^{-s} \Gamma(s)\zeta(\nu s).
$$
If we also define
\begin{equation}
\Lambda^{*}(\nu,s):=(2\pi)^{s}\pi^{-s\nu/2}\Gamma\left(\frac{s\nu}{2}\right)\Gamma\left(\frac{1}{\nu}-s\right)\Lambda(\nu;s),
\end{equation}
then
\begin{equation}
\Lambda^{*}\left(\nu;\frac{1}{\nu}-s\right)=\Lambda^{*}(\nu;s).
\end{equation}
This can be seen easily using the Riemann's functional equation for the $\zeta(s)$ function:
$$
\pi^{-s/2}\Gamma\left(\frac{s}{2}\right)\zeta(s)=\pi^{-(1-s)/2}\Gamma\left(\frac{1-s}{2}\right)\zeta\left(1-s\right).
$$
We have 
$$
\Lambda^*\left(\nu;\frac{s}{\nu}\right)=\pi^{-s/2}\Gamma\left(\frac{s}{2}\right)\Gamma\left(\frac{1-s}{\nu}\right)\zeta\left(s\right)\Gamma\left(\frac{s}{\nu}\right)=
$$
\begin{equation}
=\xi(s)\Gamma\left(\frac{1-s}{\nu}\right)\Gamma\left(\frac{s}{\nu}\right)=\xi(1-s)\Gamma\left(\frac{1-s}{\nu}\right)\Gamma\left(\frac{s}{\nu}\right).
\end{equation}
Also
$$
\Lambda^{*}\left(\nu,\frac{1-s}{\nu}\right)=\pi^{-(1-s)/2}\Gamma\left(\frac{1-s}{2}\right)\Gamma\left(\frac{s}{\nu}\right)\Gamma\left(\frac{1-s}{\nu}\right)\zeta(1-s)=
$$
\begin{equation}
=\xi(1-s)\Gamma\left(\frac{s}{\nu}\right)\Gamma\left(\frac{1-s}{\nu}\right).
\end{equation}
From (343) and (344), we get (342).\\  
\\

Set now
$$
\phi_{\nu}(\chi;z)=\sum^{\infty}_{n=1}\chi(n)q^{n^{\nu}}\textrm{, }q=e(z)\textrm{, }Im(z)>0,
$$
where $\chi(n)$ is a primitive character to the modulus $k$ 
and
$$
L\left(\chi;s\right):=\sum^{\infty}_{n=1}\frac{\chi(n)}{n^s}.
$$
Then set
$$
\Lambda_{\nu}(\chi;s):=\int^{\infty}_{0}\phi_{\nu}\left(\chi;it\right)t^{s-1}dt=\sum^{\infty}_{n=1}\chi(n)\int^{\infty}_{0}e^{-2\pi n^{\nu}t}t^{s-1}dt=
$$
\begin{equation}
=(2\pi)^{-s}\Gamma(s)L\left(\chi;s\nu\right).
\end{equation}
The functional equation of $L\left(\chi;s\right)$ setting
\begin{equation}
L_1\left(\chi;s\right):=\left(\frac{\pi}{k}\right)^{-(s+a)/2}\Gamma\left(\frac{s+a}{2}\right)L(\chi;s),
\end{equation}
is
\begin{equation}
L_1\left(\overline{\chi};1-s\right)=\frac{i^ak^{1/2}}{\tau(\chi)}L_1(\chi;s).
\end{equation}
Hence\\
\\
\textbf{Theorem 65.}\\
If we set
\begin{equation}
\Lambda^{*}\left(\nu,\chi;s\right):=(2\pi)^{s}\sec\left(\frac{a\pi s\nu}{2}\right)\left(\frac{\pi}{k}\right)^{-s\nu/2}\Gamma\left(\frac{s\nu}{2}\right)\Gamma\left(\frac{1}{\nu}-s\right)\Lambda_{\nu}(\chi;s),
\end{equation}
then the functional equation of $\Lambda^{*}(\nu,\chi;s)$ is\\
\begin{equation}
\Lambda^{*}\left(\nu,\overline{\chi};\frac{1}{\nu}-s\right)=C_{\{a,k\}}\Lambda^{*}\left(\nu,\chi;s\right),
\end{equation}
where
\begin{equation}
C_{\{a,k\}}=\frac{i^ak^{1/2}}{\tau(\chi)},
\end{equation}
with $\tau(\chi)$ being the Gauss sum
\begin{equation}
\tau(\chi)=\sum^{k}_{n=1}\chi(n)\exp\left(2\pi i n/k\right)\textrm{, }|\tau(\chi)|=k^{1/2}.
\end{equation}  
The number $a$ is 
\begin{equation}
a=\left\{\begin{array}{cc}
0\textrm{, if }\chi(-1)=1,\\
1\textrm{, if }\chi(-1)=-1
\end{array}\right\}.
\end{equation}
\\

In view of [23] there exists the next\\
\\
\textbf{Theorem 66.}\\ 
For every function $X(x)$ there exists $P(x)$ such that $X(x)$ satisfies the DE
\begin{equation}
X'(x)\pm 2^{4/3}x^{-2/3}\left(1-x^2\right)^{-1/3}P\left(X(x)\right)=0.
\end{equation}
Then for the function $Y(z)$ such that $Y(z)=X\left(m^{*}\left(2\sqrt{N}z\right)\right)$, where $m^{*}(z)$ is the elliptic singular modulus defined as
\begin{equation}
m^{*}(z)=\left(\frac{\theta_2(e^{i\pi z})}{\theta_3(e^{i\pi z})}\right)^2,
\end{equation}
holds
\begin{equation}
Y'(z)=\pm 4\pi i\eta\left(z\right)^4P\left(Y(z)\right),
\end{equation}
where $\eta(z)$ is the Dedekind's eta function. Aslo if $a_n$ are defined as
\begin{equation}
\frac{1}{P(z)}=\sum^{\infty}_{n=0}a_nq^n,
\end{equation}
then
$$
e\left(a_0Y(z)\right)\prod^{\infty}_{n=1}\left(1-e\left(nY(z)\right)\right)^{-1/n\sum_{d|n}a_d\mu\left(n/d\right)}=
$$
\begin{equation}
=C\exp\left(8\pi^2\int^{z}_{i\infty}\eta\left(w\right)^4dw\right).
\end{equation}
Moreover for the function $Y(z)$ holds
\begin{equation}
\int^{Y(z)}_{Y(i\infty)}\frac{1}{P(t)}dt=\mp 4\pi i \int^{z}_{i\infty}\eta(t)^4dt 
\end{equation}
and
\begin{equation}
Y\left(-\frac{1}{z}\right)=F_{reg}^{(-1)}\left(1-F_{reg}(Y(z))\right),
\end{equation}
where
\begin{equation}
F_{reg}(x)=\frac{1}{\int^{c'}_{c}\frac{dt}{P(t)}}\int^{x}_{c}\frac{dt}{P(t)}.
\end{equation}
\\

Assume now the function 
\begin{equation}
\theta^{\{\nu\}}_{\chi}(z)=\sum^{\infty}_{n=1}\chi(n)q^{n^{\nu}}\textrm{, }q=e(z)\textrm{, }Im(z)>0.
\end{equation}
In this case for the $\theta_{\chi}^{\{\nu\}}(z)$ holds:
\begin{equation}
\exp\left(2\pi i \int_{i\infty}^{z}\theta_{\chi}^{\{\nu\}}(w)dw\right)=e^{(1-\chi(1))q}\prod^{\infty}_{n=1}\left(1-q^n\right)^{-C_{\nu}(\chi;n)/n},
\end{equation}
where
\begin{equation}
C_{\nu}(\chi,n)=\mu(n)+\chi(n_{\nu}(n))\mu(n^{*}_{\nu}(n)).
\end{equation}
Also we have that exists functions $X(x)$,$P(x)$ such that $Y(z)=X\left(m^{*}\left(2z\right)\right)$, with $c=Y(i\infty)$ and
$$
\exp\left(2\pi i \int_{Y(i\infty)}^{Y(z)}\theta_{\chi}^{\{\nu\}}(w)dw\right)=
$$
$$
=C_1\cdot e\left(a_0Y(z)\right)\prod^{\infty}_{n=1}\left(1-e\left(nY(z)\right)\right)^{-1/n\sum_{d|n}a_d\mu\left(n/d\right)}=
$$
\begin{equation}
=C_2\cdot\exp\left(8\pi^2\int^{z}_{i\infty}\eta\left(w\right)^4dw\right).
\end{equation}
In this case $a_0=0$ and $a_n=\sum_{d|n}C_{\nu}(\chi,d)$ if $n\geq 1$. Hence there exist integer $k$ such that
$$
2\pi i\int^{z}_{c}\theta^{\{\nu\}}_{\chi}(w)dw=C_0+8\pi^2\int^{Y^{(-1)}(z)}_{i\infty}\eta(w)^4dw+2\pi i k\textrm{, }k\in\textbf{Z}.
$$
Hence if we define $Y^{\{\nu\}}_{\chi}(z)=Y(z)$, then
\begin{equation}
\int^{Y^{\{\nu\}}_{\chi}(z)}_{c}\theta^{\{\nu\}}_{\chi}(w)dw+4\pi i\int^{z}_{i\infty}\eta(w)^4dw=C_0.
\end{equation}
The function $Y^{\{\nu\}}_{\chi}(z)$ attached to $\theta^{\{\nu\}}_{\chi}(z)$ in this way have the property
\begin{equation}
Y^{\{\nu\}}_{\chi}\left(-\frac{1}{z}\right)=F_{reg}^{(-1)}\left(1-F_{reg}\left(Y^{\{\nu\}}_{\chi}(z)\right)\right),
\end{equation}
where
\begin{equation}
F_{reg}(z)=\frac{1}{\int^{c'}_{c}\theta^{\{\nu\}}_{\chi}(w)dw}\int^{z}_{c}\theta^{\{\nu\}}_{\chi}(w)dw
\end{equation}
and $c=Y(i\infty)$, $c'=Y(i0)$. Hence from Theorem 50 (i) pp 53, we get
$$
2\pi i \int^{z}_{c}\theta^{\{\nu\}}_{\chi}(t)dt=\theta^{\{\nu\}}_{\chi}\left(Z\left(Y^{{\{\nu\}}(-1)}_{\chi}(z)\right)\right)+c_1\Leftrightarrow
$$
$$
\exp\left(2\pi i \int^{Y^{\{\nu\}}_{\chi}(z)}_{c}\theta^{\{\nu\}}_{\chi}(t)dt\right)=\exp\left(C_1+\theta^{\{\nu\}}_{\chi}(Z(z))\right)=
$$
$$
=\exp\left(C_0+8\pi^2\int^{z}_{i\infty}\eta(w)^4dw\right).
$$
Hence
$$
\theta^{\{\nu\}}_{\chi}\left(Z(z)\right)=C_1+8\pi^2\int^{z}_{i\infty}\eta(w)^4dw.
$$
Also
$$
\frac{1}{P(z)}=\sum^{\infty}_{n=1}\left(\sum_{d|n}C_{\nu}(\chi,d)\right)q^n\Rightarrow
$$
$$
F_{reg}(z)=C+\frac{1}{2\pi i}\sum^{\infty}_{n=1}\left(\frac{1}{n}\sum_{d|n}C_{\nu}(\chi,d)\right)q^{n}
$$
and
$$
\frac{1}{P(z)}=\theta^{\{\nu\}}_{\chi}(z)+C_1.
$$

\section{Borcherds products}

\textbf{Theorem 67.}(see [27])\\ 
Let $X(n)=\left(n|G\right)$ be the Jacobi symbol and $G=2^mg_1^{m_1}g_2^{m_2}\ldots g_s^{m_s}$, with  $m,m_1,\dots,m_s$ non negative integers, $m\neq1$ and $g_1<g_2<\ldots<g_s$ are primes congruent to  $1(\textrm{mod}4)$. Then exist rational $A$ such that
\begin{equation}
R(q):=q^A\prod^{\infty}_{n=1}(1-q^n)^{\left(n|G\right)}=q^A\prod^{\left[\frac{G-1}{2}\right]}_{j=1}\theta\left(\frac{G}{2},\frac{G}{2}-j;q\right)^{\left(j|G\right)},
\end{equation}
where 
\begin{equation}
A=\sum^{\left[\frac{G-1}{2}\right]}_{j=1}\left(-\frac{j}{2}+\frac{j^2}{2G}+\frac{G}{12}\right)\left(j|G\right),
\end{equation}
\begin{equation}
\theta(a,b;q)=\theta_{a,b}(q):=\sum^{\infty}_{n=-\infty}(-1)^nq^{an^2+bn}.
\end{equation}
Also $R\left(e^{-\pi\sqrt{r}}\right)$ is algebraic, when $r$ positive rational.\\
\\

From Theorem 28 we get
$$
\exp\left(\sum^{\infty}_{n=1}\chi_1(n)q^n\right)=\prod^{\infty}_{n=1}(1-q^n)^{\left(n|G\right)},
$$
where
$$
\chi_1(n)=-\frac{1}{n}\sum_{d|n}(d|G)d.
$$
Hence
$$
R(q):=q^A\exp\left(-\sum^{\infty}_{n=1}\frac{q^n}{n}\sum_{d|n}(d|G)d\right)=q^A\prod^{\infty}_{n=1}(1-q^n)^{(n|G)}=
$$
$$
=q^A\prod^{\left[\frac{G-1}{2}\right]}_{j=1}\theta\left(\frac{G}{2},\frac{G}{2}-j;q\right)^{\left(j|G\right)}.
$$
and we have the next\\
\\
\textbf{Theorem 68.}\\
Let $q=e(z)$, $Im(z)>0$. Then\\ 
i) The function $\frac{d}{dz}\log R(q)$ is a wight 2 modular form in $M_2(\Gamma(G),\chi)$, where $\chi(n)=(n|G)$.\\ Hence also\\
ii) The function
\begin{equation}
F(z)=2\pi i A+\sum^{\left[\frac{G-1}{2}\right]}_{j=1}(j|G)\frac{\theta^{*}\left(\frac{G}{2},\frac{G}{2}-j;q\right)}{\theta\left(\frac{G}{2},\frac{G}{2}-j;q\right)},
\end{equation} 
where $\theta^*(a,b;q)=\frac{d}{dz}\theta(a,b;q)$, is also a weight 2 modular form in $M_2(\Gamma(G),\chi)$, with $\chi(n)=(n|G)$.\\
\\

Also (see [26]) when $\chi(n)$ is a Dirichlet character modulo $k$ and $\nu$ a positive integer such that $\chi(-1)=(-1)^{\nu}$, then (assuming also that $q=e(z)$) we get using Theorem 31:
$$
f(z):=q^{-c_{\nu}(\chi)}\prod^{\infty}_{n=1}(1-q^n)^{\chi(n)n^{\nu-2}}=
$$
\begin{equation}
=\exp\left(-2 \pi i z c_{\nu}(\chi)-\sum^{\infty}_{n=1}\frac{q^n}{n}\sum_{d|n}\chi(d)d^{{\nu}-1}\right).
\end{equation}
Hence we have the next:\\
\\
\textbf{Theorem 69.}\\
When $\nu$ is integer greater than 1 and $\chi(n)$ a Dirichlet character $\textrm{mod}k$ with $\chi(-1)=(-1)^{\nu}$, the function $\frac{d}{dz}\log f(z)$, (where $f$ is that of (372)), is a modular form of weight $\nu$ in $M_{\nu}(\Gamma(k),\chi)$, where $c_{\nu}(\chi)=\frac{1}{2}L(1-\nu,\chi)$, $L(s,\chi)=\sum^{\infty}_{n=1}\chi(n)n^{-s}$.\\
\\
\textbf{Corollary 70.}\\
Let $v(z)=R(q)$, where $R(q)$ denotes the Rogers-Ramanujan continued fraction:
\begin{equation}
R(q)=\frac{q^{1/5}}{1+}\frac{q^1}{1+}\frac{q^2}{1+}\frac{q^3}{1+}\ldots.
\end{equation}
Then $\frac{d}{dz}\log(v(z))$ is a weght 2 modular form in $M_2(\Gamma(5),\chi)$, where $\chi(n)=(n|5)$.\\
\\
\textbf{Theorem 71.}\\
The Ramanujan quantity (see [25]) is defined for $a,b,p$ positive integers, $a<b<a+b<p$ as
\begin{equation}
R(a,b;p;q):=q^{-(a-b)/2+(a^2-b^2)/(2p)}\frac{(q^a;q^p)_{\infty}(q^{p-a};q^p)_{\infty}}{(q^b;q^p)_{\infty}(q^{p-b};q^p)_{\infty}}.
\end{equation}
Then the function $u(z)=\frac{d}{dz}\log R(a,b;p;e(z))$, $Im(z)>0$, is a modular form in $M_2(\Gamma(p),X)$, where $X$ is defined as
\begin{equation}
X(a,b;p;n)=\left\{\begin{array}{cc}
1\textrm{, if }n\equiv a(\textrm{mod}p)\\
-1\textrm{, if }n\equiv p-a(\textrm{mod}p)\\
1\textrm{, if }n\equiv b(\textrm{mod}p)\\
-1\textrm{, if }n\equiv p-b(\textrm{mod}p)\\
0\textrm{, else }
\end{array}\right\}.
\end{equation}  
Also we can write
\begin{equation}
R(a,b;p;q)=q^A\prod^{\infty}_{n=1}(1-q^n)^{X(a,b;p;n)},
\end{equation}
where $A=-(a-b)/2+(a^2-b^2)/(2p)$.\\
\\
\textbf{Theorem 72.}\\
Assume the function
\begin{equation}
R_{\nu}(\chi,q)=q^A\prod^{\infty}_{n=1}(1-q^n)^{\chi(n)n^{\nu-2}},
\end{equation}
where $\chi(n)$ is real Dirichlet character $\textrm{mod}k$ with $\chi(-1)=(-1)^{\nu}$, $\nu$ is integer greater than 1 and $A=\frac{1}{2}L(\chi,1-\nu)$. Then the function 
\begin{equation}
v_{\nu}(z)=\frac{d}{dz}\log R_{\nu}\left(\chi,e(z)\right)\textrm{, }Im(z)>0
\end{equation}
is a modular form in $M_{\nu}(\Gamma(k),\chi)$.\\
\\
\textbf{Corollary 73.}\\
Suppose $\chi$ is any real Dirichlet character $\textrm{mod}k$ with $\chi(-1)=1$ and
\begin{equation}
A=\sum^{\left[\frac{k-1}{2}\right]}_{j=1}\left(-\frac{j}{2}+\frac{j^2}{2k}+\frac{k}{12}\right)\chi(j).
\end{equation}
Then the function
\begin{equation}
R(\chi,q)=q^A\prod^{\infty}_{n=1}\left(1-q^n\right)^{\chi(n)}\textrm{, }q=e(z)
\end{equation}
i) is algebraic in all points $q=e(z)$, such that $z=r_1+i\sqrt{r_2}$, $r_1\in\textbf{Q}$ and $r_2\in\textbf{Q}^{*}_{+}$.\\ 
ii) The function
\begin{equation}
\frac{d}{dz}\log R(\chi,e(z))\in M_{2}(\Gamma(k),\chi).
\end{equation}
\\
\textbf{Remark.} Actually all Dirichlet characters with $\chi(-1)=1$ are complex Moebius periodic (see [27]).\\
\\
\textbf{Theorem 74.}\\
Assume the Dirichlet character $\textrm{mod}k$, $\chi$ and an integer $\nu>1$, with $\chi(-1)=(-1)^{\nu}$. Then the function
\begin{equation}
v_0(z)=v_0(\chi,z)=R(\chi,q)=q^A\prod^{\infty}_{n=1}(1-q^n)^{\chi(n)n^{\nu-2}}\textrm{, }q=e(z),
\end{equation}
$Im(z)>0$, have logarithmic derivative
\begin{equation}
v(z)=\frac{d}{dz}\log v_0(z)=-2\pi i A+\pi\sum^{\infty}_{n=1}\frac{\chi(n)n^{\nu-1}e^{in\pi z}}{\sin(n\pi z)},
\end{equation}
where $A=\frac{1}{2}L(\chi,1-\nu)$ and thus this derivative is a modular form in $M_{\nu}(\Gamma(k),\chi)$.

\section{Quadratic residues and representations in certain forms}

Assume $n,a$ are integers with $n>1$ and the congruence
\begin{equation}
x^2\equiv a(\textrm{mod}n)
\end{equation}
and $p_1,p_2,\dots,p_r$ are all the odd prime divisors of $n$. Then we know that (357) have solution iff all 
\begin{equation}
x^2\equiv a(\textrm{mod}p_i)\textrm{, }i=1,2,\ldots,r
\end{equation}
have solution, $a\equiv1(\textrm{mod}4)$ with $4|n$ and $a\equiv 1(\textrm{mod}8)$ with $8|n$. The number of solutions is $2^{r+u}$, where $u=0,1,2$ resp. if $4$ not devides $n$, $4|n$ exactly, $8|n$.\\
Also it is known that the equation
\begin{equation}
x^2=a+py,
\end{equation}
where $p$ is prime have solution iff the Jacobi symbol $(a|p)=0,1$.\\
Assume now the quadratic residue equation 
\begin{equation}
x^2=a+n z\textrm{, }(a,n)=1
\end{equation}
and let $\textrm{res}(a,n)$ be the number of its solutions. Then the equation
\begin{equation}  
x^2=-ky^2+nz\textrm{, }(y,n)=1\textrm{, }(k,n)=1,
\end{equation}
have $r_{k,n}(y)=\textrm{res}(-ky^2,n)$ solutions. The arithmetical function $r_{k,n}(y)$, when $(k,n)=1$, is a Dirichlet character $\textrm{mod}(\textrm{rad}(n))$ multiplied by a suitable constant $c$ and $\textrm{rad}(n)$ is the radical function of $n$ which is given by
\begin{equation}
\textrm{rad}(n)=\prod_{p|n}p.
\end{equation}
Actualy $r_{k,n}(y)=c\cdot\chi_{1,\textrm{\scriptsize rad\normalsize}(n)}(y)$ is the principal Dirichlet character $\textrm{mod}( \textrm{rad}(n))$.\\
Hence given the binary quadratic congruence
\begin{equation}
x^2+ky^2\equiv0(\textrm{mod}n)\textrm{, with }(k,n)=1\textrm{ and }(y,n)=1,
\end{equation}
we have for $y=1,2,\ldots,\textrm{rad}(n)$, respectively $\textrm{res}(-ky^2,n)$ solutions. The set of these solutions is $\textrm{rad}(n)$ periodic. Hence it is suficient to consider (390) only for 
\begin{equation}
A=\{y\in\textbf{N}:1\leq y\leq \textrm{rad}(n)\}.
\end{equation}
The totality of the solutions of (390) in $A$ is
\begin{equation}
c\sum^{\textrm{\scriptsize rad\normalsize}(n)}_{m=1}\chi_{\{1,\textrm{\scriptsize rad\normalsize}(n)\}}(m)=c\cdot \phi(\textrm{rad}(n)),
\end{equation}
where $\phi$ is Euler's$-$totient arithmetical function. The constant $c$ is 0, or an integer of the form $2^m$.\\
\\
\textbf{Example.}\\ 
Assume a positive integer $n$ and a prime number $p>n$. Then exists constant $c$ such that the congruence
\begin{equation}
x^2+py^2\equiv0(\textrm{mod}n!)\textrm{, }(y,n!)=1,
\end{equation}
have $c\cdot\phi\left(\prod_{p\leq n}p\right)=c\prod_{p\leq n}(p-1)$ solutions. For $n=10$ the first two primes such that $c\neq 0$ are $p=311$ and $p=479$.\\
\\
\textbf{Example.}\\
Assume the integers $a,k,\nu$, with $a>1$, $\nu>2$ and $k$ such that $(k,a^{\nu})=1$. Then the congruence
\begin{equation}
x^2+ky^2\equiv 0(\textrm{mod}a^{\nu}),
\end{equation}
have $c\cdot \phi\left(\textrm{rad}\left(a^{\nu}\right)\right)=c\cdot \phi(\textrm{rad}(a))$ solutions.\\
\\
\textbf{Theorem 75.}\\
Let $p,q$ be prime numbers, with $q$ odd. Then the congruence
\begin{equation}
x^2+py^2\equiv0(\textrm{mod}q)\textrm{, }(y,q)=1,
\end{equation}
have $2c(p,q)\cdot(q-1)$ solutions (congruences) in the period $1\leq y\leq q-1$. The arithmetical function $c$ is: $c(p,q)=1$ if $(-p|q)=0,1$ and $c(p,q)=0$, if $(-p|q)=-1$.\\
\\
\textbf{Lemma 76.}\\
Let $p,q$ be prime numbers, with $q>2$. If $(-p|q)=-1$, then equation
\begin{equation}
x^2+p y^2=kq\textrm{, }(y,q)=1\textrm{, }x,y\in\textbf{Z}\textrm{, }k\in\textbf{N},
\end{equation}
is imposible.\\
\\
\textbf{Note.}\\
The Lemma can be extended in the case where $p$ is prime number and $n$ any positive integer with $(n,2)=1$. Then if $(-p|n)=-1$, then the equation
\begin{equation}
x^2+p y^2=n\textrm{, }(y,n)=1\textrm{, }x,y\in\textbf{Z},
\end{equation}
is imposible.\\ 
\\
More general:\\
\\
\textbf{Lemma 77.}\\
Let $a,b,n$ be integers with $a,b,n\geq1$ and $(n,2)=1$. If $(-ab|n)=-1$, then the equation
\begin{equation}
ax^2+by^2=n\textrm{, }x,y\in\textbf{Z},
\end{equation}
is imposible.\\
\\
\textbf{Notes.}\\ 
Let $a,b,n,t$ be positive integers. Then if $(n,2)=(t,2)=1$ and $(-ab|t)=\pm1$ the equation
$$
ax^2+by^2=tn,
$$
is imposible when $(-ab|n)=\mp1$.\\
\\
\textbf{Theorem 78.}\\
Let $t,n$ be integers with $t$ positive prime of the form $t\equiv 3(\textrm{mod}4)$. Then the moduli set $\overline{t}=\{1,2,\ldots,t\}$ is spliten into three disjoint classes $S_1$, $S_{-1}$, $S_0$, with $\overline{t}=S_1\cup S_{-1}\cup S_0$, where $S_i=\{n\in\overline{t}:(-t|n)=i\}$.\\ 
\textbf{1)} If $n\in S_{-1}$,  then the equation  
$$
x^2+ty^2=n,\eqno{(eq)}
$$
have no solution.\\ 
\textbf{2)} If $n^2\equiv k_1(\textrm{mod}t)$, $k_1\in S_{1}$ then 
$$
x^2+ty^2=n^2
$$
have always solution.\\
\textbf{3)} If $n$ have representation $n=n_1^2n_2$, with $n_1\geq1$ and $\mu(n_2)\neq 0$, (which is unique) and $(-t|n_2)=-1$, then we have no solution.

\section{A summation formula}

In view of [17] we have shown the next\\
\\
\textbf{Theorem 79.}\\
Let $f(z)$ be a function of the form
\begin{equation}
f(z)=\sum^{\infty}_{n=1}\chi(n)z^{n^{\nu}}\textrm{, }|z|<1,
\end{equation}
where $\nu$ is positive integer, $\nu>1$. Then if $\sum^{\infty}_{k=1}|a_k|<\infty$, we have the next summation formula
\begin{equation}
\sum^{n}_{k=1}f(a_k)=-\sum^{\infty}_{m=1}\log\left(\prod^{n}_{k=1}\left(1-a_k^m\right)\right)\chi(n_{\nu}(m))\epsilon(n^{*}_{\nu}(m)).
\end{equation}
\\
\textbf{Corollary 80.}\\
Let $f(z)$ be a function of the form
\begin{equation}
f(z)=\sum^{\infty}_{n=1}\chi(n)z^{n^{\nu}}\textrm{, }|z|<1,
\end{equation}
where $\nu$ integer, $\nu\geq 2$ and $\chi(n)$ any arithmetical function. Then if $s>1$, it holds
\begin{equation}
\sum_{p-prime}f\left(\frac{1}{p^s}\right)=\sum^{\infty}_{m=1}\log\left(\zeta(sm)\right)\chi\left(n_{\nu}(m)\right) \epsilon(n^{*}_{\nu}(m)),
\end{equation}
where $\zeta(s)$ is the Riemann's zeta function.\\
\\
\textbf{Corollary 81.}\\
If we define
\begin{equation}
(a;q)_n=\prod^{n}_{k=1}(1-aq^k),
\end{equation}
then
\begin{equation}
\exp\left(\sum^{n}_{k=1}f\left(aq^k\right)\right)=\prod^{\infty}_{m=1}\left(a^m;q^m\right)_{n}^{-\chi(n_{\nu}(m)) \epsilon(n^{*}_{\nu}(m))},
\end{equation}
where
\begin{equation}
f(z)=\sum^{\infty}_{n=1}\chi(n)z^{n^{\nu}}\textrm{, }|z|<1.
\end{equation}
\\
\textbf{Applications.}\\
\textbf{1)} If $\chi(1)=1$ and $0$ else, then
\begin{equation}
\sum_{p-prime}\frac{1}{p^s}=\sum^{\infty}_{m=1}\log\left(\zeta(sm)\right)\frac{\mu\left(m\right)}{m}.
\end{equation}
\textbf{2)} If $\chi(1)=1$ and $\chi(n)=0$, $\forall n\in\{2,3,\ldots\}$ and $\nu$ integer, $\nu>1$, $|q|<1$, $\left(aq^j\right)^m\neq 1$, then $f(z)=z$ and
\begin{equation}
\exp\left(\sum^{n}_{k=1}aq^k\right)=\prod^{\infty}_{m=1}\left(a^m;q^m\right)_{n}^{-\frac{\mu\left(m\right)}{m}}.
\end{equation}
\textbf{3)} 
\begin{equation}
\exp\left(\sum^{\infty}_{n=1}f(q^n)\right)=\prod^{\infty}_{n=1}\left(1-q^n\right)^{-\sum_{d|n}\chi(n_{\nu}(d))\epsilon(n^{*}_{\nu}(d))},
\end{equation}
where
$$
f(q)=\sum^{\infty}_{n=1}\chi(n)q^{n^{\nu}}\textrm{, }|q|<1.
$$
\\
\textbf{Theorem 82.}\\
Assume that 
\begin{equation}
f(z)=\sum^{\infty}_{n=1}\chi(n)z^{n^{\nu}}\textrm{, }|z|<1,
\end{equation}
and $\rho_1,\rho_2,\ldots$ are roots of the holomorphic function $g(z)$ in $\textbf{C}$. Assume also that $g(0)=1$. If $\sum^{\infty}_{k=1}\frac{1}{\left|\rho_k\right|}<\infty$, then
$$
\exp\left(\sum^{\infty}_{k=1}f\left(\frac{z}{\rho_k}\right)\right)=\prod^{\infty}_{k=1}\prod^{k-1}_{j=0}g\left(z\zeta_{j,k}\right)^{-\chi\left(n^{(1)}_{\nu}(k)\right) \epsilon\left(n^{(2)}_{\nu}(k)\right)},
$$
where $\zeta_{j,k}=e^{2\pi i j/k}$.\\
\textbf{Remark.}\\
If $f\equiv g$, then we get the next formula
$$
\exp\left(\sum^{\infty}_{k=1}f\left(\frac{z}{\rho_k}\right)\right)=\prod^{\infty}_{k=1}\prod^{k-1}_{j=0}f\left(z\zeta_{j,k}\right)^{-\chi\left(n^{(1)}_{\nu}(k)\right) \epsilon\left(n^{(2)}_{\nu}(k)\right)},
$$
where 
$$
f(z)=\sum^{\infty}_{n=1}\chi(n)z^{n^{\nu}}\textrm{, }|z|<1
$$
and $\rho_1,\rho_2,\ldots$, are the roots of $f(z)$ in $\textbf{C}$.\\
\\
\textbf{Proof.}\\
Assume the holomorphic complex function $g(z)$, with roots $\rho_1,\rho_2,\ldots$ and $g(0)=1$. Since $\sum^{\infty}_{k=1}\frac{1}{|\rho_k|}<\infty$, according to Weierstrass$-$Handamard factorization theorem we can write
$$
g(z)=\prod^{\infty}_{n=1}\left(1-\frac{z}{\rho_n}\right).
$$
Hence according to Theorem 79 we get
$$
\exp\left(\sum^{\infty}_{k=1}f\left(\frac{z}{\rho_k}\right)\right)=\prod^{\infty}_{m=1}\prod^{\infty}_{k=1}\left(1-\left(\frac{z}{\rho_k}\right)^m\right)^{-\chi(n_{\nu}(m))\epsilon(n^{*}_{\nu}(m))}=
$$ 
$$
=\prod^{\infty}_{m=1}\prod^{m-1}_{j=0}g\left(z\zeta_{j,m}\right)^{-\chi(n_{\nu}(m))\epsilon(n^{*}_{\nu}(m))}.
$$

\[
\]

\centerline{\textbf{References}}

\[
\]

[1]: M. Abramowitz, I.A. Stegun. 'Handbook of Mathematical Functions'. Dover Publications, New York, (1972).\\

[2]: T. Apostol. 'Introduction to Analytic Number Theory'. Springer Verlag, New York, Berlin, Heidelberg, Tokyo, (1974).\\

[3]: J.V. Armitage, W.F. Eberlein. 'Elliptic Functions'. Cambridge University Press, (2006).\\

[4]: J.M. Borwein, P.B. Borwein. 'Pi and the AGM'. John Wiley and Sons, Inc. New York, Chichester, Brisbane, Toronto, Singapore, (1987).\\

[5]: L.E. Dickson. 'History of the Theory of Numbers, Vol2: Diophantine Analysis'. Dover. New York, (2005).\\

[6]: G.H. Hardy. 'Ramanujan Twelve Lectures on Subjects Suggested by his Life and Work, 3rd ed.' Chelsea. New York, (1999).\\

[7]: M.D. Hirschhorn. 'Three classical results on representations of a number'. Seminaire Lotharingien de Combinatoire, (42), (1999).\\  

[8]: C.G.J. Jacobi. 'Fundamenta Nova Functionum Ellipticarum'. Werke I, 49-239, (1829).\\ 

[9]: William J. LeVeque. 'Fundamentals of Number Theory'. Dover Publications. New York, (1996).\\

[10]: E.T. Whittaker, G.N. Watson. 'A course on Modern Analysis'. Cambridge U.P, (1927).\\

[11]: Ken Ono. 'Representations of Integers as Sums of Squares'. Journal of Number Theory. (95), 253-258. (2002).\\

[12]: Ila Varma. 'Sums of Squares, Modular Forms, and Hecke Characters'. Master thesis. Mathematisch Institut, Universiteit Leiden, June 18 (2010).\\ 

[13]: Konrad Knopp. 'Theory and Applications of Infinite Series'. Dover Publications, Inc. New York, (1990).\\

[14]: Bruce C. Berndt. 'Ramanujan`s Notebooks Part III'. Springer Verlag, New York, (1991).\\

[15]: G.E. Andrews, 'Number Theory'. Dover Publications, New York, (1994).\\ 

[16]: G.H. Hardy. 'On the expression of a number as the sum of two squares'. Quart. J. Math. (46), (1915).\\

[17]: N. Bagis. 'Some New Results on Sums of Primes'. Mathematical Notes, Vol. 90, No. 1, pp 10-19, (2011).\\

[18]: N. Bagis. 'Some Results on Infinite Series and Divisor Sums'.\\arXiv:0912.48152v2 [math.GM], (2014).\\

[19]: H. Iwaniec, E. Kowalski. 'Analytic Number Theory'. American Mathematical Society. Colloquium Publications. Vol 53. Providence, Rhode Island, (2004).\\

[20]: Carlos J. Moreno, Samuel S. Wagstaff,JR. 'Sums of Squares of Integers'. Chapman and Hall/CRC. Taylor and Francis Group. Boca Raton, London, New York, (2006).\\

[21]: N.G. De Bruijn. 'Asymptotic Methods in Analysis'. Dover Publications, Inc. New York, (1981).\\

[22]: B.C. Berndt. 'Ramanujan`s Notebooks Part II'. Springer-Verlag, New York. 1989.\\

[23]: N. Bagis. 'On the Complete Evaluation of Jacobi Theta Functions'. arXiv:1503.01141v4 [math.GM], (2021).\\

[24]: N. Bagis. 'On Generalized Integrals and Ramanujan-Jacobi Special Functions'. arXiv:1309.7247v3 [math.GM], (2015).
\\

[25]: N. Bagis. 'Generalizations of Ramanujan Continued Fractions'. arXiv:1107.2393v2 [math.GM], (2012).\\

[26]: D. Zagier. 'Elliptic Modular Forms and Their Applications'. page stored at the web.\\

[27]: N.D. Bagis, M.L. Glasser. 'Conjectures on the evaluation of certain functions with algebraic properties'. Journal of Number Theory. 155 (2015), 63-84.\\

\end{document}